\documentclass[leqno,draft]{amsart}
\usepackage{amssymb,enumerate, latexsym, mathrsfs, mdwlist}

\usepackage{appendix}

\swapnumbers
\theoremstyle{definition}
\newtheorem{nul}{}[section]
\newtheorem{dfn}[nul]{Definition}

\newtheorem{ntn}[nul]{Notation}
\newtheorem{exm}[nul]{Example}

\newtheorem{wrn}[nul]{Warning}
\newtheorem*{dfn*}{Definition}
\newtheorem*{axm*}{Axiom}
\newtheorem*{ntn*}{Notation}
\newtheorem*{exm*}{Example}
\newtheorem*{exr*}{Exercise}
\newtheorem*{int*}{Intuition}
\newtheorem*{qst*}{Question}

\theoremstyle{plain}

\newtheorem{thm}[nul]{Theorem}
\newtheorem{prp}[nul]{Proposition}
\newtheorem{lem}[nul]{Lemma}

\newtheorem{cor}{Corollary}[nul]
\newtheorem*{thm*}{Theorem}
\newtheorem*{mainthm*}{Universal Additivity Theorem}
\newtheorem*{structurethm*}{Structure Theorem}
\newtheorem*{prp*}{Proposition}
\newtheorem*{cor*}{Corollary}
\newtheorem*{lem*}{Lemma}
\newtheorem*{cnj*}{Conjecture}

\numberwithin{equation}{nul}

\DeclareMathOperator{\Aut}{Aut}

\DeclareMathOperator{\colim}{colim}

\DeclareMathOperator{\et}{\acute{e}t}

\DeclareMathOperator{\Exc}{Exc}

\DeclareMathOperator{\Fun}{Fun}
\DeclareMathOperator{\hocolim}{hocolim}
\DeclareMathOperator{\id}{id}
\DeclareMathOperator{\Ind}{Ind}

\DeclareMathOperator{\Map}{Map}
\DeclareMathOperator{\Mor}{Mor}

\DeclareMathOperator{\Spec}{Spec}

\renewcommand{\AA}{\mathbf{A}}

\newcommand{\CC}{\mathbf{C}}
\newcommand{\DD}{\mathbf{D}}

\newcommand{\FF}{\mathbf{F}}

\newcommand{\KK}{\mathbf{K}}
\newcommand{\LL}{\mathbf{L}}

\newcommand{\RR}{\mathbf{R}}
\renewcommand{\SS}{\mathbf{S}}

\newcommand{\VV}{\mathbf{V}}

\newcommand{\XX}{\mathbf{X}}

\newcommand{\Cat}{\mathbf{Cat}}

\newcommand{\Kan}{\mathbf{Kan}}

\newcommand{\Mack}{\mathbf{Mack}}
\newcommand{\Mod}{\mathbf{Mod}}

\newcommand{\Perf}{\mathbf{Perf}}

\renewcommand{\Pr}{\mathbf{Pr}}
\newcommand{\QCoh}{\mathbf{QCoh}}

\newcommand{\Set}{\mathbf{Set}}
\newcommand{\Sp}{\mathbf{Sp}}
\newcommand{\Trip}{\mathbf{Trip}}

\newcommand{\Wald}{\mathbf{Wald}}

\newcommand{\add}{\mathit{add}}

\newcommand{\coh}{\textit{coh}}

\newcommand{\disj}{\mathit{disj}}
\newcommand{\eff}{\textit{eff}}

\newcommand{\op}{\mathit{op}}

\newcommand{\coloneq}{\mathrel{\mathop:}=}

\makeatletter 
\def\revddots{\mathinner{\mkern1mu\raise\p@ 
\vbox{\kern7\p@\hbox{.}}\mkern2mu 
\raise4\p@\hbox{.}\mkern2mu\raise7\p@\hbox{.}\mkern1mu}} 
\makeatother

\newcommand{\upsidedown}[1]{\textrm{\raisebox{\depth}{\rotatebox{180}{$#1$}}}}

\usepackage{tikz}
\usetikzlibrary{matrix,arrows,decorations}

\pgfdeclaredecoration{single line}{initial}{\state{initial}[width=\pgfdecoratedpathlength-1sp]{\pgfpathmoveto{\pgfpointorigin}}\state{final}{\pgfpathlineto{\pgfpointorigin}}} 

\newcommand{\fromto}[2]{{#1}\ \tikz[baseline]\draw[>=stealth,->](0,0.5ex)--(0.5,0.5ex);\ {#2}}
\newcommand{\into}[2]{{#1}\ \tikz[baseline]\draw[>=stealth,right hook->](0,0.5ex)--(0.5,0.5ex);\ {#2}}

\newcommand{\cofto}[2]{{#1}\ \tikz[baseline]\draw[>=stealth,>->](0,0.5ex)--(0.5,0.5ex);\ {#2}}
\newcommand{\fibto}[2]{{#1}\ \tikz[baseline]\draw[>=stealth,->>](0,0.5ex)--(0.5,0.5ex);\ {#2}}
\newcommand{\equivto}[2]{{#1}\ \tikz[baseline]\draw[>=stealth,->,font=\scriptsize,inner sep=0.5pt](0,0.5ex)--node[above]{$\sim$}(0.5,0.5ex);\ {#2}}
\newcommand{\goesto}[2]{{#1}\ \tikz[baseline]\draw[|->](0,0.5ex)--(0.5,0.5ex);\ {#2}}
\newcommand{\adjunct}[4]{{#1}\colon{#2}\ \begin{tikzpicture}[baseline] \draw[>=stealth,->] (0,1ex) -- (0.75,1ex); \draw[>=stealth,->] (0.75,0.25ex) -- (0,0.25ex); \end{tikzpicture}\ {#3}\colon{#4}}

\renewcommand{\to}{\ \tikz[baseline]\draw[>=stealth,->](0,0.5ex)--(0.5,0.5ex);\ }
\newcommand{\ot}{\ \tikz[baseline]\draw[>=stealth,<-](0,0.5ex)--(0.5,0.5ex);\ }

\title[Spectral {M}ackey functors and equivariant algebraic {$K$}-theory {(I)}]{Spectral {M}ackey functors and equivariant algebraic {$K$}-theory {(I)}}

\author{Clark Barwick}
\address{Massachusetts Institute of Technology, Department of Mathematics, Building 2, 77 Massachusetts Avenue, Cambridge, MA 02139-4307, USA}
\email{clarkbar@gmail.com}

\dedicatory{For my dear friend Dan Kan.}

\begin{document}

\begin{abstract} \emph{Spectral Mackey functors} are homotopy-coherent versions of ordinary Mackey functors as defined by Dress. We show that they can be described as excisive functors on a suitable $\infty$-category, and we use this to show that universal examples of these objects are given by algebraic $K$-theory.

More importantly, we introduce the \emph{unfurling} of certain families of Waldhausen $\infty$-categories bound together with suitable adjoint pairs of functors; this construction completely solves the homotopy coherence problem that arises when one wishes to study the algebraic $K$-theory of such objects as spectral Mackey functors.

Finally, we employ this technology to lay the foundations of equivariant stable homotopy theory for profinite groups; the lack of such foundations has been a serious impediment to progress on the conjectures of Gunnar Carlsson. We also study fully functorial versions of $A$-theory, upside-down $A$-theory, and the algebraic $K$-theory of derived stacks.
\end{abstract}

\maketitle

\setcounter{tocdepth}{1}
\tableofcontents

%----------------------------------------------------------------------%

\setcounter{section}{-1}

\section{Summary} This paper lays the foundations of what might be called \emph{axiomatic derived representation theory}. Inspired by Bert Guillou and Peter May \cite{guillou1,guillou2,guillou3} and Dmitry Kaledin \cite{MR2918295}, we construct here a very general homotopy theory of \emph{spectral Mackey functors} --- families of spectra equipped with operations that mirror the restriction and induction operations found in ordinary Mackey functors. Our theory of spectral Mackey functors accounts for all of the compositions of these operations, their homotopies, their homotopies between homotopies, etc.

The basic input is an $\infty$-category $C$ with two subcategories, $C_{\dag}\subset C$, whose maps we call \emph{ingressive}, and $C^{\dag}\subset C$, whose maps we call \emph{egressive}. We require that ingressive and egressive maps are stable under pullback, and we require that $C$ admit finite coproducts that act effectively as disjoint unions (Df. \ref{dfn:adequate}). We call this a \emph{disjunctive triple} $(C,C_{\dag},C^{\dag})$.

\begin{exm*} The ordinary category of finite continuous $G$-sets for a profinite group $G$ defines a disjunctive triple in which every morphism is both ingressive and egressive.
\end{exm*}

\begin{exm*} The ordinary category of varieties over a field defines a disjunctive triple in which every morphism is egressive but only flat and proper maps are ingressive.
\end{exm*}

\begin{exm*} The $\infty$-category of spaces defines a disjunctive triple in which every morphism is ingressive but only morphisms with finite (homotopy) fibers are egressive.
\end{exm*}

The names ``ingressive'' and ``egressive'' are meant to suggest functorialities: a spectral Mackey functor $M$ on $C$ should consist of a covariant functor
\begin{equation*}
M_{\star}\colon\fromto{C_{\dag}}{\Sp}
\end{equation*}
and a contravariant functor
\begin{equation*}
M^{\star}\colon\fromto{(C^{\dag})^{\op}}{\Sp},
\end{equation*}
each valued in the $\infty$-category $\Sp$ of spectra. These functors will be required to carry coproducts to wedges of spectra. These two functors will agree on objects, so that given a map $f\colon\fromto{X}{Y}$, we obtain a \emph{pullback} map
\begin{equation*}
f^{\star}\colon\fromto{M(Y)=M^{\star}(Y)}{M^{\star}(X)=M(X)}
\end{equation*}
and a \emph{pushforward} map 
\begin{equation*}
f_{\star}\colon\fromto{M(X)=M_{\star}(X)}{M_{\star}(Y)=M(Y)}.
\end{equation*}
Furthermore, the pullback and pushforward maps are required to satisfy a \emph{base change condition}; namely, for any pullback square
\begin{equation*}
\begin{tikzpicture} 
\matrix(m)[matrix of math nodes, 
row sep=4ex, column sep=4ex, 
text height=1.5ex, text depth=0.25ex] 
{X'&Y'\\ 
X&Y\\}; 
\path[>=stealth,->,font=\scriptsize] 
(m-1-1) edge node[above]{$f$} (m-1-2) 
edge node[left]{$g$} (m-2-1) 
(m-1-2) edge node[right]{$g$} (m-2-2) 
(m-2-1) edge node[below]{$f$} (m-2-2); 
\end{tikzpicture}
\end{equation*}
(with abusively named morphisms), we require that
\begin{equation*}
g^{\star}f_{\star}\simeq f_{\star}g^{\star}.
\end{equation*}

But this description won't quite do as a definition of Mackey functors. After all, the base change condition is no condition at all: the homotopies $g^{\star}f_{\star}\simeq f_{\star}g^{\star}$ are additional \emph{data}, and these data have to satisfy additional coherences, which themselves are homotopies that in turn have to satisfy further coherences, etc., etc., ad infinitum. To encode all this data efficiently, we define the \emph{effective Burnside $\infty$-category} $A^{\eff}(C,C_{\dag},C^{\dag})$ (Df. \ref{dfn:geneffBurn}). The objects of this $\infty$-category are the objects of $C$, and a morphism from $X$ to $Y$ in $A^{\eff}(C,C_{\dag},C^{\dag})$ is a span
\begin{equation*}
\begin{tikzpicture}[baseline]
\matrix(m)[matrix of math nodes, 
row sep=3ex, column sep=3ex, 
text height=1.5ex, text depth=0.25ex] 
{&U&\\ 
X&&Y,\\}; 
\path[>=stealth,->,font=\scriptsize] 
(m-1-2) edge[->>] (m-2-1) 
edge[>->] (m-2-3); 
\end{tikzpicture}
\end{equation*}
in $C$ in which $\fibto{U}{X}$ is egressive and $\cofto{U}{Y}$ is ingressive. Composition is then defined by forming pullbacks. Of course, pullbacks are only unique up to a contractible choice, so composition in $A^{\eff}(C,C_{\dag},C^{\dag})$ is only defined up to a contractible choice. This is no cause for concern, however, as this is exactly the sort of thing $\infty$-categories were designed to handle gracefully. In particular, even when $C$ is an ordinary category, $A^{\eff}(C,C_{\dag},C^{\dag})$ typically won't be.

Now the $\infty$-category $A^{\eff}(C,C_{\dag},C^{\dag})$ has direct sums (Pr. \ref{prp:Aeffdirsums} and \ref{nul:Aefftripledirsums}), which are given by the coproduct in $C$, and a \emph{spectral Mackey functor} on $(C,C_{\dag},C^{\dag})$ is a functor
\begin{equation*}
\fromto{A^{\eff}(C,C_{\dag},C^{\dag})}{\Sp}
\end{equation*}
that carries this direct sum to the wedge of spectra (Df. \ref{dfn:Mackfun}). When $C$ is the category of finite $G$-sets for some finite group $G$, a spectral Mackey functor on $C$ is tantamount to a genuine $G$-equivariant spectrum. (This is a theorem of Guillou and May \cite{guillou2}.) If we replace the $\infty$-category of spectra in this discussion with an ordinary abelian category $A$, then we recover the usual notion of Mackey functors for $G$ with values in $A$ in the sense of Dress \cite{MR0360771}. So the homotopy groups of a spectral Mackey functor form an ordinary Mackey functor in abelian groups, just as one would expect. (We will actually formulate our definition in terms of general target additive $\infty$-categories.)

If $C$ is the ordinary category of finite $G$-sets for some finite group $G$, the homotopy category $hA^{\eff}(C)$ of the effective Burnside $\infty$-category is not quite what one would typically call the Burnside category. Rather, the Burnside category is obtained by forming the group completion of the Hom sets --- the ``local group completion'' --- of $A^{\eff}(C)$. (This is the significance of the word ``effective'' here; it's meant as a loose pun on the phrase ``effective divisor.'') Since our target $\infty$-categories will be additive and thus already locally group complete, however, the local group completion of $A^{\eff}(C)$ is a layer of complication we can do without.

The main result of this paper is the discovery that there is a deep connection between algebraic $K$-theory and the homotopy theory of spectral Mackey functors. We show that spectral Mackey functors can be described as excisive functors from a certain ``derived Burnside $\infty$-category'' (Lm. \ref{lem:mackisexc}). This identification allows us to use the Goodwillie differential \cite{MR2026544} to construct \emph{Mackey stabilizations} of functors from $A^{\eff}(C,C_{\dag},C^{\dag})$ to the $\infty$-category $\Kan$ of spaces (Pr. \ref{prp:Mackstabexist}). It turns out that this use of the Goodwillie differential can in important cases be related (\S \ref{sect:MackstabisK}) to the use of the Goodwillie differential in our characterization of algebraic $K$-theory \cite{K1}. This little observation permits us to express \emph{representable} Mackey functors $\SS^X$ as \emph{equivariant algebraic $K$-theory spectra}. The spectral Mackey functor $\SS$ represented by the terminal object $1$ is called the \emph{Burnside Mackey functor} (Df. \ref{dfn:BurnsideMack}); it is the analog of the sphere spectrum in this context. (In a sequel to this work, we will show that in fact it is the unit for the natural symmetric monoidal structure on Mackey functors.) Our formula for the Mackey stabilization gives us a $K$-theoretic interpretation of this Mackey functor (\S \ref{sect:BurnWaldbifib}).

More importantly, we solve the central homotopy coherence problem of equivariant algebraic $K$-theory: we prove that the algebraic $K$-theories of a family of Waldhausen $\infty$-categories \cite{K1} connected by suitable adjoint pairs of functors together define a spectral Mackey functor. This we do via an operation we call \emph{unfurling} (Df. \ref{dfn:unfurl}). The resulting structure provides a complete accounting for the functorialities enjoyed by the algebraic $K$-theory of such a family of Waldhausen $\infty$-categories.

The structure of a spectral Mackey functor is a very rich one, and it may be difficult to appreciate its utility from the abstract formalism alone. Therefore, in a sequence of appendices, we have the pleasure of studying our examples in some detail. Readers may find it useful to flip back and forth between bits of the body of the paper and these examples.
\begin{enumerate}[(A)]
\item The full subcategory of a coherent $n$-topos in the sense of Lurie ($1\leq n\leq\infty$) spanned by the coherent objects is a disjunctive triple in which every morphism is both ingressive and egressive (Ex. \ref{exm:cohntopos}). If, moreover, every coherent object can be written as a coproduct of finitely many \emph{connected objects} (Df. \ref{dfn:connectedobject}), then our $K$-theoretic description of the Burnside Mackey functor $\SS$ gives us a formula:
\begin{equation*}
\SS(1)\simeq\bigvee_{X}\Sigma_+^{\infty}B\Aut(X),
\end{equation*}
where the wedge is taken over all equivalence classes of connected objects $X$, and $\Aut(X)$ denotes the space of auto-equivalences of $X$ (Th. \ref{thm:tomDieck}). This is a very general form of the \emph{Segal--tom Dieck splitting}.
\item Suppose $G$ a profinite group. Then the (nerve of the) ordinary category of finite continuous $G$-sets (i.e., finite sets with an action of $G$ whose stabilizers are all open) is an example of the kind above. In particular, we may speak of spectral Mackey functors for $G$. When $G$ is finite, it follows from work of Bert Guillou and Peter May \cite{guillou2} that the homotopy theory of spectral Mackey functors is equivalent to the homotopy theory of $G$-equivariant spectra in the sense of Lewis--May--Steinberger \cite{MR866482}, Mandell--May \cite{MR1922205}, and Hill--Hopkins--Ravenel \cite{MR2757358,MR2906370,HHRarxiv}. For general profinite groups, we believe that our $\infty$-category of spectral Mackey functors for $G$ is the first definition of the homotopy theory of $G$-equivariant spectra. (The lack of foundations for such a subject has been a serious impediment to real progress on the conjectures of Gunnar Carlsson.) The generalized Segal--tom Dieck splitting (Th. \ref{thm:tomDieck}) gives a formula for the $G$-fixed points of the $G$-equivariant sphere spectrum:
\begin{equation*}
\SS_G^G\simeq\bigvee_{H}\Sigma_+^{\infty}B(N_GH/H),
\end{equation*}
where the wedge is taken over conjugacy classes of open subgroups $H\leq G$, and $N_GH$ denotes the normalizer of $H$ in $G$ (Pr. \ref{prp:SegtomDieckprofinite}).
\item For any space $X$, there are two kinds of conditions one might impose on a retractive spaces $X\to X'\to X$ (Nt. \ref{ntn:IXandJX}): one could demand that $X'$ is a retract of a finite CW complex (or, alternately, a finite CW complex itself); alternately, one could demand that the homotopy fibers of $\fromto{X'}{X}$ are retracts of finite CW complexes (or, alternately, finite CW complexes). The algebraic $K$-theories of these $\infty$-categories are spectra denoted $\AA(X)$ and $\upsidedown{\AA}(X)$. The former is covariantly functorial and the latter is contravariantly functorial. But each possesses additional functorialities: if $\fromto{X}{Y}$ is a map whose homotopy fibers are retracts of finite CW complexes, then we obtain umkehr maps $\fromto{\AA(Y)}{\AA(X)}$ and $\fromto{\upsidedown{\AA}(X)}{\upsidedown{\AA}(Y)}$. We can assemble these functorialities together to get Mackey functors $\AA$ (respectively, $\upsidedown{\AA}$) for the disjunctive triple given by the $\infty$-category of spaces in which every map is ingressive (resp., egressive), and only those maps whose homotopy fibers are retracts of finite CW complexes are egressive (resp., ingressive) (Nt. \ref{ntn:AandupsidedownA}). This provides a host of interesting \emph{assembly maps} (\ref{nul:Aassembly})
\begin{equation*}
\fromto{\AA^{\oplus}(X,Y)\wedge\AA(X)}{\AA(Y)}\textrm{\quad and\quad}\fromto{\AA^{\oplus}(X,Y)\wedge\upsidedown{\AA}(Y)}{\upsidedown{\AA}(X)},
\end{equation*}
where $\AA^{\oplus}(X,Y)$ is the group completion of the $E_{\infty}$ space of diagrams
\begin{equation*}
X\ot U\to Y,
\end{equation*}
where the homotopy fibers of $\fromto{U}{X}$ are finite CW complexes, and the $E_{\infty}$ structure is given by coproduct. When $X=\ast$, the maps $\SS\to\AA(\ast)\to\SS$ can be composed with these assembly maps to obtain, for any retract $U$ of a finite CW complex, maps
\begin{equation*}
\fromto{\Sigma_{+}^{\infty}\Map(U,Y)}{\AA(Y)}\textrm{\quad and\quad}\fromto{\Sigma_{+}^{\infty}\Map(U,Y)}{D\upsidedown{\AA}(Y)}.
\end{equation*}
Dually, when $Y=\ast$, we obtain, for any space $V$, maps
\begin{equation*}
\fromto{\Sigma_{+}^{\infty}\Map_{\mathrm{rc}}(X,V)}{D\AA(X)}\textrm{\quad and\quad}\fromto{\Sigma_{+}^{\infty}\Map_{\mathrm{rc}}(X,V)}{\upsidedown{\AA}(X)},
\end{equation*}
where $\Map_{\mathrm{rc}}$ denotes the space of maps whose homotopy fibers are retracts of finite CW complexes. Special cases of these maps have been constructed by Waldhausen, Malkiewich, and others. Moreover, it turns out that all of this holds when the $\infty$-category of spaces is replaced with any compactly generated $\infty$-topos whose terminal object is compact.
\item In \cite[\S 12]{K1}, we defined the algebraic $K$-theory of derived stacks, and we observed there that it was contravariantly functorial in morphisms of (nonconnective) spectral Deligne--Mumford stacks. Here we push this further by including the covariant functoriality of the algebraic $K$-theory of spectral Deligne--Mumford stacks in certain and its compatibility with the contravariant functoriality. More precisely, we construct two disjunctive triples of derived stacks. The first is the $\infty$-category of spectral Deligne--Mumford stacks in which we declare that every morphism is egressive, and a morphism is ingressive if and only if it is strongly proper, of finite Tor-amplitude, and locally of finite presentation (\ref{ntn:properness}). The second is the $\infty$-category of all flat sheaves (of spaces) on the $\infty$-category of connective $E_{\infty}$ rings in which we declare that every morphism is egressive, and a morphism is ingressive if and only if it is quasi-affine representable and \emph{perfect} in the sense that the pushforward of the structure sheaf is perfect (Pr. \ref{prp:perf}). We prove that algebraic $K$-theory is a spectral Mackey functor for each of these disjunctive triples (Nt. \ref{ntn:KofDMstacks} and Nt. \ref{ntn:KasMackonflatsheaves}). This can be thought of as a very general and very structured form of proper base change for $K$-theory. It also ensures the existence of interesting \emph{assembly maps}
\begin{equation*}
\fromto{\AA_{\mathrm{DM}}^{\oplus}(X,Y)\wedge\KK(X)}{\KK(Y)}\textrm{\quad (respectively,\ }\fromto{\AA_{\mathrm{Shv}}^{\oplus}(X,Y)\wedge\KK(X)}{\KK(Y)}\textrm{\ ),}
\end{equation*}
where $\AA_{\mathrm{DM}}^{\oplus}(X,Y)$ (resp., $\AA_{\mathrm{Shv}}^{\oplus}(X,Y)$) is the group completion of the $E_{\infty}$ space of diagrams
\begin{equation*}
X\ot U\to Y
\end{equation*}
of spectral Deligne--Mumford stacks (resp., of flat sheaves), where $\fromto{U}{Y}$ is strongly proper, of finite Tor-amplitude, and locally of finite presentation (resp., quasi-affine representable and perfect), and the $E_{\infty}$ structure is given by coproduct. When $X=\Spec\SS$, the maps $\SS\to\KK(\SS)\to\SS$ can be composed with these assembly maps to obtain, for any spectral Deligne--Mumford stack (resp., any flat sheaf) $U$, maps
\begin{equation*}
\fromto{\Sigma_{+}^{\infty}\Map_{\mathrm{pr}}(U,Y)}{\KK(Y)}\textrm{\quad and\quad}\fromto{\Sigma_{+}^{\infty}\Map_{\mathrm{perf}}(U,Y)}{\KK(Y)},
\end{equation*}
where $\Map_{\mathrm{pr}}(U,Y)$ (resp., $\Map_{\mathrm{perf}}(U,Y)$) denotes the space of morphisms $\fromto{U}{Y}$ that are strongly proper, of finite Tor-amplitude, and locally of finite presentation (resp., quasi-affine representable and perfect). Dually, when $Y=\ast$, we obtain, for any spectral Deligne--Mumford stack (resp., any flat sheaf) $V$ that is strongly proper, of finite Tor-amplitude, and locally of finite presentation (resp., quasi-affine representable and perfect) over $\Spec\SS$, maps
\begin{equation*}
\fromto{\Sigma_{+}^{\infty}\Map(V,X)}{D\KK(X)}\textrm{\quad and\quad}\fromto{\Sigma_{+}^{\infty}\Map(V,X)}{D\KK(X)}.
\end{equation*}
Restricting this to \'etale covers of a fixed (nice) scheme $X$ with a geometric point $x\in X(\Omega)$, we obtain the \emph{Galois equivariant $K$-theory spectrum} (Ex. \ref{exm:GaloisequivariantK})
\begin{equation*}
K_{\pi_1^{\et}(X,x)}(X)\colon\fromto{B_{\pi_1^{\et}(X,x)}^{\mathit{fin}}}{\Sp}.
\end{equation*}
In a very precise manner, this object encodes the failure of $K$-theory to satisfy Galois descent; this is the subject of a series of conjectures of Gunnar Carlsson \cite{carlss:repass}, which can now be formulated thanks to the foundations we provide here. We intend to study this object in great detail in future work.
\end{enumerate}

\subsection*{Related work} It is probably safe to say that others have anticipated a fully $\infty$-categorical (and thus ``model-independent'') construction of equivariant stable homotopy theory, but, as always, the devil is in the details. We believe this is the first $\infty$-categorical approach to equivariant stable homotopy theory, and it is the first to provide a complete construction of equivariant algebraic $K$-theory of families of Waldhausen $\infty$-categories of the kind described above. There are, however, a number of precursors to this paper.

In a brilliant series of papers, Bert Guillou and Peter May construct \cite{guillou1,guillou2,guillou3} a homotopy theory of spectral Mackey functors, and they show that for a finite group $G$, the homotopy theory of spectral Mackey functors for finite $G$ is equivalent to the homotopy theory of genuine $G$-spectra. Our homotopy theory of spectral Mackey functors for finite $G$-sets is easily seen to be equivalent to theirs.

Dmitry Kaledin has also developed \cite{MR2918295} a theory of ``derived Mackey functors'' (again for finite groups $G$). Our work here is a generalization of Kaledin's: the homotopy category of Mackey functors valued in the derived $\infty$-category of abelian groups is naturally equivalent to his derived category of Mackey functors.

Aderemi Kuku, extending work of Andreas Dress \cite{MR0384917}, has worked for many years on the foundations of equivariant higher algebraic $K$-theory as a construction that yields ordinary abelian group-valued Mackey functors (partly joint with Dress) \cite{MR609220,MR689366,MR750684,MR1745585,MR2175639,MR2310576,MR2259035}. For Mackey functors for finite groups, our work can be understood as lifting the target of Kuku's constructions to the $\infty$-category of spectrum-valued Mackey functors.

Work of Shimakawa \cite{MR1132161,MR1162447,MR1162448,MR1233745,MR1295582} shows that the $K$-theory of permutative categories with a suitable action of a finite group can be given the structure of a genuine $G$-spectrum. The work here amounts to the generalization of this result to the context of Waldhausen $K$-theory for $\infty$-categories with suitable $G$-actions.

Mona Merling has an alternate approach to constructing $G$-spectra from Waldhausen categories with a suitable action of a finite group $G$. It is possible that her approach and the one given here are suitably equivalent; however, it seems that the two approaches differ significantly in the details, and Merling's appears to be adapted to the technology of equivariant homotopy theory as developed by Lewis--May--Steinberger \cite{MR866482}, Mandell--May \cite{MR1922205}, and Hill--Hopkins--Ravenel \cite{MR2757358,MR2906370,HHRarxiv}.

\subsection*{Acknowledgments} I have had very helpful conversations with David Ayala and Haynes Miller about the contents of this paper and its sequels.

%----------------------------------------------------------------------%
%----------------------------------------------------------------------%
%----------------------------------------------------------------------%

\section{Preliminaries on $\infty$-categories} In general, we use the terminology from \cite{HTT,K1,K21}. We review some of the relevant notation here.

\begin{ntn} In order to deal gracefully with size issues, we'll use Grothendieck universes in this paper. In particular, we fix, once and for all, three uncountable, strongly inaccessible cardinals $\kappa_0<\kappa_1<\kappa_2$ and the corresponding universes $\VV_{\kappa_0}\in\VV_{\kappa_1}\in\VV_{\kappa_2}$. Now a set, simplicial set, category, etc., will be said to be \textbf{\emph{small}} if it is contained in the universe $\VV_{\kappa_0}$; it will be said to be \textbf{\emph{large}} if it is contained in the universe $\VV_{\kappa_1}$; and it will be said to be \textbf{\emph{huge}} if it is contained in the universe $\VV_{\kappa_2}$. We will say that a set, simplicial set, category, etc., is \textbf{\emph{essentially small}} if it is equivalent (in the appropriate sense) to a small one.
\end{ntn}

\begin{nul} The model of $\infty$-categories we will employ is Joyal's model of \emph{quasicategories}, which we will here call \textbf{\emph{$\infty$-categories}}. We refer systematically to \cite{HTT} for details about this model of higher categories. 
\end{nul}

\begin{ntn}\label{ntn:superscriptsscats} A \textbf{\emph{simplicial category}} --- that is, a category enriched in the category of simplicial sets --- will frequently be denoted with a superscript $(-)^{\Delta}$.

Suppose $\CC^{\Delta}$ a simplicial category. Then we write $(\CC^{\Delta})_{0}$ for the ordinary category given by taking the $0$-simplices of the $\Mor$ spaces. That is, $(\CC^{\Delta})_{0}$ is the category whose objects are the objects of $\CC$, and whose morphisms are given by
\begin{equation*}
(\CC^{\Delta})_{0}(x,y)\coloneq\CC^{\Delta}(x,y)_0.
\end{equation*}
If the $\Mor$ spaces of $\CC^{\Delta}$ are all fibrant, then we will often write
\begin{equation*}
\CC\textrm{\quad for the simplicial nerve\quad}N(\CC^{\Delta})
\end{equation*}
\cite[Df. 1.1.5.5]{HTT}, which is an $\infty$-category \cite[Pr. 1.1.5.10]{HTT}.
\end{ntn}

\begin{ntn}\label{ntn:interior} For any $\infty$-category $A$, there exists a simplicial subset $\iota A\subset A$, which is the largest Kan simplicial subset of $A$ \cite[1.2.5.3]{HTT}. We shall call this space the \textbf{\emph{interior $\infty$-groupoid of $A$}}. The assignment $\goesto{A}{\iota A}$ defines a right adjoint $\iota$ to the inclusion functor $u$ from Kan simplicial sets to $\infty$-categories.
\end{ntn}

\begin{ntn}\label{ntn:catofcat} The large simplicial category $\Kan^{\Delta}$ is the category of small Kan simplicial sets, with the usual notion of mapping space. The large simplicial category $\Cat^{\Delta}_{\infty}$ is defined in the following manner \cite[Df. 3.0.0.1]{HTT}. The objects of $\Cat^{\Delta}_{\infty}$ are small $\infty$-categories, and for any two $\infty$-categories $A$ and $B$, the morphism space
\begin{equation*}
\Cat^{\Delta}_{\infty}(A,B)\coloneq\iota\Fun(A,B)
\end{equation*}
is the interior $\infty$-groupoid of the $\infty$-category $\Fun(A,B)$.

Similarly, for any strongly inaccessible cardinal $\tau$, we may define the locally $\tau$-small simplicial category $\Kan(\tau)^{\Delta}$ of $\tau$-small simplicial sets and the locally $\tau$-small simplicial category $\Cat_{\infty}(\tau)^{\Delta}$ of $\tau$-small $\infty$-categories.
\end{ntn}

%----------------------------------------------------------------------%

\section{The twisted arrow $\infty$-category} We have elsewhere \cite{} spoken of the twisted arrow $\infty$-category of an $\infty$-category. Let us recall the basic facts here.

\begin{prp}\label{prp:subdivision} The following are equivalent for a functor $\theta\colon\fromto{\Delta}{\Delta}$.
\begin{enumerate}[(\ref{prp:subdivision}.1)]
\item The functor $\theta^{\op}\colon\fromto{N\Delta^{\op}}{N\Delta^{\op}}$ is cofinal in the sense of Joyal \cite[Df. 4.1.1.1]{HTT}.
\item The induced endofunctor $\theta^{\star}\colon\fromto{s\Set}{s\Set}$ on the ordinary category of simplicial sets (so that $(\theta^{\star}X)_n=X_{\theta(n)}$) carries every standard simplex $\Delta^m$ to a weakly contractible simplicial set.
\item The induced endofunctor $\theta^{\star}\colon\fromto{s\Set}{s\Set}$ on the ordinary category of simplicial sets is a left Quillen functor for the usual Quillen model structure.
\end{enumerate}
\begin{proof} By Joyal's variant of Quillen's Theorem A \cite[Th. 4.1.3.1]{HTT}, the functor $\theta^{\op}$ is cofinal just in case, for any integer $m\geq 0$, the nerve $N(\theta/\mathbf{m})$ is weakly contractible. The category $(\theta/\mathbf{m})$ is clearly equivalent to the category of simplices of $\theta^{\star}(\Delta^m)$, whose nerve is weakly equivalent to $\theta^{\star}(\Delta^m)$. This proves the equivalence of the first two conditions.

It is clear that for any functor $\theta\colon\fromto{\Delta}{\Delta}$, the induced functor $\theta^{\star}\colon\fromto{s\Set}{s\Set}$ preserves monomorphisms. Hence $\theta^{\star}$ is left Quillen just in case it preserves weak equivalences. Hence if $\theta^{\star}$ is left Quillen, then it carries the map $\equivto{\Delta^n}{\Delta^0}$ to an equivalence $\equivto{\theta^{\star}\Delta^n}{\theta^{\star}\Delta^0\cong\Delta^0}$, and, conversely, if $\theta^{\op}\colon\fromto{N\Delta^{\op}}{N\Delta^{\op}}$ is cofinal, then for any weak equivalence $\equivto{X}{Y}$, the induced map $\fromto{\theta^{\star}X}{\theta^{\star}Y}$ factors as
\begin{eqnarray}
\theta^{\star}X\simeq\hocolim_nX_{\theta(n)}&\simeq&\hocolim_nP_N\nonumber\\
&\simeq&X\nonumber\\
&\equivto{}{}&Y\nonumber\\
&\simeq&\hocolim_nY_n\nonumber\\
&\simeq&\hocolim_nY_{\theta(n)}\simeq\theta^{\star}Y,\nonumber
\end{eqnarray}
which is a weak equivalence. This proves the equivalence of the third condition with the first two.
\end{proof}
\end{prp}

\begin{nul} One may call any functor $\theta\colon\fromto{\Delta}{\Delta}$ satisfying the equivalent conditions above a \textbf{\emph{combinatorial subdivision}}. Work of Katerina Velcheva shows that in fact combinatorial subdivisions can be classified: they are all iterated joins of $\id$ and $\op$. The example in which we are interested, the join $\op\star\id$, is originally due to Segal.
\end{nul}

\begin{ntn}\label{exm:edgewise} Denote by $\epsilon\colon\fromto{\Delta}{\Delta}$ the combinatorial subdivision
\begin{equation*}
\goesto{\mathbf{[n]}}{\mathbf{[n]}^{\op}\star\mathbf{[n]}\cong\mathbf{[2n+1]}}.
\end{equation*}
Including $\mathbf{[n]}$ into either factor of the join $\mathbf{[n]}^{\op}\star\mathbf{[n]}$ (either contravariantly or covariantly) defines two natural transformations $\fromto{\op}{\epsilon}$ and $\fromto{\id}{\epsilon}$. Precomposition with $\epsilon$ induces an endofunctor $\epsilon^{\star}$ on the ordinary category of simplicial sets, together with natural transformations $\fromto{\epsilon^{\star}}{\op}$ and $\fromto{\epsilon^{\star}}{\id}$.

For any simplicial set $X$, the \textbf{\emph{edgewise subdivision}} of $X$ is the simplicial set
\begin{equation*}
\widetilde{\mathscr{O}}(X)\coloneq\epsilon^{\star}X.
\end{equation*}
That is, $\widetilde{\mathscr{O}}(X)$ is given by the formula
\begin{equation*}
\widetilde{\mathscr{O}}(X)_n=\Mor(\Delta^{n,\op}\star\Delta^n,X)\cong X_{2n+1}.
\end{equation*}
The two natural transformations described above give rise to a morphism
\begin{equation*}
\fromto{\widetilde{\mathscr{O}}(X)}{X^{\op}\times X},
\end{equation*}
functorial in $X$.
\end{ntn}

\begin{nul}\label{nul:twarrofnerveisverve} For any simplicial set $X$, the vertices of $\widetilde{\mathscr{O}}(X)$ are edges of $X$; an edge of $\widetilde{\mathscr{O}}(X)$ from $\fromto{u}{v}$ to $\fromto{x}{y}$ can be viewed as a commutative diagram (up to chosen homotopy)
\begin{equation*}
\begin{tikzpicture} 
\matrix(m)[matrix of math nodes, 
row sep=4ex, column sep=4ex, 
text height=1.5ex, text depth=0.25ex] 
{u&x\\ 
v&y\\}; 
\path[>=stealth,->,font=\scriptsize] 
(m-1-2) edge (m-1-1) 
edge (m-2-2)
(m-1-1) edge (m-2-1)
(m-2-1) edge (m-2-2); 
\end{tikzpicture}
\end{equation*}
When $X$ is the nerve of an ordinary category $C$, $\widetilde{\mathscr{O}}(X)$ is isomorphic to the nerve of the twisted arrow category of $C$ in the sense of \cite{MR705421}. When $X$ is an $\infty$-category, we are therefore inclined to call $\widetilde{\mathscr{O}}(X)$ the \textbf{\emph{twisted arrow $\infty$-category}} of $X$. This terminology is justified by the following.
\end{nul}

\begin{prp}[Lurie, \protect{\cite[Pr. 4.2.3]{DAGX}}]\label{prp:twarrisinfincat} If $X$ is an $\infty$-category, then the functor $\fromto{\widetilde{\mathscr{O}}(X)}{X^{\op}\times X}$ is a left fibration; in particular, $\widetilde{\mathscr{O}}(X)$ is an $\infty$-category.
\end{prp}

\begin{exm} To illustrate, for any object $\mathbf{p}\in\Delta$, the $\infty$-category $\widetilde{\mathscr{O}}(\Delta^p)$ is the nerve of the category
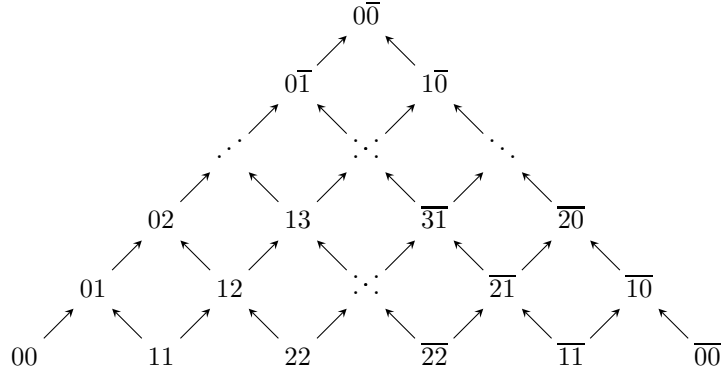
\begin{equation*}
\begin{tikzpicture} 
\matrix(m)[matrix of math nodes, 
row sep={6ex,between origins}, column sep={6ex,between origins}, 
text height=1.5ex, text depth=0.25ex] 
{&&&&&0\overline{0}&&&&&\\
&&&&0\overline{1}&&1\overline{0}&&&&\\
&&&\revddots&&\node{\ddots};\node{\revddots};&&\ddots&&&\\
&&02&&13&&\overline{3}\overline{1}&&\overline{2}\overline{0}&&\\
&01&&12&&\node{\ddots};\node{\revddots};&&\overline{2}\overline{1}&&\overline{1}\overline{0}&\\
00&&11&&22&&\overline{2}\overline{2}&&\overline{1}\overline{1}&&\overline{0}\overline{0}\\}; 
\path[>=stealth,->,font=\scriptsize] 
(m-2-5) edge (m-1-6) 
(m-2-7) edge (m-1-6)
(m-3-4) edge (m-2-5)
(m-3-6) edge (m-2-5)
(m-3-6) edge (m-2-7)
(m-3-8) edge (m-2-7)
(m-4-3) edge (m-3-4)
(m-4-5) edge (m-3-4)
(m-4-5) edge (m-3-6)
(m-4-7) edge (m-3-6)
(m-4-7) edge (m-3-8)
(m-4-9) edge (m-3-8)
(m-5-2) edge (m-4-3)
(m-5-4) edge (m-4-3)
(m-5-4) edge (m-4-5)
(m-5-6) edge (m-4-5)
(m-5-6) edge (m-4-7)
(m-5-8) edge (m-4-7)
(m-5-8) edge (m-4-9)
(m-5-10) edge (m-4-9)
(m-6-1) edge (m-5-2)
(m-6-3) edge (m-5-2)
(m-6-3) edge (m-5-4)
(m-6-5) edge (m-5-4)
(m-6-5) edge (m-5-6)
(m-6-7) edge (m-5-6)
(m-6-7) edge (m-5-8)
(m-6-9) edge (m-5-8)
(m-6-9) edge (m-5-10) 
(m-6-11) edge (m-5-10);
\end{tikzpicture}
\end{equation*}
(Here we write $\overline{n}$ for $p-n$.)
\end{exm}

%----------------------------------------------------------------------%

\section{The effective Burnside $\infty$-category} 

We now employ the edgewise subdivision to define a quasicategorical variant of the Burnside category. The essence of the idea was explored in our work on the $\infty$-categorical $Q$ construction.

\begin{ntn}\label{dfn:Rstar} For any $\infty$-category $C$, denote by $\RR_{\ast}(C)\colon\fromto{\Delta^{\op}}{s\Set}$ the functor given by the assignment
\begin{equation*}
\goesto{\mathbf{[n]}}{\iota\Fun(\widetilde{\mathscr{O}}(\Delta^n)^{\op},C)}.
\end{equation*}
\end{ntn}

\begin{prp}\label{prp:Rstarisgood} The functor $\RR_{\ast}\colon\fromto{\Cat_{\infty}^{0}}{\Fun(\Delta^{\op},s\Set)}$ carries every quasicategory to a Reedy fibrant simplicial space, and it preserves weak equivalences.
\begin{proof} We first show that for any $\infty$-category $C$, the simplicial space $\RR_{\ast}(C)$ is Reedy fibrant. This is the condition that for any monomorphism $\into{K}{L}$, the map
\begin{equation*}
\fromto{\iota\Fun(\widetilde{\mathscr{O}}(L)^{\op},C)}{\iota\Fun(\widetilde{\mathscr{O}}(K)^{\op},C)}
\end{equation*}
is a Kan fibration of simplicial sets. This follows immediately from Pr. \ref{prp:subdivision} and \cite[Lm. 3.1.3.6]{HTT}. To see that $\RR_{\ast}$ preserves weak equivalences, we note that since $\Fun(\widetilde{\mathscr{O}}(\Delta^n)^{\op},-)$ preserves weak equivalences, so does $\RR_n$.
\end{proof}
\end{prp}

\begin{dfn}\label{dfn:cartesianfunct} Suppose $C$ an $\infty$-category. For any integer $n\geq 0$, let us say that a functor $X\colon\fromto{\widetilde{\mathscr{O}}(\Delta^n)^{\op}}{C}$ is \textbf{\emph{cartesian}} if, for any integers $0\leq i\leq k\leq\ell\leq j\leq n$, the square
\begin{equation*}
\begin{tikzpicture} 
\matrix(m)[matrix of math nodes, 
row sep=4ex, column sep=4ex, 
text height=1.5ex, text depth=0.25ex] 
{X_{ij}&X_{kj}\\ 
X_{i\ell}&X_{k\ell}\\}; 
\path[>=stealth,->,font=\scriptsize] 
(m-1-1) edge (m-1-2) 
edge (m-2-1) 
(m-1-2) edge (m-2-2) 
(m-2-1) edge (m-2-2); 
\end{tikzpicture}
\end{equation*}
is a pullback. 

Write $\AA^{\eff}_{\ast}(C)\subset\RR_{\ast}(C)$ for the subfunctor in which $\AA^{\eff}_n(C)$ is the full simplicial subset of $\RR_n(C)$ spanned by the cartesian functors
\begin{equation*}
X\colon\fromto{\widetilde{\mathscr{O}}(\Delta^n)^{\op}}{C}.
\end{equation*}
Note that since any functor that is equivalent to an cartesian functor is itself cartesian, the simplicial set $\AA^{\eff}_n(C)$ is a union of connected components of $\RR_n(C)$.
\end{dfn}

\begin{prp}\label{prp:QQstarisCSS} For any $\infty$-category $C$ that admits all pullbacks, the simplicial space $\AA^{\eff}_{\ast}(C)$ is a complete Segal space.
\begin{proof} The Reedy fibrancy of $\AA^{\eff}_{\ast}(C)$ follows easily from the Reedy fibrancy of $\RR_{\ast}(C)$.

To see that $\AA^{\eff}_{\ast}(C)$ is a Segal space, it is necessary to show that for any integer $n\geq 1$, the Segal map
\begin{equation*}
\fromto{\AA^{\eff}_n(C)}{\AA^{\eff}_1(C)\times_{\AA^{\eff}_0(C)}\cdots\times_{\AA^{\eff}_0(C)}\AA^{\eff}_1(C)}
\end{equation*}
is an equivalence. Let $L_n$ denote the ordinary category
\begin{equation*}
00\ot 01\to 11\ot 12\to\cdots (n-1)(n-1)\ot (n-1)n\to nn.
\end{equation*}
The target of the Segal map can then be identified with the maximal Kan complex contained in the $\infty$-category
\begin{equation*}
\Fun(NL_n,C).
\end{equation*}
The Segal map is therefore an equivalence by the uniqueness of limits in $\infty$-categories \cite[Pr. 1.2.12.9]{HTT}.

Finally, to check that $\AA^{\eff}_{\ast}(C)$ is complete, let $E$ be the nerve of the contractible ordinary groupoid with two objects; then completeness is equivalent to the assertion that the Rezk map
\begin{equation*}
\fromto{\AA^{\eff}_0(C)}{\lim_{\mathbf{[n]}\in(\Delta/E)^{\op}}\AA^{\eff}_n(C)}
\end{equation*}
is a weak equivalence. The source of this map can be identified with $\iota C$; its target can be identified with the full simplicial subset of 
\begin{equation*}
\iota\Fun(\widetilde{\mathscr{O}}(E)^{\op},C)
\end{equation*}
spanned by those functors $X\colon\fromto{\widetilde{\mathscr{O}}(E)^{\op}}{C}$ such that for any simplex $\fromto{\Delta^n}{E}$, the induced functor $\fromto{\widetilde{\mathscr{O}}(\Delta^n)^{\op}}{C}$ is ambigressive. Note that the twisted arrow category of the contractible ordinary groupoid with two objects is the contractible ordinary groupoid with four objects. Consequently, the image of any functor $X\colon\fromto{\widetilde{\mathscr{O}}(E)^{\op}}{C}$ is contained in $\iota C$. Thus the target of the Rezk map can be identified with $\iota\Fun(\widetilde{\mathscr{O}}(E)^{\op},C)$ itself, and the Rezk map is an equivalence.
\end{proof}
\end{prp}

\begin{nul} It is now clear that $\AA^{\eff}_{\ast}$ defines a functor of $\infty$-categories
\begin{equation*}
\AA^{\eff}_{\ast}\colon\fromto{\Cat_{\infty}^{\mathrm{lex}}}{\mathbf{CSS}},
\end{equation*}
where $\Cat_{\infty}^{\mathrm{lex}}\subset\Cat_{\infty}$ is the subcategory consisting of $\infty$-categories with all finite limits and left exact functors between them, and where $\mathbf{CSS}\subset\Fun(\Delta^{\op},\Kan)$ is the full subcategory spanned by complete Segal spaces.
\end{nul}

Joyal and Tierney show that the functor that carries a simplicial space $\XX$ to the simplicial set whose $n$-simplices are the vertices of $\XX_n$ induces an equivalence of $\infty$-categories $\fromto{\mathbf{CSS}}{\Cat_{\infty}}$. This leads us to the following definition.

\begin{dfn} For any $\infty$-category $C$ that admits all pullbacks, denote by $A^{\eff}(C)$ the $\infty$-category whose $n$-simplices are vertices of $\AA^{\eff}_n(C)$, i.e., cartesian functors $\fromto{\widetilde{\mathscr{O}}(\Delta^n)^{\op}}{C}$. We may call this the \textbf{\emph{effective Burnside $\infty$-category of $C$}}.

This defines a functor of $\infty$-categories
\begin{equation*}
A^{\eff}\colon\fromto{\Cat_{\infty}^{\mathrm{lex}}}{\Cat_{\infty}}. 
\end{equation*}
\end{dfn}

\begin{nul} For any $\infty$-category $C$ that admits all pullbacks, an $n$-simplex of $A^{\eff}(C)$ is a diagram
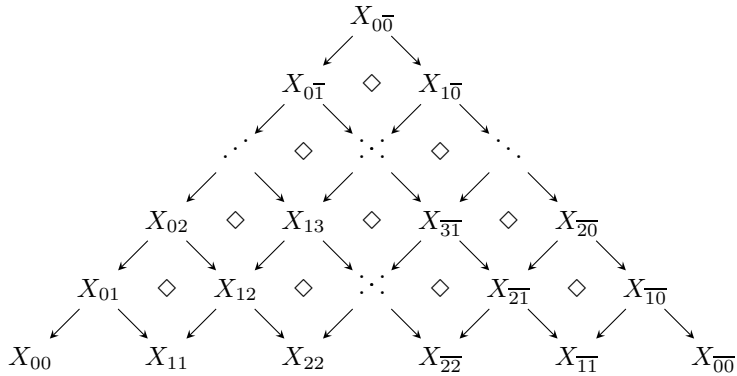
\begin{equation*}
\begin{tikzpicture} 
\matrix(m)[matrix of math nodes, 
row sep={6ex,between origins}, column sep={6ex,between origins}, 
text height=1.5ex, text depth=0.25ex] 
{&&&&&X_{0\overline{0}}&&&&&\\
&&&&X_{0\overline{1}}&\Diamond&X_{1\overline{0}}&&&&\\
&&&\revddots&\Diamond&\node{\ddots};\node{\revddots};&\Diamond&\ddots&&&\\
&&X_{02}&\Diamond&X_{13}&\Diamond&X_{\overline{3}\overline{1}}&\Diamond&X_{\overline{2}\overline{0}}&&\\
&X_{01}&\Diamond&X_{12}&\Diamond&\node{\ddots};\node{\revddots};&\Diamond&X_{\overline{2}\overline{1}}&\Diamond&X_{\overline{1}\overline{0}}&\\
X_{00}&&X_{11}&&X_{22}&&X_{\overline{2}\overline{2}}&&X_{\overline{1}\overline{1}}&&X_{\overline{0}\overline{0}}\\}; 
\path[>=stealth,<-,font=\scriptsize] 
(m-2-5) edge (m-1-6) 
(m-2-7) edge (m-1-6)
(m-3-4) edge (m-2-5)
(m-3-6) edge (m-2-5)
(m-3-6) edge (m-2-7)
(m-3-8) edge (m-2-7)
(m-4-3) edge (m-3-4)
(m-4-5) edge (m-3-4)
(m-4-5) edge (m-3-6)
(m-4-7) edge (m-3-6)
(m-4-7) edge (m-3-8)
(m-4-9) edge (m-3-8)
(m-5-2) edge (m-4-3)
(m-5-4) edge (m-4-3)
(m-5-4) edge (m-4-5)
(m-5-6) edge (m-4-5)
(m-5-6) edge (m-4-7)
(m-5-8) edge (m-4-7)
(m-5-8) edge (m-4-9)
(m-5-10) edge (m-4-9)
(m-6-1) edge (m-5-2)
(m-6-3) edge (m-5-2)
(m-6-3) edge (m-5-4)
(m-6-5) edge (m-5-4)
(m-6-5) edge (m-5-6)
(m-6-7) edge (m-5-6)
(m-6-7) edge (m-5-8)
(m-6-9) edge (m-5-8)
(m-6-9) edge (m-5-10) 
(m-6-11) edge (m-5-10);
\end{tikzpicture}
\end{equation*}
of $C$ in which every square is a pullback. Here we write $\overline{n}$ for $p-n$.

Another way of describing $A^{\eff}(C)$ is as follows. The objects of $A^{\eff}(C)$ are precisely those of $C$. Between objects $X$ and $Y$, the space of maps is given by
\begin{equation*}
\Map_{A^{\eff}(C)}(X,Y)\simeq\iota C_{/\{X,Y\}},
\end{equation*}
where $\{X,Y\}$ denotes the diagram $\fromto{\{x,y\}}{C}$ from the discrete simplicial set $\{x,y\}$ to $C$ that carries $x$ to $X$ and $y$ to $Y$. Composition
\begin{equation*}
\fromto{\iota C_{/\{X,Y\}}\times\iota C_{/\{Y,Z\}}}{\iota C_{/\{X,Z\}}}
\end{equation*}
is defined, up to coherent homotopy, by pullback $-\times_Y-$.
\end{nul}

\begin{nul}\label{nul:ordinaryBurnside} Note that the traditional Burnside category is distinct from $A^{\eff}(C)$ in two ways. First, when forming the traditional Burnside category, one begins by studying \emph{isomorphism classes} of spans between objects. This is to ensure that one obtains a nice set of maps for which the pullback construction is sensible. In our effective Burnside $\infty$-category, we do not pass to isomorphism classes. Rather, we are content to take the entire \emph{space} of spans between objects as our mapping space. We again use pullback to define composition, and we lose no sleep over the fact that pullbacks are only defined up to coherent equivalence, since composition in any $\infty$-category is only required to be defined up to coherent equivalence in the first place. The \textbf{\emph{ordinary effective Burnside category}} of an ordinary category $C$ may be identified with the homotopy category $hA^{\eff}(NC)$ of $A^{\eff}(NC)$.

Second, the ordinary Burnside category is usually defined as the ``local group completion'' of this ordinary effective Burnside category. Then Mackey functors are then defined as additive functors from this Burnside category to, say, the category of abelian groups. This is overkill: if the target is already group complete, one knows already what additive functors from the Burnside category will be in terms of the category before group completion. The group completion is a relatively minor procedure for ordinary categories, but for $\infty$-categories, group completion is serious business. Indeed, if $\FF$ is the ordinary category of finite sets, then when one forms the local group completion $A(N\FF)$ of $A^{\eff}(N\FF)$, the space of endomorphisms on the one point set becomes
\begin{equation*}
\Map_{A(N\FF)}(\ast,\ast)\simeq QS^0,
\end{equation*}
by Barratt--Priddy--Quillen. To avoid such complications, we happily stick with the effective Burnside $\infty$-category.
\end{nul}

\begin{ntn} The two natural transformations $\fromto{\varepsilon^{\star}}{\op}$ and $\fromto{\varepsilon^{\star}}{\id}$ induce two natural transformations
\begin{equation*}
(\cdot)_{\star}\colon\fromto{\id}{A^{\eff}}\textrm{\quad and\quad}(\cdot)^{\star}\colon\fromto{\op}{A^{\eff}}.
\end{equation*}
For any morphism $f\colon\fromto{U}{V}$ of an $\infty$-category $C$ that admits all pullbacks, one thus obtains morphisms
\begin{equation*}
f_{\star}\colon\fromto{U}{V}\textrm{\quad and\quad}f^{\star}\colon\fromto{V}{U}.
\end{equation*}
For any pullback square
\begin{equation*}
\begin{tikzpicture} 
\matrix(m)[matrix of math nodes, 
row sep=4ex, column sep=4ex, 
text height=1.5ex, text depth=0.25ex] 
{U&X\\ 
V&Y,\\}; 
\path[>=stealth,->,font=\scriptsize] 
(m-1-1) edge node[above]{$i$} (m-1-2) 
edge node[left]{$f$} (m-2-1) 
(m-1-2) edge node[right]{$g$} (m-2-2) 
(m-2-1) edge node[below]{$j$} (m-2-2); 
\end{tikzpicture}
\end{equation*}
one obtains a homotopy
\begin{equation*}
g^{\star}\circ j_{\star}\simeq i_{\star}\circ f^{\star}\colon\fromto{V}{X}.
\end{equation*}
\end{ntn}

\begin{ntn}\label{ntn:dualityfunctorD} Additionally, there is a self-duality equivalence
\begin{equation*}
\equivto{\widetilde{\mathscr{O}}^{\op}}{\widetilde{\mathscr{O}}},
\end{equation*}
whence we have a natural equivalence
\begin{equation*}
D\colon\equivto{A^{\eff,\op}}{A^{\eff}}
\end{equation*}
such that the diagram
\begin{equation*}
\begin{tikzpicture} 
\matrix(m)[matrix of math nodes, 
row sep=4ex, column sep=4ex, 
text height=1.5ex, text depth=0.25ex] 
{&A^{\eff,\op}&\\
\id&&\op\\
&A^{\eff}&\\}; 
\path[>=stealth,->,font=\scriptsize] 
(m-1-2) edge[inner sep=0.75pt] node[right]{$D$} (m-3-2) 
(m-2-1) edge node[above left]{$(\cdot)^{\star,\op}$} (m-1-2) 
edge node[below left]{$(\cdot)_{\star}$} (m-3-2)
(m-2-3) edge node[above right]{$(\cdot)_{\star}^{\op}$} (m-1-2)
edge node[below right]{$(\cdot)^{\star}$} (m-3-2); 
\end{tikzpicture}
\end{equation*}
commutes.
\end{ntn}

%----------------------------------------------------------------------%
%----------------------------------------------------------------------%

\section{Disjunctive $\infty$-categories} An $\infty$-category with all finite limits is disjunctive if the coproduct acts effectively as a disjoint union. The effective Burnside $\infty$-category $A^{\eff}(C)$ of a disjunctive $\infty$-category $C$ has the peculiar property that the initial object of $C$ becomes a zero object in $A^{\eff}(C)$, and the coproduct in $C$ becomes both the coproduct and the product in $A^{\eff}(C)$. This permits us to regard the effective Burnside $\infty$-category as somewhat ``algebraic'' in nature.

\begin{dfn} \label{item:directsums} Suppose $C$ is an $\infty$-category. Then $C$ is said to \textbf{\emph{admit direct sums}} if the following conditions hold.
\begin{enumerate}[(\ref{item:directsums}.1)]
\item The $\infty$-category $C$ is pointed.
\item The $\infty$-category $C$ has all finite products and coproducts.
\item For any finite set $I$ and any $I$-tuple $(X_i)_{i\in I}$ of objects of $C$, the map
\begin{equation*}
\fromto{\coprod X_I}{\prod X_I}
\end{equation*}
in $hC$ --- given by the maps $\phi_{ij}\colon\fromto{X_i}{X_j}$, where $\phi_{ij}$ is zero unless $i=j$, in which case it is the identity --- is an isomorphism.
\end{enumerate}
If $C$ admits finite direct sums, then for any finite set $I$ and any $I$-tuple $(X_i)_{i\in I}$ of objects of $C$, we denote by $\bigoplus X_I$ the product (or, equivalently, the coproduct) of the $X_i$.

If $C$ admits direct sums, then $C$ will be said to be \textbf{\emph{additive}} if its homotopy category $hC$ is additive. Denote by $\Cat_{\infty}^{\add}\subset\Cat_{\infty}(\kappa_1)$ the subcategory consisting of locally small additive $\infty$-categories and functors between them that preserve direct sums.
\end{dfn}

We are mostly interested in $\infty$-categories $C$ such that the $\infty$-category $A^{\eff}(C)$ admit direct sums. To ensure this, we introduce the following class of $\infty$-categories.

\begin{dfn} An $\infty$-category will be called \textbf{\emph{disjunctive}} if it admits all finite limits and finite coproducts and if, in addition, finite coproducts are disjoint and universal \cite[\S 6.1.1, (ii) and (iii)]{HTT}.

Equivalently, an $\infty$-category $C$ that admits all finite limits and all finite coproducts is disjunctive just in case, for any finite set $I$ and any collection $\{X_i\}_{i\in I}$ of objects of $C$, the natural functor
\begin{equation*}
\fromto{\prod_{i\in I}C_{/X_i}}{C_{/\coprod_{i\in I}X_i}}
\end{equation*}
given by the coproduct is an equivalence of $\infty$-categories. Its inverse is given informally by the assignment
\begin{equation*}
\goesto{U}{(U\times_{\coprod_{j\in I}X_j}X_i)_{i\in I}}
\end{equation*}

Let us denote by $\Cat_{\infty}^{\disj}\subset\Cat_{\infty}$ the subcategory whose objects are (small) disjunctive $\infty$-categories and whose morphisms are those functors that preserve pullbacks and finite coproducts.
\end{dfn}

In ordinary category theory, it may be more customary to refer to categories with the properties described above with the portmanteau ``lextensive.'' We won't be doing that.

\begin{prp}\label{prp:Aeffdirsums} If $C$ is a disjunctive $\infty$-category, then the $\infty$-category $A^{\eff}(C)$ admits direct sums.
\begin{proof} We show that the natural functor $(\cdot)_{\star}\colon\fromto{C}{A^{\eff}(C)}$ preserves coproducts. The result will then follow from the self-duality $D\colon A^{\eff}(C)^{\op}\simeq A^{\eff}(C)$. Unwinding the definitions, the claim that $(\cdot)_{\star}$ preserves coproducts amounts to the following claim: for any object $Y$ of $C$, the space
\begin{equation*}
\iota C_{/\{\varnothing,Y\}}
\end{equation*}
is contractible (which follows directly from \cite[Lm. 6.1.3.6]{HTT}), and for any objects $X$ and $X'$ of $C$, the map
\begin{equation*}
\fromto{\iota C_{/\{X\sqcup Y,Z\}}}{\iota C_{/\{X,Z\}}\times\iota C_{/\{Y,Z\}}}
\end{equation*}
given informally by the assignment
\begin{equation*}
\goesto{W}{(W\times_{X\sqcup Y}X,W\times_{X\sqcup Y}Y)}
\end{equation*}
is an equivalence. We claim that the map
\begin{equation*}
\fromto{\iota C_{/\{X,Z\}}\times\iota C_{/\{Y,Z\}}}{\iota C_{/\{X\sqcup Y,Z\}}}
\end{equation*}
given informally by the assignment
\begin{equation*}
\goesto{(U,V)}{U\sqcup V}
\end{equation*}
is a homotopy inverse. Indeed, the statement that
\begin{equation*}
W\simeq(W\times_{X\sqcup Y}X)\sqcup(W\times_{X\sqcup Y}Y)
\end{equation*}
follows from the universality of finite coproducts, and the statement that
\begin{equation*}
U\simeq(U\sqcup V)\times_{X\sqcup Y}X\textrm{\quad and\quad}V\simeq(U\sqcup V)\times_{X\sqcup Y}Y
\end{equation*}
follows from the identifications
\begin{equation*}
X\simeq X\times_{X\sqcup Y}X\textrm{,\qquad}\varnothing\simeq X\times_{X\sqcup Y}Y\textrm{\quad and\quad}Y\simeq Y\times_{X\sqcup Y}Y,
\end{equation*}
all of which follow easily from the disjointness and universality of coproducts.
\end{proof}
\end{prp}

Note that we do not quite use the full strength of disjunctivity here. It would have been enough to assume only that $C$ admits pullbacks and finite coproducts and that finite coproducts are disjoint and universal. However, we will use the further condition that $C$ admits a terminal object and hence all finite products when we study the Burnside Mackey functor.

%----------------------------------------------------------------------%

\section{Disjunctive triples} In a little while we will define Mackey functors on a disjunctive $\infty$-category $C$ (valued in spectra, say) as direct-sum preserving functors $M$ from the effective Burnside category $A^{\eff}(C)$. This means that for any object $X$ of $C$, we'll have an associated spectrum $M(X)$, and for any morphism $\fromto{X}{Y}$ of $C$, we'll have both a morphism $\fromto{M(X)}{M(Y)}$ and a morphism $\fromto{M(Y)}{M(X)}$. So a Mackey functor will splice together a covariant functor and a contravariant functor.

However, it is not always reasonable to expect both covariance and contravariance for all morphisms simultaneously. Rather, one may wish instead to specify classes of morphisms in which one has covariance and contravariance. This leads to the notion of a disjunctive triple.

\begin{nul} Recall \cite[Df. 1.11]{K1} that a \textbf{\emph{pair}} of $\infty$-categories $(C,C_{\dag})$ consists of an $\infty$-category $C$ and a subcategory $C_{\dag}$ \cite[\S 1.2.11]{HTT} that contains all the equivalences.

A \textbf{\emph{triple}} of $\infty$-categories is an $\infty$-category equipped with two pair structures. That is, a triple $(C,C_{\dag},C^{\dag})$ consists of an $\infty$-category $C$ and two subcategories
\begin{equation*}
C_{\dag},C^{\dag}\subset C,
\end{equation*}
each of which contains all the equivalences. We call morphisms of $C_{\dag}$ \textbf{\emph{ingressive}} and morphisms of $C^{\dag}$ \textbf{\emph{egressive}}.
\end{nul}

\begin{dfn}\label{dfn:adequate} A triple $(C,C_{\dag},C^{\dag})$ of $\infty$-categories is said to be \textbf{\emph{adequate}} if the following conditions obtain.
\begin{enumerate}[(\ref{dfn:adequate}.1)]
\item\label{item:ambipbsexist} For any ingressive morphism $\cofto{Y}{X}$ and any egressive morphism $\fibto{X'}{X}$, there exists a pullback square
\begin{equation*}
\begin{tikzpicture} 
\matrix(m)[matrix of math nodes, 
row sep=4ex, column sep=4ex, 
text height=1.5ex, text depth=0.25ex] 
{Y'&X'\\ 
Y&X.\\}; 
\path[>=stealth,->,font=\scriptsize] 
(m-1-1) edge (m-1-2) 
edge (m-2-1) 
(m-1-2) edge[->>] (m-2-2) 
(m-2-1) edge[>->] (m-2-2); 
\end{tikzpicture}
\end{equation*}
\item\label{item:ambipbsareambi} In any pullback square
\begin{equation*}
\begin{tikzpicture} 
\matrix(m)[matrix of math nodes, 
row sep=4ex, column sep=4ex, 
text height=1.5ex, text depth=0.25ex] 
{Y'&X'\\ 
Y&X.\\}; 
\path[>=stealth,->,font=\scriptsize] 
(m-1-1) edge node[above]{$f'$} (m-1-2) 
edge (m-2-1) 
(m-1-2) edge (m-2-2) 
(m-2-1) edge node[below]{$f$} (m-2-2); 
\end{tikzpicture}
\end{equation*}
if $f$ is ingressive (respectively, egressive), then so is $f'$.
\suspend{enumerate}
We will say that an adequate triple $(C,C_{\dag},C^{\dag})$ is a \textbf{\emph{disjunctive triple}} if the following further conditions obtain.
\resume{enumerate}[{[(\ref{dfn:adequate}.1)]}]
\item The $\infty$-category $C$ admits finite coproducts.
\item\label{item:compatiblewithcoprods} The class of ingressive morphisms and the class of egressive morphisms are each \textbf{\emph{compatible with coproducts}} in the following sense. First, any morphism from an initial object is both ingressive and egressive. Second, for any objects $X,Y,Z\in C$, a morphism $\fromto{X\sqcup Y}{Z}$ is ingressive (respectively, egressive) just in case both the restrictions $\fromto{X}{Z}$ and $\fromto{Y}{Z}$ are so.
\item Suppose $I$ and $J$ finite sets. Suppose that for any pair $(i,j)\in I\times J$, we are given a pullback square
\begin{equation*}
\begin{tikzpicture} 
\matrix(m)[matrix of math nodes, 
row sep=4ex, column sep=4ex, 
text height=1.5ex, text depth=0.25ex] 
{X'_{ij}&Y'_j\\ 
X_i&Y\\}; 
\path[>=stealth,->,font=\scriptsize] 
(m-1-1) edge[>->] (m-1-2) 
edge[->>] (m-2-1) 
(m-1-2) edge[->>] (m-2-2) 
(m-2-1) edge[>->] (m-2-2); 
\end{tikzpicture}
\end{equation*}
a pullback square in which $\cofto{X_i}{Y}$ is ingressive and $\fibto{Y'_j}{Y}$ is egressive. Then the resulting square
\begin{equation*}
\begin{tikzpicture} 
\matrix(m)[matrix of math nodes, 
row sep=6ex, column sep=4ex, 
text height=1.5ex, text depth=0.25ex] 
{\coprod_{(i,j)\in I\times J}X'_{ij}&\coprod_{j\in J}Y'_j\\ 
\coprod_{i\in I}X_i&Y\\}; 
\path[>=stealth,->,font=\scriptsize] 
(m-1-1) edge[>->] (m-1-2) 
edge[->>] (m-2-1) 
(m-1-2) edge[->>] (m-2-2) 
(m-2-1) edge[>->] (m-2-2); 
\end{tikzpicture}
\end{equation*}
is also a pullback square.
\end{enumerate}

If $(C,C_{\dag},C^{\dag})$ is a disjunctive triple, then we shall call a pullback square
\begin{equation*}
\begin{tikzpicture} 
\matrix(m)[matrix of math nodes, 
row sep=6ex, column sep=4ex, 
text height=1.5ex, text depth=0.25ex] 
{Y'&Y\\ 
X'&X\\}; 
\path[>=stealth,->,font=\scriptsize] 
(m-1-1) edge[>->] node[above]{$j$} (m-1-2) 
edge[->>] node[left]{$p'$} (m-2-1) 
(m-1-2) edge[->>] node[right]{$p$} (m-2-2) 
(m-2-1) edge[>->] node[below]{$i$} (m-2-2); 
\end{tikzpicture}
\end{equation*}
of $C$ \textbf{\emph{disjunctive}} if $i$ (and hence also $j$) is ingressive and $p$ (and hence also $p'$) is egressive.

Now a \textbf{\emph{functor of disjunctive triples}} is a functor of triples
\begin{equation*}
f\colon\fromto{(C,C_{\dag},C^{\dag})}{(D,D_{\dag},D^{\dag})}
\end{equation*}
(i.e., a functor $\fromto{C}{D}$ that carries ingressive morphisms to ingressive morphisms and egressive morphisms to egressive morphisms) such that $f$ preserves finite coproducts and ambigressive pullbacks.
\end{dfn}

\begin{ntn} Let us write $\Trip_{\infty}^{\disj}$ for the subcategory of the $\infty$-category $\Trip_{\infty}$ of small triples of $\infty$-categories whose objects are disjunctive triples and whose morphisms are functors of disjunctive triples.
\end{ntn}

\begin{exm} Of course every disjunctive $\infty$-category $C$ admits its \textbf{\emph{maximal triple structure}} $(C,C,C)$. Hence anything said of disjunctive triples specializes to disjunctive $\infty$-categories.
\end{exm}

\begin{dfn}\label{dfn:ambigressfunct} Suppose $(C,C_{\dag},C^{\dag})$ a triple of $\infty$-categories. For any integer $n\geq 0$, let us say that a functor $X\colon\fromto{\widetilde{\mathscr{O}}(\Delta^n)^{\op}}{C}$ is \textbf{\emph{ambigressive cartesian}} (relative to the triple structure) if, for any integers $0\leq i\leq k\leq\ell\leq j\leq n$, the square
\begin{equation*}
\begin{tikzpicture} 
\matrix(m)[matrix of math nodes, 
row sep=4ex, column sep=4ex, 
text height=1.5ex, text depth=0.25ex] 
{X_{ij}&X_{kj}\\ 
X_{i\ell}&X_{k\ell}\\}; 
\path[>=stealth,->,font=\scriptsize] 
(m-1-1) edge[->>] (m-1-2) 
edge[>->] (m-2-1) 
(m-1-2) edge[->>] (m-2-2) 
(m-2-1) edge[>->] (m-2-2); 
\end{tikzpicture}
\end{equation*}
is a pullback in which the morphisms $\cofto{X_{ij}}{X_{kj}}$ and $\cofto{X_{i\ell}}{X_{k\ell}}$ are ingressive, and the morphisms $\fibto{X_{ij}}{X_{i\ell}}$ and $\fibto{X_{kj}}{X_{k\ell}}$ are egressive.

Recall the functor $\RR_{\ast}(C)\colon\fromto{\Delta^{\op}}{\Kan}$ from Nt. \ref{dfn:Rstar}. Write $\AA^{\eff}_{\ast}(C,C_{\dag},C^{\dag})\subset\RR_{\ast}(C)$ for the subfunctor in which $\AA^{\eff}_n(C,C_{\dag},C^{\dag})$ is the full simplicial subset of $\RR_n(C)$ spanned by the cartesian functors $X\colon\fromto{\widetilde{\mathscr{O}}(\Delta^n)^{\op}}{C}$. Note that since any functor that is equivalent to an cartesian functor is itself cartesian, the simplicial set $\AA^{\eff}_n(C,C_{\dag},C^{\dag})$ is a union of connected components of $\RR_n(C)$.
\end{dfn}

The proof of the following is virtually identical to that of Pr. \ref{prp:QQstarisCSS}.
\begin{prp} For any adequate triple of $\infty$-categories $(C,C_{\dag},C^{\dag})$, the simplicial space $\AA^{\eff}_{\ast}(C,C_{\dag},C^{\dag})$ is a complete Segal space.
\end{prp}

\begin{dfn}\label{dfn:geneffBurn} For any adequate triple of $\infty$-categories $(C,C_{\dag},C^{\dag})$, denote by $A^{\eff}(C,C_{\dag},C^{\dag})$ the $\infty$-category whose $n$-simplices are vertices of $\AA^{\eff}_n(C,C_{\dag},C^{\dag})$, i.e., ambigressive cartesian functors $\fromto{\widetilde{\mathscr{O}}(\Delta^n)^{\op}}{C}$. We may call this the \textbf{\emph{effective Burnside $\infty$-category of $(C,C_{\dag},C^{\dag})$}}.
\end{dfn}

\begin{nul}\label{nul:Aefftripledirsums} Suppose $(C,C_{\dag},C^{\dag})$ a locally small adequate triple of $\infty$-categories $(C,C_{\dag},C^{\dag})$. Here's an alternate way to go about defining $A^{\eff}(C,C_{\dag},C^{\dag})$. Write $\mathscr{P}(C)\coloneq\Fun(C^{\op},\Kan)$ for the usual $\infty$-category of presheaves of spaces. We may describe (an equivalent version of)
\begin{equation*}
A^{\eff}(C,C_{\dag},C^{\dag})\subset A^{\eff}(\mathscr{P}(C))
\end{equation*}
as the subcategory whose objects are those functors that are representable, in which a morphism $\fromto{F}{G}$ of $A^{\eff}(\mathscr{P}(C))$ exhibited as a diagram
\begin{equation*}
\begin{tikzpicture}[baseline]
\matrix(m)[matrix of math nodes, 
row sep=3ex, column sep=3ex, 
text height=1.5ex, text depth=0.25ex] 
{&H&\\ 
F&&G,\\}; 
\path[>=stealth,->,font=\scriptsize] 
(m-1-2) edge (m-2-1) 
edge (m-2-3); 
\end{tikzpicture}
\end{equation*}
lies in $A^{\eff}(C,C_{\dag},C^{\dag})$ just in case the morphism of $C$ representing the morphism $\fromto{H}{F}$ of $\mathscr{P}(C)$ is egressive and the morphism of $C$ representing the morphism $\fromto{H}{G}$ of $\mathscr{P}(C)$ is ingressive. Since the pullback of an ingressive (respectively, egressive) morphism along an egressive (resp., ingressive) morphism is ingressive (resp., egressive) of $C$, and since the Yoneda embedding preserves all limits that exist in $C$, it follows that $A^{\eff}(C,C_{\dag},C^{\dag})$ is indeed a subcategory.

Note that for any disjunctive triple $(C,C_{\dag},C^{\dag})$, the subcategory
\begin{equation*}
A^{\eff}(C,C_{\dag},C^{\dag})\subset A^{\eff}(\mathscr{P}(C))
\end{equation*}
is closed under direct sums.
\end{nul}

\begin{ntn} Suppose $(C,C_{\dag},C^{\dag})$ an adequate triple. The two natural transformations
\begin{equation*}
(\cdot)_{\star}\colon\fromto{\id}{A^{\eff}}\textrm{\quad and\quad}(\cdot)^{\star}\colon\fromto{\op}{A^{\eff}}
\end{equation*}
restrict to yield functors
\begin{equation*}
(\cdot)_{\star}\colon\fromto{C_{\dag}}{A^{\eff}(C,C_{\dag},C^{\dag})}\textrm{\quad and\quad}(\cdot)^{\star}\colon\fromto{C^{\dag,\op}}{A^{\eff}(C,C_{\dag},C^{\dag})}
\end{equation*}
Consequently, for any ingressive morphism $f\colon\cofto{U}{V}$ of $C$, one obtains a morphism $f_{\star}\colon\fromto{U}{V}$, and for any egressive morphism $f\colon\fibto{U}{V}$ of $C$, one obtains a morphism $f^{\star}\colon\fromto{V}{U}$. Additionally, for any pullback square
\begin{equation*}
\begin{tikzpicture} 
\matrix(m)[matrix of math nodes, 
row sep=4ex, column sep=4ex, 
text height=1.5ex, text depth=0.25ex] 
{U&X\\ 
V&Y\\}; 
\path[>=stealth,->,font=\scriptsize] 
(m-1-1) edge[>->] node[above]{$i$} (m-1-2) 
edge[->>] node[left]{$f$} (m-2-1) 
(m-1-2) edge[->>] node[right]{$g$} (m-2-2) 
(m-2-1) edge[>->] node[below]{$j$} (m-2-2); 
\end{tikzpicture}
\end{equation*}
in which $i,j$ are ingressive and $f,g$ are egressive, one obtains a homotopy
\begin{equation*}
g^{\star}\circ j_{\star}\simeq i_{\star}\circ f^{\star}\colon\fromto{V}{X}.
\end{equation*}
\end{ntn}

\begin{ntn} Suppose $(C,C_{\dag},C^{\dag})$ a disjunctive triple. Then the self-duality equivalence
\begin{equation*}
\equivto{\widetilde{\mathscr{O}}^{\op}}{\widetilde{\mathscr{O}}},
\end{equation*}
induces the natural equivalence
\begin{equation*}
D\colon\equivto{A^{\eff}(C,C^{\dag},C_{\dag})^{\op}}{A^{\eff}(C,C_{\dag},C^{\dag})}
\end{equation*}
such that the diagram
\begin{equation*}
\begin{tikzpicture} 
\matrix(m)[matrix of math nodes, 
row sep=4ex, column sep=4ex, 
text height=1.5ex, text depth=0.25ex] 
{&A^{\eff}(C,C^{\dag},C_{\dag})^{\op}&\\
C_{\dag}&&C^{\dag,\op}\\
&A^{\eff}(C,C_{\dag},C^{\dag})&\\}; 
\path[>=stealth,->,font=\scriptsize] 
(m-1-2) edge[inner sep=0.75pt] node[right]{$D$} (m-3-2) 
(m-2-1) edge node[above left]{$(\cdot)^{\star,\op}$} (m-1-2) 
edge node[below left]{$(\cdot)_{\star}$} (m-3-2)
(m-2-3) edge node[above right]{$(\cdot)_{\star}^{\op}$} (m-1-2)
edge node[below right]{$(\cdot)^{\star}$} (m-3-2); 
\end{tikzpicture}
\end{equation*}
commutes.
\end{ntn}

\begin{nul} A functor of disjunctive triples $f\colon\fromto{(C,C_{\dag},C^{\dag})}{(D,D_{\dag},D^{\dag})}$ induces a functor $\fromto{A^{\eff}(C,C_{\dag},C^{\dag})}{A^{\eff}(D,D_{\dag},D^{\dag})}$ that preserves direct sums.
\end{nul}

%----------------------------------------------------------------------%

\section{Mackey functors}

\begin{dfn}\label{dfn:Mackfun} Suppose $E$ an additive $\infty$-category, and suppose $(C,C_{\dag},C^{\dag})$ a disjunctive triple. Then a \textbf{\emph{Mackey functor on $(C,C_{\dag},C^{\dag})$}} valued in $E$ is a functor
\begin{equation*}
M\colon\fromto{A^{\eff}(C,C_{\dag},C^{\dag})}{E}
\end{equation*}
that preserves direct sums. If $C$ itself is a disjunctive $\infty$-category, then a \textbf{\emph{Mackey functor on $C$}} is nothing more than a Mackey functor on the maximal triple $(C,C,C)$.
\end{dfn}

\begin{exm} When $E$ is the nerve of an ordinary additive category, a Mackey functor $M\colon\fromto{A^{\eff}(C,C_{\dag},C^{\dag})}{E}$ factors in an essentially unique fashion through the homotopy category $hA^{\eff}(C,C_{\dag},C^{\dag})$ and then its local group completion (obtained by taking the Grothendieck group of the Hom sets). Hence the notion of Mackey functor described here subsumes the one defined by Dress \cite{MR0360771}.

Note that some authors define ``ordinary'' Mackey functors as functors on the local group completion of the \emph{opposite} category $hA^{\eff}(C,C_{\dag},C^{\dag})^{\op}$. This is just a matter of convention, as the duality functor provides the equivalence $A^{\eff}(C,C_{\dag},C^{\dag})^{\op}\simeq A^{\eff}(C,C^{\dag},C_{\dag})$.
\end{exm}

\begin{ntn} Suppose $E$ an additive $\infty$-category, and suppose $(C,C_{\dag},C^{\dag})$ a small disjunctive triple. Then we denote by
\begin{equation*}
\Mack(C,C_{\dag},C^{\dag};E)\subset\Fun(A^{\eff}(C,C_{\dag},C^{\dag}),E)
\end{equation*}
the full subcategory spanned by the Mackey functors. This is covariantly functorial for additive functors in $E$. For any functor $f\colon\fromto{(C,C_{\dag},C^{\dag})}{(D,D_{\dag},D^{\dag})}$ of disjunctive triples, we have an induced functor
\begin{equation*}
f^{\star}\colon\fromto{\Mack(D,D_{\dag},D^{\dag},E)}{\Mack(C,C_{\dag},C^{\dag},E)}.
\end{equation*}
These functors fit together to yield a functor
\begin{equation*}
\Mack\colon\fromto{\Trip_{\infty}^{\disj,\op}\times\Cat_{\infty}^{\add}}{\Cat_{\infty}(\kappa_1)}.
\end{equation*}

If $C$ is a disjunctive $\infty$-category, then $\Mack(C,E)=\Mack(C,C,C;E)$.
\end{ntn}

In fact, let's see that the functor $\Mack$ is valued in $\Cat_{\infty}^{\add}$.
\begin{prp} For any disjunctive triple $(C,C_{\dag},C^{\dag})$ and for any additive $\infty$-category $E$, the $\infty$-category $\Mack(C,C_{\dag},C^{\dag};E)$ is additive.
\begin{proof} It is easy to see that $\Fun(A^{\eff}(C,C_{\dag},C^{\dag}),E)$ is additive. The full subcategory
\begin{equation*}
\Mack(C,C_{\dag},C^{\dag};E)\subset\Fun(A^{\eff}(C,C_{\dag},C^{\dag}),E)
\end{equation*}
is closed under finite direct sums by noting that the constant functor at a zero object clearly preserves direct sums, and for any Mackey functors $M$ and $N$, the functor $M\oplus N$ carries zero objects to zero objects, and
\begin{eqnarray}
(M\oplus N)(X\oplus Y)&=&M(X\oplus Y)\oplus N(X\oplus Y)\nonumber\\
&\simeq&M(X)\oplus M(Y)\oplus N(X)\oplus N(Y)\nonumber\\
&\simeq&(M\oplus N)(X)\oplus (M\oplus N)(Y)\nonumber
\end{eqnarray}
for any objects $X,Y\in A^{\eff}(C,C_{\dag},C^{\dag})$.
\end{proof}
\end{prp}

Perhaps surprisingly, Mackey functors are closed under all limits and colimits.
\begin{prp} For any disjunctive triple $(C,C_{\dag},C^{\dag})$ and for any additive $\infty$-category $E$ that admits all limits (respectively, all colimits), the full subcategory $\Mack(C,C_{\dag},C^{\dag};E)\subset\Fun(A^{\eff}(C,C_{\dag},C^{\dag}),E)$ is closed under limits (resp., under colimits).
\begin{proof} We will prove the statement about colimits. The statement about limits will then follow from consideration of the equivalence
\begin{equation*}
\Mack(C,C_{\dag},C^{\dag};E)^{\op}\simeq\Mack(C,C^{\dag},C_{\dag};E^{\op}).
\end{equation*}
We have already seen that $\Mack(C,C_{\dag},C^{\dag};E)\subset\Fun(A^{\eff}(C,C_{\dag},C^{\dag}),E)$ is closed under finite coproducts; it therefore remains to show that if $\Lambda$ is a sifted $\infty$-category and if $M\colon\fromto{\Lambda}{\Mack(C,C_{\dag},C^{\dag};E)}$ is a diagram of Mackey functors, then the colimit $M_{\infty}=\colim_{\alpha\in\Lambda}M_{\alpha}$ in $\Fun(A^{\eff}(C,C_{\dag},C^{\dag}),E)$ is again a Mackey functor. For this, suppose $X,Y\in A^{\eff}(C,C_{\dag},C^{\dag})$, and observe that
\begin{equation*}
M_{\infty}(X\oplus Y)\simeq\colim_{\alpha\in\Lambda}M_{\alpha}(X\oplus Y)\simeq\colim_{\alpha\in\Lambda}M_{\alpha}(X)\oplus M_{\alpha}(Y).
\end{equation*}
Now since $\Lambda$ is sifted, we further have
\begin{equation*}
\colim_{\alpha\in\Lambda}M_{\alpha}(X)\oplus M_{\alpha}(Y)\simeq\colim_{\alpha,\beta\in\Lambda}M_{\alpha}(X)\oplus M_{\beta}(Y).
\end{equation*}
Since
\begin{equation*}
\colim_{\alpha,\beta\in\Lambda}M_{\alpha}(X)\oplus M_{\beta}(Y)\simeq\colim_{\alpha\in\Lambda}M_{\alpha}(X)\oplus\colim_{\beta\in\Lambda}M_{\beta}(Y),
\end{equation*}
we obtain
\begin{equation*}
M_{\infty}(X\oplus Y)\simeq M_{\infty}(X)\oplus M_{\infty}(Y).\qedhere
\end{equation*}
\end{proof}
\end{prp}

\begin{cor} For any disjunctive triple $(C,C_{\dag},C^{\dag})$ and for any additive $\infty$-category $E$ that admits all limits (respectively, all colimits), limits (resp., colimits) in $\Mack(C,C_{\dag},C^{\dag};E)$ are computed objectwise.
\end{cor}

\begin{exm} For any disjunctive triple $(C,C_{\dag},C^{\dag})$, the effective Burnside $\infty$-category admits a local group completion, which is a universal target for Mackey functors. This is the \textbf{\emph{Burnside $\infty$-category}}. More precisely, there exists an additive category $A(C,C_{\dag},C^{\dag})$ and a Mackey functor
\begin{equation*}
\fromto{A^{\eff}(C,C_{\dag},C^{\dag})}{A(C,C_{\dag},C^{\dag})}
\end{equation*}
with the following universal property. For any additive $\infty$-category $E$, the functor
\begin{equation*}
\fromto{\Fun^{\add}(A(C,C_{\dag},C^{\dag}),E)}{\Mack(C,C_{\dag},C^{\dag};E)}
\end{equation*}
is an equivalence, where $\Fun^{\add}(A(C,C_{\dag},C^{\dag}),E)\subset\Fun(A(C,C_{\dag},C^{\dag}),E)$ is the full subcategory of additive functors.

This follows from the fact that the functor
\begin{equation*}
\Mack(C,C_{\dag},C^{\dag};-)\colon\fromto{\Cat_{\infty}^{\add}}{\Cat_{\infty}^{\add}}
\end{equation*}
preserves all limits, which in turn follows from the fact that limits in $\Cat_{\infty}^{\add}$ are computed in the $\infty$-category $\Cat_{\infty}$. We leave the details to the reader, as we will not need the Burnside $\infty$-category itself in this paper.
\end{exm}

Mackey functors valued in an additive $\infty$-category $E$ will inherit duality functors. To illustrate, we focus particularly on the case of Mackey functors valued in finite spectra.

\begin{nul} Suppose $(C,C_{\dag},C^{\dag})$ a disjunctive triple. Then the equivalence
\begin{equation*}
D\colon\equivto{A^{\eff}(C,C_{\dag},C^{\dag})}{A^{\eff}(C,C^{\dag},C_{\dag})^{\op}}
\end{equation*}
induces an equivalence
\begin{equation*}
\Mack(C,C_{\dag},C^{\dag};E)\simeq\Mack(C,C^{\dag},C_{\dag};E^{\op})^{\op}
\end{equation*}
for any additive $\infty$-category $E$.
\end{nul}

\begin{dfn} Suppose $(C,C_{\dag},C^{\dag})$ a disjunctive triple. The (Spanier--Whitehead) duality functor $\equivto{\Sp^{\omega}}{\Sp^{\omega,\op}}$ for finite spectra can be composed with the equivalence above to yield an equivalence
\begin{equation*}
\Mack(C,C_{\dag},C^{\dag};\Sp^{\omega})^{\op}\simeq\Mack(C,C_{\dag},C^{\dag};\Sp^{\omega,\op})^{\op}\simeq\Mack(C,C^{\dag},C_{\dag};\Sp^{\omega}).
\end{equation*}
The image of a Mackey functor
\begin{equation*}
M\colon\fromto{A^{\eff}(C,C_{\dag},C^{\dag})}{\Sp^{\omega}}
\end{equation*}
under this equivalence is the \textbf{\emph{dual Mackey functor}}
\begin{equation*}
M^{\vee}\colon\fromto{A^{\eff}(C,C^{\dag},C_{\dag})}{\Sp^{\omega}}.
\end{equation*}
\end{dfn}

We conclude this subsection by remarking on the potential \emph{covariant} dependence of $\Mack(C,C_{\dag},C^{\dag};E)$ on $C$.

\begin{ntn} Suppose $E$ a presentable additive $\infty$-category, suppose $(C,C_{\dag},C^{\dag})$ and $(D,D_{\dag},D^{\dag})$ disjunctive $\infty$-categories, and suppose
\begin{equation*}
f\colon\fromto{(C,C_{\dag},C^{\dag})}{(D,D_{\dag},D^{\dag})}
\end{equation*}
a functor of disjunctive triples. Then the induced functor
\begin{equation*}
f^{\star}\colon\fromto{\Mack(D,D_{\dag},D^{\dag};E)}{\Mack(C,C_{\dag},C^{\dag};E)}
\end{equation*}
preserves limits (since they are computed objectwise), whence it follows from the Adjoint Functor Theorem \cite{HTT} that it admits a left adjoint $f_!$. These left adjoints fit together to yield a functor
\begin{equation*}
\Mack\colon\fromto{\Trip_{\infty}^{\disj}\times\Pr^{\LL,\add}}{\Pr^{\LL,\add}},
\end{equation*}
where, in the notation of \cite[Df. 5.5.3.1]{HTT},
\begin{equation*}
\Pr^{\LL,\add}\coloneq\Pr^{\LL}\cap\Cat_{\infty}^{\add},
\end{equation*}
the $\infty$-category of presentable additive $\infty$-categories, whose morphisms are left adjoints.
\end{ntn}

%----------------------------------------------------------------------%

\section{Mackey stabilization}

\begin{dfn}\label{dfn:Mackstab} Suppose $(C,C_{\dag},C^{\dag})$ a disjunctive triple, suppose $E$ a presentable $\infty$-category in which filtered colimits are left exact \cite[Df. 7.3.4.2]{HTT}, and suppose
\begin{equation*}
f\colon\fromto{A^{\eff}(C,C_{\dag},C^{\dag})}{E}\textrm{\quad and\quad}F\colon\fromto{A^{\eff}(C,C_{\dag},C^{\dag})}{\Sp(E)}
\end{equation*}
two functors. Then a natural transformation
\begin{equation*}
\eta\colon\fromto{f}{\Omega^{\infty}\circ F}
\end{equation*}
will be said to exhibit $F$ as the \textbf{\emph{Mackey stabilization}} of $f$ if $F$ is a Mackey functor, and if, for any Mackey functor $M\colon\fromto{A^{\eff}(C,C_{\dag},C^{\dag})}{\Sp(E)}$, the map
\begin{equation*}
\fromto{\Map_{\Mack(C,C_{\dag},C^{\dag};\Sp(E))}(F,M)}{\Map_{\Fun(A^{\eff}(C,C_{\dag},C^{\dag}),E)}(f,\Omega^{\infty}\circ M)}
\end{equation*}
induced by $\eta$ is an equivalence.
\end{dfn}

We shall now show that Mackey stabilization  exist and are computable by means of a Goodwillie derivative. In particular, we will show that the functor
\begin{equation*}
\Omega^{\infty}\circ-\colon\fromto{\Mack(C,C_{\dag},C^{\dag};\Sp(E))}{\Fun(A^{\eff}(C,C_{\dag},C^{\dag}),E)}
\end{equation*}
admits a left adjoint.

\begin{ntn}\label{ntn:DAC} Suppose $(C,C_{\dag},C^{\dag})$ a disjunctive triple. Then we write
\begin{equation*}
\mathrm{D}A(C,C_{\dag},C^{\dag})\coloneq\mathscr{P}_{\Sigma}(A^{\eff}(C,C_{\dag},C^{\dag}))
\end{equation*}
for the nonabelian derived $\infty$-category of $A^{\eff}(C,C_{\dag},C^{\dag})$ \cite{HTT}. This $\infty$-category admits all colimits, and it comes equipped with a fully faithful functor
\begin{equation*}
j\colon\into{A^{\eff}(C,C_{\dag},C^{\dag})}{\mathrm{D}A(C,C_{\dag},C^{\dag})}.
\end{equation*}
These are uniquely characterized by \emph{either} of the following conditions.
\begin{enumerate}[(\ref{ntn:DAC}.1)]
\item For any $\infty$-category $E$ that admits all colimits, the restriction functor
\begin{equation*}
\fromto{\Fun^{\mathrm{L}}(\mathrm{D}A(C,C_{\dag},C^{\dag}),E)}{\Fun(A^{\eff}(C,C_{\dag},C^{\dag}),E)}
\end{equation*}
induced by $j$ is fully faithful, and its essential image is spanned by those functors $\fromto{A^{\eff}(C,C_{\dag},C^{\dag})}{E}$ that preserve finite coproducts.
\item For any $\infty$-category $E$ with all sifted colimits, the restriction functor
\begin{equation*}
\fromto{\Fun_{\mathscr{G}}(\mathrm{D}A(C,C_{\dag},C^{\dag}),E)}{\Fun(A^{\eff}(C,C_{\dag},C^{\dag}),E)}
\end{equation*}
induced by $j$ is an equivalence, where $\mathscr{G}$ denotes the class of small sifted simplicial sets.
\end{enumerate}
\end{ntn}

\noindent The following is now immediate from \cite{HA}.

\begin{lem}\label{lem:mackisexc} Suppose $(C,C_{\dag},C^{\dag})$ a disjunctive triple, and suppose $E$ a presentable $\infty$-category in which filtered colimits are left exact. Then the inclusion functor
\begin{equation*}
\into{A^{\eff}(C,C_{\dag},C^{\dag})}{\mathrm{D}A(C,C_{\dag},C^{\dag})}
\end{equation*}
and the $0$-th space functor $\Omega^{\infty}\colon\fromto{\Sp(E)}{E}$ induce equivalences
\begin{equation*}
\begin{tikzpicture} 
\matrix(m)[matrix of math nodes, 
row sep=4ex, column sep=4ex, 
text height=1.5ex, text depth=0.25ex] 
{\Fun^{\mathrm{L}}(\mathrm{D}A(C,C_{\dag},C^{\dag}),\Sp(E))&\Mack(C,C_{\dag},C^{\dag};\Sp(E))\\ 
\Exc_{\mathscr{F}}(\mathrm{D}A(C,C_{\dag},C^{\dag}),\Sp(E))&\Exc_{\mathscr{F}}(\mathrm{D}A(C,C_{\dag},C^{\dag}),E)\\}; 
\path[>=stealth,->,,inner sep=0.5pt,font=\scriptsize] 
(m-1-1) edge node[above]{$\sim$} (m-1-2) 
edge[-,double distance=1.5pt] (m-2-1) 
(m-2-1) edge node[below]{$\sim$} (m-2-2); 
\end{tikzpicture}
\end{equation*}
where $\Exc_{\mathscr{F}}$ denotes the $\infty$-category of ($1$-)excisive functor that preserve all filtered colimits.
\end{lem}

Now we are well positioned to obtain Mackey stabilizations.

\begin{prp}\label{prp:Mackstabexist} Suppose $(C,C_{\dag},C^{\dag})$ a disjunctive triple, and suppose $E$ a presentable $\infty$-category in which filtered colimits are left exact. Then any functor $f\colon\fromto{A^{\eff}(C,C_{\dag},C^{\dag})}{E}$ admits a Mackey stabilization. In particular, the functor
\begin{equation*}
\Omega^{\infty}\circ-\colon\fromto{\Mack(C,C_{\dag},C^{\dag};\Sp(E))}{\Fun(A^{\eff}(C,C_{\dag},C^{\dag}),E)}
\end{equation*}
admits a left adjoint.
\begin{proof} Compose the equivalences of the previous lemma with the $1$-excisive approximation functor $P_1\colon\fromto{\Fun_{\mathscr{F}}(\mathrm{D}A(C,C_{\dag},C^{\dag}),E)}{\Exc_{\mathscr{F}}(\mathrm{D}A(C,C_{\dag},C^{\dag}),E)}$, which is left adjoint to the inclusion. Employing the equivalences of the previous lemma, a left adjoint to the inclusion functor
\begin{equation*}
\into{\Mack(C,C_{\dag},C^{\dag};\Sp(E))\simeq\Exc_{\mathscr{F}}(\mathrm{D}A(C,C_{\dag},C^{\dag}),E)}{\Fun_{\mathscr{F}}(\mathrm{D}A(C,C_{\dag},C^{\dag}),E)}.
\end{equation*}
This left adjoint may now be composed with the inclusion
\begin{equation*}
\into{\Fun(A^{\eff}(C,C_{\dag},C^{\dag}),E)\simeq\Fun_{\mathscr{G}}(\mathrm{D}A(C,C_{\dag},C^{\dag}),E)}{\Fun_{\mathscr{F}}(\mathrm{D}A(C,C_{\dag},C^{\dag}),E)}
\end{equation*}
to obtain the desired Mackey stabilization functor.
\end{proof}
\end{prp}

Happily, Tom Goodwillie has provided us with a formula for the $1$-excisive approximation \cite[Cnstr. 7.1.1.27]{HA}. Hence for any functor $f\colon\fromto{A^{\eff}(C,C_{\dag},C^{\dag})}{E}$, we obtain a formula for its Mackey stabilization as a $1$-excisive functor $\Omega^{\infty}\circ F\colon\fromto{\mathrm{D}A(C,C_{\dag},C^{\dag})}{E}$:
\begin{equation*}
\Omega^{\infty}\circ F\simeq\underset{n\geq0}{\colim}\ \Omega^n\circ\overline{f}\circ\Sigma_{\mathrm{D}A(C,C_{\dag},C^{\dag})}^n,
\end{equation*}
where $\overline{f}\colon\fromto{\mathrm{D}A(C,C_{\dag},C^{\dag})}{E}$ is the essentially unique functor that preserves sifted colimits such that
\begin{equation*}
\overline{f}|_{A^{\eff}(C,C_{\dag},C^{\dag})}=f.
\end{equation*}
Consequently, we have the task of studying the suspension functor $\Sigma_{\mathrm{D}A(C,C_{\dag},C^{\dag})}$ on $\mathrm{D}A(C,C_{\dag},C^{\dag})$. In particular, we are interested in its values on objects of the effective Burnside $\infty$-category. For any object $X\in A^{\eff}(C,C_{\dag},C^{\dag})$, we have a simplicial object $B_{\ast}(0,X,0)$ given by
\begin{equation*}
\goesto{\mathbf{n}}{X^{n}},
\end{equation*}
whose geometric realization in $\mathrm{D}A(C,C_{\dag},C^{\dag})$ is $\Sigma_{\mathrm{D}A(C,C_{\dag},C^{\dag})}X$. Since $\overline{f}$ preserves geometric realizations, we find
\begin{equation*}
\Omega^{\infty}F(X)\simeq\underset{n\geq0}{\colim}\ \Omega^n\left(\underset{(\mathbf{k}_1,\dots,\mathbf{k}_n)\in(N\Delta^{\op})^n}{\colim}\ f\left(X^{k_1+\cdots+k_n}\right)\right).
\end{equation*}

In one important class of cases, it follows immediately from Segal's delooping machine that passage to the colimit is unnecessary.
\begin{prp} Suppose $(C,C_{\dag},C^{\dag})$ a disjunctive triple, and suppose $E$ an $\infty$-topos. If $f\colon\fromto{A^{\eff}(C,C_{\dag},C^{\dag})}{E}$ is a functor that preserves finite products, the Mackey stabilization $F$ of $f$ is defined by the formula
\begin{equation*}
\Omega^{\infty}F(X)\simeq\Omega\big|B_{\ast}(\ast,f(X),\ast)\big|_{N\Delta^{\op}},
\end{equation*}
where $B_{\ast}(\ast,f(X),\ast)$ is the simplicial object $\goesto{\mathbf{k}}{f(X)^k}$, and $|\cdot|_{N\Delta^{\op}}$ denotes geometric realization.
\end{prp}

%----------------------------------------------------------------------%

\section{Representable spectral Mackey functors and assembly morphisms} The Mackey stabilization is useful for constructing universal examples of Mackey functors.

\begin{dfn}\label{dfn:BurnsideMack} Suppose $(C,C_{\dag},C^{\dag})$ a disjunctive triple, and suppose $X$ and object of $C$. Then the Mackey stabilization of the functor $\fromto{A^{\eff}(C,C_{\dag},C^{\dag})}{\Kan}$ corepresented by $X$ will be denoted
\begin{equation*}
\SS_{(C,C_{\dag},C^{\dag})}^X\colon\fromto{A^{\eff}(C,C_{\dag},C^{\dag})}{\Sp}.
\end{equation*}
(We will drop the subscript and write $\SS^X$ when the chosen disjunctive triple is clear from the context.) This is the \textbf{\emph{Mackey functor corepresented by $X$}}

We will call the Mackey functor corepresented by the terminal object $1$ the \textbf{\emph{Burnside Mackey functor}}. In this case, we drop the superscript and write simply $\SS_{(C,C_{\dag},C^{\dag})}$.
\end{dfn}

The following is now an immediate consequence of the universal property of the Mackey stabilization.
\begin{prp} Suppose $(C,C_{\dag},C^{\dag})$ a disjunctive triple, and suppose $X$ and object of $C$. Then the corepresentable Mackey functor has the universal property that for any Mackey functor $M\colon\fromto{A^{\eff}(C,C_{\dag},C^{\dag})}{\Sp}$, there is an identification
\begin{equation*}
\Map_{\Mack(C,C_{\dag},C^{\dag};\Sp)}(\SS^X,M)\simeq\Omega^{\infty}M(X),
\end{equation*}
functorial in $M$.
\end{prp}

\begin{dfn}\label{dfn:assembly} Suppose $(C,C_{\dag},C^{\dag})$ a disjunctive triple, suppose $X$ and object of $C$, and suppose $M\colon\fromto{A^{\eff}(C,C_{\dag},C^{\dag})}{\Sp}$ a Mackey functor. Then the identity functor $\fromto{M(X)}{M(X)}$ defines a morphism of Mackey functors
\begin{equation*}
\fromto{\SS^X}{F(M(X),M)},
\end{equation*}
where the target is the composite
\begin{equation*}
A^{\eff}(C,C_{\dag},C^{\dag})\ \tikz[baseline]\draw[>=stealth,->,font=\scriptsize](0,0.5ex)--node[above]{$M$}(0.5,0.5ex);\ \Sp\ \tikz[baseline]\draw[>=stealth,->,font=\scriptsize](0,0.5ex)--node[above]{$F(M(X),-)$}(1.5,0.5ex);\ \Sp.
\end{equation*}
For any object $Y\in C$, we call the corresponding morphism
\begin{equation*}
\alpha\colon\fromto{\SS^X(Y)\wedge M(X)}{M(Y)}
\end{equation*}
the \textbf{\emph{assembly morphism}} for $M$, $X$, and $Y$.
\end{dfn}

%----------------------------------------------------------------------%

\section{Mackey stabilization via algebraic $K$-theory}\label{sect:MackstabisK} Let us discuss a key circumstance in which we can express the Mackey stabilization (Df. \ref{dfn:Mackstab}) of a functor in terms of the additivization presented in \cite[Df. 7.9]{K1}: we are interested in the situation in which a functor is given by composing a Mackey functor valued in Waldhausen $\infty$-categories with a suitable theory in the sense of \cite[Df. 7.1]{K1}.
\begin{dfn} Suppose $\phi\colon\fromto{\Wald_{\infty}}{\mathscr{E}}$ a pre-additive theory \cite[Df. 7.11]{K1}. Then we will say that a Waldhausen $\infty$-category $\mathscr{C}$ is \textbf{\emph{$\phi$-split}} if, for any integer $m\geq 0$, the functors $\fromto{\mathscr{F}_m(\mathscr{C})}{\mathscr{C}}$ and $\fromto{\mathscr{F}_m(\mathscr{C})}{\mathscr{S}_m(\mathscr{C})}$ induce an equivalence
\begin{equation*}
\equivto{\phi(\mathscr{F}_m(\mathscr{C}))}{\phi(\mathscr{C})\times\phi(\mathscr{S}_m(\mathscr{C}))}.
\end{equation*}
\end{dfn}

\begin{prp} Suppose $(C,C_{\dag},C^{\dag})$ a disjunctive triple, and suppose
\begin{equation*}
\mathscr{X}\colon\fromto{A^{\eff}(C,C_{\dag},C^{\dag})}{\Wald_{\infty}}
\end{equation*}
a Mackey functor valued in the $\infty$-category of Waldhausen $\infty$-categories. Suppose, additionally, that $\mathscr{E}$ is an $\infty$-topos and $\phi\colon\fromto{\Wald_{\infty}}{\mathscr{E}}$ a pre-additive theory such that for any object $s\in C$, the Waldhausen $\infty$-category is $\phi$-split. Then the Mackey stabilization of the composite $\phi\circ\mathscr{X}$ is given by
\begin{equation*}
\goesto{s}{\DD\phi(\mathscr{X}(s))},
\end{equation*}
where $\DD\phi\colon\fromto{\Wald_{\infty}}{\Sp(\mathscr{E})}$ denote the canonical lift of the additivization of \cite[Cor. 7.6.1]{K1}.
\begin{proof} Let us extend $\phi\circ\mathscr{X}$ to a functor $\fromto{\mathrm{D}A(C,C_{\dag},C^{\dag})}{\mathscr{E}}$ and compute the $1$-excisive approximation. The Mackey stabilization of $\phi\circ\mathscr{X}$ is then the spectrum-valued lift of the $1$-excisive approximation of the composite
\begin{equation*}
\mathrm{D}A(C,C_{\dag},C^{\dag})\ \tikz[baseline]\draw[>=stealth,->,font=\scriptsize](0,0.5ex)--node[above]{$\LL\mathscr{X}$}(0.75,0.5ex);\ \mathrm{D}(\Wald_{\infty})\ \tikz[baseline]\draw[>=stealth,->,font=\scriptsize](0,0.5ex)--node[above]{$\Phi$}(0.5,0.5ex);\ \mathscr{E},
\end{equation*}
where $\Phi$ is the left derived functor of $\phi$ \cite[Df. 4.14]{K1}, and $\LL\mathscr{X}$ is the essentially unique colimit-preserving functor whose restriction to $A^{\eff}(C,C_{\dag},C^{\dag})$ is $\mathscr{X}$. Then one has (by, for example, \cite[Rk. 7.1.1.30]{HA})
\begin{equation*}
P_1(\Phi\circ\LL\mathscr{X})\simeq P_1(\Phi)\circ\LL\mathscr{X};
\end{equation*}
hence the Mackey stabilization of $\phi\circ\mathscr{X}$ is given by the spectrum-valued lift of the functor
\begin{equation*}
\goesto{s}{\colim_m\Omega^m\Phi\Sigma^m\mathscr{X}(s)}.
\end{equation*}

Here $\Sigma$ denotes the suspension in $\mathrm{D}(\Wald_{\infty})$. It can be computed by means of a bar construction:
\begin{equation*}
\Sigma\mathscr{Y}\simeq|B_{\ast}(0,\mathscr{Y},0)|_{N\Delta^{\op}},
\end{equation*}
where $B_n(0,\mathscr{Y},0)\simeq\mathscr{Y}^n$. Consequently, one may compute $\Phi\circ\Sigma$ as a geometric realization:
\begin{equation*}
\Phi\Sigma(\mathscr{Y})\simeq|B_{\ast}(\ast,\Phi(\mathscr{Y}),\ast)|_{N\Delta^{\op}}.
\end{equation*}

Now let us assume further that for any object $s\in S$, the Waldhausen $\infty$-category $\mathscr{X}(s)$ is $\phi$-split.
In this case, one has
\begin{equation*}
|B_{\ast}(\ast,\phi(\mathscr{X}(s)),\ast)|_{N\Delta^{\op}}\simeq\Phi(\mathscr{S}\mathscr{X}(s)),
\end{equation*}
whence $\Omega\Phi\Sigma(\mathscr{X}(s))$ is the additivization \cite[Df. 7.9]{K1} of $\phi$ applied to $\mathscr{X}(s)$. It follows that the colimit stabilizes, and the result follows.
\end{proof}
\end{prp}

\begin{cor} Suppose $(C,C_{\dag},C^{\dag})$ a disjunctive triple, and suppose
\begin{equation*}
\mathscr{X}\colon\fromto{A^{\eff}(C,C_{\dag},C^{\dag})}{\Wald_{\infty}}
\end{equation*}
a Mackey functor valued in the $\infty$-category of Waldhausen $\infty$-categories. If each Waldhausen $\infty$-category $\mathscr{X}(s)$ is $\iota$-split, then the Mackey stabilization of $\iota\circ\mathscr{X}$ is given by $\KK\circ\mathscr{X}$, where $\KK$ denotes connective algebraic $K$-theory.
\end{cor}

%----------------------------------------------------------------------%

\section{Waldhausen bicartesian fibrations} We have already seen that there is a close relationship between algebraic $K$-theory and spectral Mackey functors. The inputs required there were Mackey functors valued in Waldhausen $\infty$-categories. Unfortunately, in nature, these Mackey functors tend not to appear with all of their coherences splayed out. Instead, the most interesting examples are found \emph{furled} --- as fibrations that exhibit both covariant functoriality and contravariant functoriality along with a compatibility between the two in certain situations. We call these fibrations \emph{Waldhausen bicartesian fibrations}. In this section, we define this notion, and in the next, we show how to \emph{unfurl} these fibrations to extract the desired Mackey functors valued in Waldhausen $\infty$-categories.

\begin{nul}\label{nul:basechange} Suppose $p\colon\fromto{X}{S}$ a cartesian and cocartesian fibration. Then for any morphism $f\colon\fromto{s}{t}$ of $S$, one has and adjoint pair of functors
\begin{equation*}
\adjunct{f_{!}}{X_s}{X_t}{f^{\star}}.
\end{equation*}
For any square
\begin{equation}\label{eqn:ambipb}
\begin{tikzpicture}[baseline]
\matrix(m)[matrix of math nodes, 
row sep=4ex, column sep=4ex, 
text height=1.5ex, text depth=0.25ex] 
{s&s'\\ 
t&t',\\}; 
\path[>=stealth,->,font=\scriptsize] 
(m-1-1) edge node[above]{$i$} (m-1-2) 
edge node[left]{$q$} (m-2-1) 
(m-1-2) edge node[right]{$q'$} (m-2-2) 
(m-2-1) edge node[below]{$j$} (m-2-2); 
\end{tikzpicture}
\end{equation}
the unit $\eta\colon\fromto{\id}{j^{\star}\circ j_!}$ induces a natural transformation
\begin{equation*}
q^{\star}\ \tikz[baseline]\draw[>=stealth,->,font=\scriptsize](0,0.5ex)--node[above]{$q^{\star}(\eta)$}(0.5,0.5ex);\ q^{\star}\circ j^{\star}\circ j_!\simeq i^{\star}\circ {q'}^{\star}\circ j_!,
\end{equation*}
which is adjoint to a natural transformation
\begin{equation}\label{eqn:basechange}
\fromto{i_!\circ q^{\star}}{{q'}^{\star}\circ j_!},
\end{equation}
which we call the \textbf{\emph{base change natural transformation}}. Equivalently, we may construct this natural transformation by using the counit $\epsilon\colon\fromto{q_!\circ q^{\star}}{\id}$ to define
\begin{equation*}
q'_{!}\circ i_{!}\circ q^{\star}\simeq j_{!}\circ f_{!}\circ q^{\star}\ \tikz[baseline]\draw[>=stealth,->,font=\scriptsize](0,0.5ex)--node[above]{$j_{!}(\epsilon)$}(0.5,0.5ex);\ j_!;
\end{equation*}
its right adjoint is then the base change natural transformation \eqref{eqn:basechange}.

When the base change natural transformation is an equivalence, then one says that the square
\begin{equation*}
\begin{tikzpicture}[baseline]
\matrix(m)[matrix of math nodes, 
row sep=4ex, column sep=4ex, 
text height=1.5ex, text depth=0.25ex] 
{X_s&X_{s'}\\ 
X_t&X_{t'},\\}; 
\path[>=stealth,->,font=\scriptsize] 
(m-1-1) edge node[above]{$i_{!}$} (m-1-2) 
edge node[left]{$q_{!}$} (m-2-1) 
(m-1-2) edge node[right]{$q'_{!}$} (m-2-2) 
(m-2-1) edge node[below]{$j_{!}$} (m-2-2); 
\end{tikzpicture}
\textrm{\qquad(respectively, the square\quad}
\begin{tikzpicture}[baseline]
\matrix(m)[matrix of math nodes, 
row sep=4ex, column sep=4ex, 
text height=1.5ex, text depth=0.25ex] 
{X_{t'}&X_{s'}\\ 
X_t&X_{s},\\}; 
\path[>=stealth,->,font=\scriptsize] 
(m-1-1) edge node[above]{${q'}^{\star}$} (m-1-2) 
edge node[left]{$j^{\star}$} (m-2-1) 
(m-1-2) edge node[right]{$i^{\star}$} (m-2-2) 
(m-2-1) edge node[below]{$q^{\star}$} (m-2-2); 
\end{tikzpicture}
\textrm{\quad)}
\end{equation*}
is \textbf{\emph{right adjointable}} (resp., \textbf{\emph{left adjointable}}) \cite{HTT}. Apparently this is sometimes called the \emph{Beck--Chevalley condition}.

If we only assume $p$ an inner fibration, we can make sense of the base change natural transformation \eqref{eqn:basechange} in the presence of a small amount of extra information. Of course one needs to know that the functors $i_{!}$, $j_{!}$, $q^{\star}$, and $q'^{\star}$ all exist, and in order to construct \eqref{eqn:basechange}, it is enough to assume \emph{either} that the functors $q_{!}$ and $q'_{!}$ exist \emph{or} that the functors $i^{\star}$ and $j^{\star}$. That is, if $\sigma\colon\fromto{\Delta^1\times\Delta^1}{S}$ is given by the square \eqref{eqn:ambipb}, then it suffices to assume only \emph{one} of the following.
\begin{enumerate}[(\ref{nul:basechange}.1)]\addtocounter{enumi}{2}
\item The functor
\begin{equation*}
\fromto{X\times_{S}(\Delta^{1}\times\Delta^{1})}{\Delta^{1}\times\Delta^{1}}
\end{equation*}
is a cocartesian fibration, and the functor
\begin{equation*}
\fromto{X\times_{S}(\Delta^{1}\times\partial\Delta^{1})}{\Delta^{1}\times\partial\Delta^{1}}
\end{equation*}
is a cartesian fibration.
\item The functor
\begin{equation*}
\fromto{X\times_{S}(\Delta^{1}\times\Delta^{1})}{\Delta^{1}\times\Delta^{1}}
\end{equation*}
is a cartesian fibration, and the functor
\begin{equation*}
\fromto{X\times_{S}(\partial\Delta^{1}\times\Delta^{1})}{\partial\Delta^{1}\times\Delta^{1}}
\end{equation*}
is a cocartesian fibration.
\end{enumerate}

We will apply this idea in the case that \eqref{eqn:ambipb} is a square of a disjunctive triple in which the vertical morphisms $q,q'$ of  are egressive and the horizontal morphisms $i,j$ are ingressive.
\end{nul}

\begin{dfn}\label{dfn:leftrightcomplete} A triple $(C,C_{\dag},C^{\dag})$ is said to be \textbf{\emph{left complete}} if $C_{\dag}\subset C^{\dag}$ and \textbf{\emph{right complete}} if $C^{\dag}\subset C_{\dag}$.
\end{dfn}

In the examples of disjunctive triples of greatest interest to us, it is often the case that either every map is egressive or every map is ingressive. In particular, most of the examples of interest to us are either left complete or right complete.

Our equivariant $K$-theory will take as input assignments of Waldhausen $\infty$-categories to objects of disjunctive triples $(C,C_{\dag},C^{\dag})$ that are covariant in ingressive morphisms and contravariant in egressive morphisms. We will insist that for morphisms that are both ingressive and egressive, the resulting functors are adjoint. Finally, we will assume that for pullback squares of egressive morphisms along ingressive morphisms, the base change natural transformation is an equivalence. But for the base change natural transformation to make sense, we must assume that $(C,C_{\dag},C^{\dag})$ is either left or right complete. This leads us to the following.

\begin{dfn}\label{dfn:Waldbicart} Suppose $(C,C_{\dag},C^{\dag})$ an adequate triple that is either left complete or right complete. An inner fibration $p\colon\fromto{X}{C}$ is said to be \textbf{\emph{adequate over the triple $(C,C_{\dag},C^{\dag})$}} if the following conditions obtain.
\begin{enumerate}[(\ref{dfn:Waldbicart}.1)]
\item\label{item:adequatecofib} For any ingressive morphism $f\colon\cofto{s}{t}$ and any object $x\in X_s$, there exists a $p$-cocartesian edge $\fromto{x}{y}$ covering $f$. In particular, the functor
\begin{equation*}
p_{\dag}\colon\fromto{X\times_CC_{\dag}}{C_{\dag}}
\end{equation*}
is a cocartesian fibration.
\item\label{item:adequatefib} For any egressive morphism $f\colon\fibto{s}{t}$ and any object $y\in X_t$, there exists a $p$-cartesian edge $\fromto{x}{y}$ covering $f$. In particular, the functor
\begin{equation*}
p^{\dag}\colon\fromto{X\times_CC^{\dag}}{C^{\dag}}
\end{equation*}
is a cartesian fibration.
\item\label{item:ambiBC} For any ambigressive pullback square
\begin{equation*}
\begin{tikzpicture} 
\matrix(m)[matrix of math nodes, 
row sep=4ex, column sep=4ex, 
text height=1.5ex, text depth=0.25ex] 
{s&s'\\ 
t&t',\\}; 
\path[>=stealth,->,font=\scriptsize] 
(m-1-1) edge[>->] node[above]{$i$} (m-1-2) 
edge[->>] node[left]{$q$} (m-2-1) 
(m-1-2) edge[->>] node[right]{$q'$} (m-2-2) 
(m-2-1) edge[>->] node[below]{$j$} (m-2-2); 
\end{tikzpicture}
\end{equation*}
the base change natural transformation
\begin{equation*}
\fromto{i_!\circ q^{\star}}{{q'}^{\star}\circ j_!}
\end{equation*}
is an equivalence.
\suspend{enumerate}

Suppose now that $(C,C_{\dag},C^{\dag})$ is a disjunctive triple that is either left complete or right complete. A \textbf{\emph{Waldhausen bicartesian fibration over the triple $(C,C_{\dag},C^{\dag})$}}
\begin{equation*}
q\colon\fromto{\mathscr{X}}{C}
\end{equation*}
is a functor of pairs $\fromto{\mathscr{X}}{C^{\flat}}$ that enjoys the following properties.
\resume{enumerate}[{[(\ref{dfn:Waldbicart}.1)]}]
\item The underlying functor $q\colon\fromto{\mathscr{X}}{C}$ is an adequate inner fibration over the triple $(C,C_{\dag},C^{\dag})$.
\item For any ingressive morphism $\eta\colon\cofto{s}{t}$, the induced functor $\eta_{!}\colon\fromto{\mathscr{X}_s}{\mathscr{X}_t}$ carries cofibrations to cofibrations, and it is an exact functor
\begin{equation*}
\fromto{(\mathscr{X}_s,\mathscr{X}_s\times_{\mathscr{X}}\mathscr{X}_{\dag})}{(\mathscr{X}_t,\mathscr{X}_t\times_{\mathscr{X}}\mathscr{X}_{\dag})}
\end{equation*}
of Waldhausen $\infty$-categories.
\item For any egressive morphism $\eta\colon\fibto{s}{t}$, the induced functor $\eta^{\star}\colon\fromto{\mathscr{X}_t}{\mathscr{X}_s}$ carries cofibrations to cofibrations, and it is an exact functor
\begin{equation*}
\fromto{(\mathscr{X}_t,\mathscr{X}_t\times_{\mathscr{X}}\mathscr{X}_{\dag})}{(\mathscr{X}_s,\mathscr{X}_s\times_{\mathscr{X}}\mathscr{X}_{\dag})}
\end{equation*}
of Waldhausen $\infty$-categories.
\item For any finite set $I$ and any collection $\{s_i\ |\ i\in I\}$ of objects of $C$ indexed by the elements of $I$ with coproduct $s$, the functors
\begin{equation*}
j_i^{\star}\colon\fromto{\mathscr{X}_s}{\mathscr{X}_{s_i}}
\end{equation*}
induced by the inclusions $j_i\colon\into{s_i}{s}$ together exhibit $\mathscr{X}_s$ as the direct sum $\bigoplus_{i\in I}\mathscr{X}_{s_i}$.
\end{enumerate}
\end{dfn}

The following lemma is just an unwinding of the relevant definitions, but it will come in handy later.

\begin{lem}\label{lem:equivambiBC} Suppose $(C,C_{\dag},C^{\dag})$ an adequate triple that is either left or right complete. If $p\colon\fromto{X}{C}$ is an inner fibration satisfying conditions  (\ref{dfn:Waldbicart}.\ref{item:adequatecofib}) and  (\ref{dfn:Waldbicart}.\ref{item:adequatefib}), then condition (\ref{dfn:Waldbicart}.\ref{item:ambiBC}) is equivalent to the condition that for any square
\begin{equation*}
\begin{tikzpicture} 
\matrix(m)[matrix of math nodes, 
row sep=4ex, column sep=4ex, 
text height=1.5ex, text depth=0.25ex] 
{x&x'\\ 
y&y'\\}; 
\path[>=stealth,->,font=\scriptsize] 
(m-1-1) edge[>->] node[above]{$\eta$} (m-1-2) 
edge[->>] node[left]{$\alpha$} (m-2-1) 
(m-1-2) edge[->>] node[right]{$\alpha'$} (m-2-2) 
(m-2-1) edge[>->] node[below]{$\theta$} (m-2-2); 
\end{tikzpicture}
\end{equation*}
of $\mathscr{X}$ that covers an ambigressive pullback square
\begin{equation*}
\begin{tikzpicture} 
\matrix(m)[matrix of math nodes, 
row sep=4ex, column sep=4ex, 
text height=1.5ex, text depth=0.25ex] 
{s&s'\\ 
t&t'\\}; 
\path[>=stealth,->,font=\scriptsize] 
(m-1-1) edge[>->] node[above]{$i$} (m-1-2) 
edge[->>] node[left]{$q$} (m-2-1) 
(m-1-2) edge[->>] node[right]{$q'$} (m-2-2) 
(m-2-1) edge[>->] node[below]{$j$} (m-2-2); 
\end{tikzpicture}
\end{equation*}
of $C$, if $\alpha$ and $\alpha'$ are $p$-cartesian and $\theta$ is $p$-cocartesian, then $\eta$ is also $p$-cocartesian. 
\end{lem}

%----------------------------------------------------------------------%

\section{Unfurling} Here is the central construction of this paper. A Waldhausen bicartesian fibration over a left or right complete disjunctive triple $(C,C_{\dag},C^{\dag})$ has all the elements that we might look for in a Mackey functor valued in Waldhausen $\infty$-categories: there's a covariant functor $\fromto{C_{\dag}}{\Wald_{\infty}}$ and a contravariant functor $\fromto{(C^{\dag})^{\op}}{\Wald_{\infty}}$, and the two are glued via base change equivalences. Unfortunately, these data are not displayed in a fashion that makes it easy to spot the functoriality in the effective Burnside $\infty$-category $A^{\eff}(C,C_{\dag},C^{\dag})$. In order to extract something that is visibly functorial in the effective Burnside $\infty$-category, we must perform an operation, which we call \emph{unfurling}. When we unfurl a Waldhausen bicartesian fibration $\fromto{\mathscr{X}}{C}$ for the triple $(C,C_{\dag},C^{\dag})$, we end up with a Waldhausen cocartesian fibration $\fromto{\Upsilon(\mathscr{X}/(C,C_{\dag},C^{\dag}))}{A^{\eff}(C,C_{\dag},C^{\dag})}$, which we may then straighten into a Mackey functor on $(C,C_{\dag},C^{\dag})$ valued in Waldhausen $\infty$-categories.

\begin{ntn} Suppose $(C,C_{\dag},C^{\dag})$ a triple, and suppose $p\colon\fromto{X}{C}$ an inner fibration. Denote by $X_{\dag}\subset X\times_{C}C_{\dag}$ (respectively, $X^{\dag}\subset X\times_{C}C^{\dag}$) the subcategory containing all the objects whose morphisms are $p$-cocartesian (resp., $p$-cartesian). 
\end{ntn}

\begin{prp} Suppose $(C,C_{\dag},C^{\dag})$ an adequate triple (Df. \ref{dfn:adequate}), and suppose $p\colon\fromto{X}{C}$ an adequate inner fibration over $(C,C_{\dag},C^{\dag})$ (Df. \ref{dfn:Waldbicart}). Then the triples $(X,X\times_{C}C_{\dag},X^{\dag})$ and $(X,X_{\dag},X\times_CC^{\dag})$ are adequate as well.
\begin{proof} We show that $(X,X\times_{C}C_{\dag},X^{\dag})$ is adequate; the case of  $(X,X_{\dag},X\times_CC^{\dag})$ is dual. Suppose $\sigma\colon\fromto{\Delta^1\times\Delta^1}{C}$ an ambigressive pullback
\begin{equation*}
\begin{tikzpicture} 
\matrix(m)[matrix of math nodes, 
row sep=4ex, column sep=4ex, 
text height=1.5ex, text depth=0.25ex] 
{s'&t'\\ 
s&t,\\}; 
\path[>=stealth,->,font=\scriptsize] 
(m-1-1) edge[>->] (m-1-2) 
edge[->>] (m-2-1) 
(m-1-2) edge[->>] (m-2-2) 
(m-2-1) edge[>->] (m-2-2); 
\end{tikzpicture}
\end{equation*}
and suppose $\cofto{x}{y}$ a morphism covering $\cofto{s}{t}$ and $\fibto{y'}{y}$ a $p$-cartesian edge covering $\fromto{t'}{t}$. Then there exists a $p$-cartesian edge $\fibto{x'}{x}$ covering $\fromto{s'}{s}$, and after filling in the inner horn $x'\to x\to y$ and the outer horn $x'\to y\ot y'$ (using the $p$-cartesianness of $\fibto{y'}{y}$), we obtain a pullback square
\begin{equation*}
\begin{tikzpicture} 
\matrix(m)[matrix of math nodes, 
row sep=4ex, column sep=4ex, 
text height=1.5ex, text depth=0.25ex] 
{x'&y'\\ 
x&y\\}; 
\path[>=stealth,->,font=\scriptsize] 
(m-1-1) edge[>->] (m-1-2) 
edge[->>] (m-2-1) 
(m-1-2) edge[->>] (m-2-2) 
(m-2-1) edge[>->] (m-2-2); 
\end{tikzpicture}
\end{equation*}
covering $\sigma$.
\end{proof}
\end{prp}

\begin{dfn}\label{dfn:unfurl} Suppose $(C,C_{\dag},C^{\dag})$ an adequate triple, and suppose $p\colon\fromto{X}{C}$ an adequate inner fibration over $(C,C_{\dag},C^{\dag})$. Then the \textbf{\emph{unfurling}} of $p$ is the $\infty$-category
\begin{equation*}
\Upsilon(X/(C,C_{\dag},C^{\dag}))\coloneq A^{\eff}(X,X\times_{C}C_{\dag},X^{\dag}).
\end{equation*}
Composition with $p$ defines a natural map
\begin{equation*}
\Upsilon(p)\colon\fromto{\Upsilon(X/(C,C_{\dag},C^{\dag}))}{A^{\eff}(C,C_{\dag},C^{\dag})}.
\end{equation*}
\end{dfn}

We'll prove the following brace of lemmas in the next section.

\begin{lem}\label{lem:unfurlinnerfib} Suppose $(C,C_{\dag},C^{\dag})$ an adequate triple, and suppose $p\colon\fromto{X}{C}$ an adequate inner fibration over $(C,C_{\dag},C^{\dag})$. Then
\begin{equation*}
\Upsilon(p)\colon\fromto{\Upsilon(X/(C,C_{\dag},C^{\dag}))}{A^{\eff}(C,C_{\dag},C^{\dag})}
\end{equation*}
is an inner fibration.
\end{lem}

\begin{lem}\label{lem:unfurlcocartfib} Suppose $(C,C_{\dag},C^{\dag})$ an adequate triple, and suppose $p\colon\fromto{X}{C}$ an adequate inner fibration over $(C,C_{\dag},C^{\dag})$. An edge $f\colon\fromto{y}{z}$ of $\Upsilon(X/(C,C_{\dag},C^{\dag}))$ is $\Upsilon(p)$-cocartesian if it is represented as a span
\begin{equation*}
\begin{tikzpicture}[baseline]
\matrix(m)[matrix of math nodes, 
row sep=3ex, column sep=3ex, 
text height=1.5ex, text depth=0.25ex] 
{&u&\\ 
y&&z,\\}; 
\path[>=stealth,->,font=\scriptsize,inner sep=1.5pt] 
(m-1-2) edge[->>] node[above left]{$\phi$} (m-2-1) 
edge[>->] node[above right]{$\psi$} (m-2-3); 
\end{tikzpicture}
\end{equation*}
in which $\phi$ is $p$-cartesian over an egressive morphism and $\psi$ is $p$-cocartesian over an ingressive morphism.
\end{lem}

The following is now immediate.
\begin{prp} Suppose $(C,C_{\dag},C^{\dag})$ an adequate triple that is either left or right complete, and suppose
\begin{equation*}
p\colon\fromto{\mathscr{X}}{C}
\end{equation*}
an adequate inner fibration over $(C,C_{\dag},C^{\dag})$. Then the unfurling $\Upsilon(p)$ is a cocartesian fibration.
\end{prp}

More particularly, for any disjunctive triple $(C,C_{\dag},C^{\dag})$ that is either left or right complete (Df. \ref{dfn:leftrightcomplete}), and for any Waldhausen bicartesian fibration $p\colon\fromto{\mathscr{X}}{C}$ over $(C,C_{\dag},C^{\dag})$ (Df. \ref{dfn:Waldbicart}), it follows that $\Upsilon(p)$ is a cocartesian fibration. Moreover, for any object $s$ of $C$, consider the functor $\fromto{\mathscr{X}_s}{\Upsilon(\mathscr{X}/(C,C_{\dag},C^{\dag}))_s}$ induced by the natural transformation $\fromto{\varepsilon^{\star}}{\op}$. This functor is the identity on objects, and it is easy to see that it is fully faithful. Now for any $\fromto{s}{t}$ of $A^{\eff}(C,C_{\dag},C^{\dag})$ represented as a span
\begin{equation*}
\begin{tikzpicture}[baseline]
\matrix(m)[matrix of math nodes, 
row sep=3ex, column sep=3ex, 
text height=1.5ex, text depth=0.25ex] 
{&u&\\ 
s&&t,\\}; 
\path[>=stealth,->,font=\scriptsize,inner sep=1.5pt] 
(m-1-2) edge[->>] node[above left]{$f$} (m-2-1) 
edge[>->] node[above right]{$g$} (m-2-3); 
\end{tikzpicture}
\end{equation*}
in $C$, the induced functor
\begin{equation*}
\fromto{\mathscr{X}_s\simeq\Upsilon(\mathscr{X}/(C,C_{\dag},C^{\dag}))_s}{\Upsilon(\mathscr{X}/(C,C_{\dag},C^{\dag}))_t\simeq\mathscr{X}_t}
\end{equation*}
is equivalent to $g_{!}\circ f^{\star}$. In particular, it is exact, whence we have the following.

\begin{prp} Suppose $(C,C_{\dag},C^{\dag})$ a disjunctive triple that is either left or right complete, and suppose
\begin{equation*}
p\colon\fromto{\mathscr{X}}{C}
\end{equation*}
a Waldhausen bicartesian fibration over $(C,C_{\dag},C^{\dag})$. Then the unfurling $\Upsilon(p)$ is a Waldhausen cocartesian fibration.
\end{prp}
\noindent Now condition (\ref{dfn:Waldbicart}.2) immediately implies the main feature of unfurlings.
\begin{thm} Suppose $(C,C_{\dag},C^{\dag})$ a disjunctive triple, and suppose
\begin{equation*}
p\colon\fromto{\mathscr{X}}{C}
\end{equation*}
a Waldhausen bicartesian fibration over $(C,C_{\dag},C^{\dag})$. Then a functor
\begin{equation*}
\mathscr{M}_p\colon\fromto{A^{\eff}(C,C_{\dag},C^{\dag})}{\Wald_{\infty}}
\end{equation*}
that classifies the unfurling $\Upsilon(p)$ is a Mackey functor.
\end{thm}
\noindent One may compose a functor classifying the unfurling of a Waldhausen bicartesian fibration with the delooping $\fromto{\Wald_{\infty}}{\Sp}$ of any additive theory to obtain the following.
\begin{cor} Suppose $(C,C_{\dag},C^{\dag})$ a disjunctive triple, and suppose
\begin{equation*}
p\colon\fromto{\mathscr{X}}{C}
\end{equation*}
a Waldhausen bicartesian fibration over $(C,C_{\dag},C^{\dag})$. Suppose
\begin{equation*}
\mathscr{M}_p\colon\fromto{A^{\eff}(C,C_{\dag},C^{\dag})}{\Wald_{\infty}}
\end{equation*}
a functor that classifies the unfurling $\Upsilon(p)$. Then for any additive theory
\begin{equation*}
F\colon\fromto{\Wald_{\infty}}{\mathscr{E}}
\end{equation*}
in the sense of \cite[Df. 7.1]{K1}, the composition
\begin{equation*}
\FF\circ\mathscr{M}_p\colon\fromto{A^{\eff}(C,C_{\dag},C^{\dag})}{\Sp(\mathscr{E})},
\end{equation*}
where $\FF\colon\fromto{\Wald_{\infty}}{\Sp(\mathscr{E})}$ is the canonical delooping of $F$ \cite[Cor. 7.6.1]{K1}, is a Mackey functor.
\end{cor}

In particular, we see that the algebraic $K$-theory of a Waldhausen bicartesian fibration naturally organizes itself into a Mackey functor valued in spectra.

%----------------------------------------------------------------------%

\section{Horn filling in effective Burnside $\infty$-categories}\label{sect:fillmyhorn} In both Lm. \ref{lem:unfurlinnerfib} and Lm. \ref{lem:unfurlcocartfib}, we are interested in filling horns in effective Burnside $\infty$-categories. These correspond to extensions along inclusions of the form $\into{\widetilde{\mathscr{O}}(\Lambda^m_k)^{\op}}{\widetilde{\mathscr{O}}(\Delta^m)^{\op}}$ that enjoy certain properties. In this section, we construct a filtration that provides a general strategy for constructing the desired extensions, and we use this to prove Lms. \ref{lem:unfurlinnerfib} and \ref{lem:unfurlcocartfib}.

The reader uninterested in such nitty-gritty may be forgiven for skipping this section; however, one must acknowledge that this is where the rubber meets the road. Ultimately, it is the combinatorics of simplices that allow us to solve homotopy-coherence problems.

\begin{ntn} In this section, let $(C,C_{\dag},C^{\dag})$ and $(D,D_{\dag},D^{\dag})$ denote two adequate triples, and let $p\colon\fromto{(C,C_{\dag},C^{\dag})}{(D,D_{\dag},D^{\dag})}$ be a functor of triples that carries ambigressive pullbacks to ambigressive pullbacks. Write $p_{\dag}$ for the restriction $\fromto{C_{\dag}}{D_{\dag}}$, and write $p^{\dag}$ for the restriction $\fromto{C^{\dag}}{D^{\dag}}$.
\end{ntn}

Here is what we will prove.
\begin{thm}\label{thm:pain} Assume that the underlying functor $\fromto{C}{D}$ is an inner fibration. Then the induced functor
\begin{equation*}
A^{\eff}(p)\colon\fromto{A^{\eff}(C,C_{\dag},C^{\dag})}{A^{\eff}(D,D_{\dag},D^{\dag})}
\end{equation*}
is an inner fibration.

Furthermore, assume the following.
\begin{enumerate}[(\ref{thm:pain}.1)]
\item For any ingressive morphism $g\colon\cofto{s}{t}$ of $D$ and any object $x\in C_s$, there exists an ingressive morphism $f\colon\cofto{x}{y}$ of $C$ covering $g$ that is both $p$-cocartesian and $p_{\dag}$-cocartesian.
\item Suppose $\sigma$ a commutative square
\begin{equation*}
\begin{tikzpicture} 
\matrix(m)[matrix of math nodes, 
row sep=4ex, column sep=4ex, 
text height=1.5ex, text depth=0.25ex] 
{x'&y'\\ 
x&y,\\}; 
\path[>=stealth,->,font=\scriptsize] 
(m-1-1) edge[>->] node[above]{$f'$} (m-1-2) 
edge[->>] node[left]{$\phi$} (m-2-1) 
(m-1-2) edge[->] node[right]{$\psi$} (m-2-2) 
(m-2-1) edge[>->] node[below]{$f$} (m-2-2); 
\end{tikzpicture}
\end{equation*}
of $C$ such that the square $f(\sigma)$ is an ambigressive pullback in $D$, the morphism $f'$ is ingressive, the morphism $\phi$ is egressive, and the morphism $f$ is $p$-cocartesian. Then $f'$ is $p$-cocartesian if and only if the square is an ambigressive pullback (and in particular $\psi$ is egressive).
\end{enumerate}
Then an edge $f\colon\fromto{y}{z}$ of $A^{\eff}(C,C_{\dag},C^{\dag})$ is $A^{\eff}(p)$-cocartesian if it is represented as a span
\begin{equation*}
\begin{tikzpicture}[baseline]
\matrix(m)[matrix of math nodes, 
row sep=3ex, column sep=3ex, 
text height=1.5ex, text depth=0.25ex] 
{&u&\\ 
y&&z,\\}; 
\path[>=stealth,->,font=\scriptsize,inner sep=1.5pt] 
(m-1-2) edge[->>] node[above left]{$\phi$} (m-2-1) 
edge[>->] node[above right]{$\psi$} (m-2-3); 
\end{tikzpicture}
\end{equation*}
in which $\phi$ is egressive and $p$-cartesian and $\psi$ is ingressive and $p$-cocartesian.
\end{thm}

\noindent The proof will occupy the entirety of this section. This implies both Lm. \ref{lem:unfurlinnerfib} and Lm. \ref{lem:unfurlcocartfib} as special cases. (Note that (\ref{thm:pain}.2) is an immediate consequence of the Beck--Chevvaley condition (\ref{dfn:Waldbicart}.\ref{item:ambiBC}).) This result may also be used to give an alternative argument that $A^{\eff}(C,C_{\dag},C^{\dag})$ is an $\infty$-category.

\begin{ntn} Suppose $m\geq 2$ and suppose $0\leq k<m$ (note that we are excluding the case $k=m$), and suppose we are given a commutative square
\begin{equation*}
\begin{tikzpicture} 
\matrix(m)[matrix of math nodes, 
row sep=6ex, column sep=4ex, 
text height=1.5ex, text depth=0.25ex] 
{\Lambda^m_k&A^{\eff}(C,C_{\dag},C^{\dag}).\\ 
\Delta^m&A^{\eff}(D,D_{\dag},D^{\dag}),\\}; 
\path[>=stealth,->,font=\scriptsize] 
(m-1-1) edge (m-1-2) 
edge (m-2-1) 
(m-1-2) edge node[right]{$A^{\eff}(p)$} (m-2-2) 
(m-2-1) edge (m-2-2); 
\end{tikzpicture}
\end{equation*}
and we seek a lift $\fromto{\Delta^m}{A^{\eff}(C,C_{\dag},C^{\dag})}$. This corresponds to a (solid) commutative square
\begin{equation*}
\begin{tikzpicture} 
\matrix(m)[matrix of math nodes, 
row sep=4ex, column sep=4ex, 
text height=1.5ex, text depth=0.25ex] 
{\widetilde{\mathscr{O}}(\Lambda^m_k)^{\op}&C\\ 
\widetilde{\mathscr{O}}(\Delta^m)^{\op}&D\\}; 
\path[>=stealth,->,font=\scriptsize] 
(m-1-1) edge node[above]{$g$} (m-1-2) 
edge (m-2-1) 
(m-1-2) edge node[right]{$p$} (m-2-2) 
(m-2-1) edge node[below]{$h$} (m-2-2)
(m-2-1) edge[dotted,inner sep=0.5] node[above left]{$\overline{g}$} (m-1-2); 
\end{tikzpicture}
\end{equation*}
of simplicial sets in which $h$ is ambigressive cartesian and, for every $i\neq k$, the restriction
\begin{equation*}
\widetilde{\mathscr{O}}(\Delta^{\{0,\dots,\hat\imath,\dots,m\}})^{\op}\subset\widetilde{\mathscr{O}}(\Lambda^m_k)^{\op}\ \tikz[baseline]\draw[>=stealth,->,font=\scriptsize](0,0.5ex)--node[above]{$g$}(0.5,0.5ex);\ C
\end{equation*}
is ambigressive cartesian. Our objective then becomes to construct a (dotted) lift
\begin{equation*}
\overline{g}\colon\fromto{\widetilde{\mathscr{O}}(\Delta^m)^{\op}}{C}
\end{equation*}
that is ambigressive cartesian.
\end{ntn}

\begin{dfn} Let us call a $m$-simplex
\begin{equation*}
i_0j_0\to i_1j_1\to\cdots\to i_mj_m
\end{equation*}
of $\widetilde{\mathscr{O}}(\Delta^m)^{\op}$ \textbf{\emph{completely factored}} if $i_r-i_{r-1}+j_{r-1}-j_r=1$ for each $1\leq r\leq m$. Note that any complete $m$-simplex is in particular nondegenerate, and $i_{0}=0$ and $j_{0}=m$.
\end{dfn}

\begin{nul} Completely factored $m$-simplices are essentially the same as walks on the poset $\widetilde{\mathscr{O}}([\mathbf{m}])^{\op}$ that begin at the point $0m$ and end at a point of the form $pp$. We may therefore parametrize the completely factored $m$-simplices as follows. For each integer $0\leq N\leq 2^m-1$, let $\sigma(N)$ be the unique completely factored $m$-simplex
\begin{equation*}
i_0j_0\to i_1j_1\to\cdots\to i_mj_m
\end{equation*}
such that if $d_r$ is the $r$-th binary digit of $N$ (read left to right, so that $N=\sum_{s=1}^{m}2^{m-s}d_s$), then
\begin{equation*}
d_r=\begin{cases}
0&\textrm{if }i_{r-1}=i_{r};\\
1&\textrm{if }j_{r-1}=j_r.
\end{cases}
\end{equation*}
We order these simplices accordingly. Hence $\sigma(N)$ is the $m$-simplex
\begin{equation*}
i_0j_0\to i_1j_1\to\cdots\to i_mj_m
\end{equation*}
with
\begin{equation*}
i_r=\sum_{s=1}^rd_s\textrm{\quad and\quad}j_r=(m-r)+\sum_{s=1}^rd_s.
\end{equation*}
Fig. \ref{fig:completelyfactored} shows the completely factored $5$-simplex $\sigma(01101)=\sigma(13)\subset\widetilde{\mathscr{O}}(\Delta^m)^{\op}$.

\begin{figure}
\begin{tikzpicture} 
\matrix(m)[matrix of math nodes, 
row sep={6ex,between origins}, column sep={6ex,between origins}, 
text height=1.5ex, text depth=0.25ex] 
{&&&&&\textcolor{red}{05}&&&&&\\
&&&&\textcolor{red}{04}&&15&&&&\\
&&&03&&\textcolor{red}{14}&&25&&&\\
&&02&&13&&\textcolor{red}{24}&&35&&\\
&01&&12&&\textcolor{red}{23}&&34&&45&\\
00&&11&&22&&\textcolor{red}{33}&&44&&55\\}; 
\path[>=stealth,<-,font=\scriptsize] 
(m-2-5) edge[color=red,thick] (m-1-6) 
(m-2-7) edge (m-1-6)
(m-3-4) edge (m-2-5)
(m-3-6) edge[color=red,thick] (m-2-5)
(m-3-6) edge (m-2-7)
(m-3-8) edge (m-2-7)
(m-4-3) edge (m-3-4)
(m-4-5) edge (m-3-4)
(m-4-5) edge (m-3-6)
(m-4-7) edge[color=red,thick] (m-3-6)
(m-4-7) edge (m-3-8)
(m-4-9) edge (m-3-8)
(m-5-2) edge (m-4-3)
(m-5-4) edge (m-4-3)
(m-5-4) edge (m-4-5)
(m-5-6) edge (m-4-5)
(m-5-6) edge[color=red,thick] (m-4-7)
(m-5-8) edge (m-4-7)
(m-5-8) edge (m-4-9)
(m-5-10) edge (m-4-9)
(m-6-1) edge (m-5-2)
(m-6-3) edge (m-5-2)
(m-6-3) edge (m-5-4)
(m-6-5) edge (m-5-4)
(m-6-5) edge (m-5-6)
(m-6-7) edge[color=red,thick] (m-5-6)
(m-6-7) edge (m-5-8)
(m-6-9) edge (m-5-8)
(m-6-9) edge (m-5-10) 
(m-6-11) edge (m-5-10);
\end{tikzpicture}
\caption{The completely factored $5$-simplex $\sigma(01101)=\sigma(13)\subset\widetilde{\mathscr{O}}(\Delta^5)^{\op}$, drawn in red as a walk from $05$ to $33$. The juts of this $5$-simplex are $3$ and $5$. There are no crossings away from $1$; if $k\neq1$, the only crossing away from $1$ is $2$.}\label{fig:completelyfactored}
\end{figure}
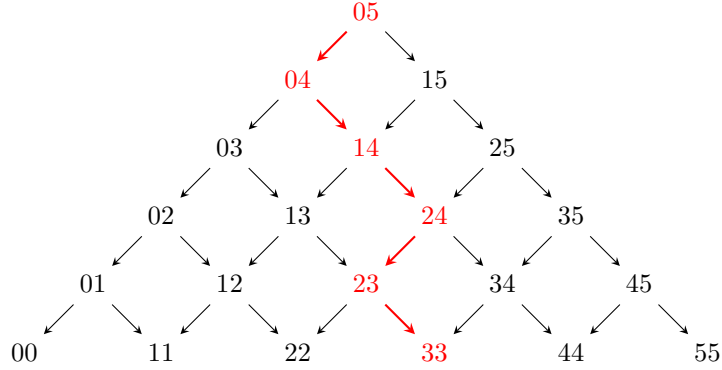

Now for any integer $0\leq N\leq 2^m$, write
\begin{equation*}
P_N(k)\coloneq\widetilde{\mathscr{O}}(\Lambda^m_k)^{\op}\cup\bigcup_{0\leq K<N}\sigma(K);
\end{equation*}
this provides a filtration
\begin{equation*}
\widetilde{\mathscr{O}}(\Lambda^m_k)^{\op}=P_{0}(k)\subset\cdots\subset P_{2^m}(k)=\widetilde{\mathscr{O}}(\Delta^m)^{\op}.
\end{equation*}
Our aim is to find conditions that permit us to extend $g$ along this filtration.
\end{nul}

We proceed to analyze the intersections $\sigma(N)\cap P_N(k)$ as a subset of $\sigma(N)\cong\Delta^m$. We will find that each $\sigma(N)\cap P_N(k)$ is in fact a union of faces of $\Delta^m$. There are two kinds of faces that will appear in this intersection. To describe these, let's introduce some simplifying notation.

\begin{ntn} For any nonempty totally ordered finite set $T$ and any element $j\in T$, write $\Delta^{\hat{\jmath}}$ for the face
\begin{equation*}
\Delta^{T-\{j\}}\subset\Delta^T.
\end{equation*}

More generally, for any ordered subsets $S\subset T$, write $\Lambda^T_S\subset\Delta^T$ for the union of all the faces (i.e., $\#T$-simplices) of $\Delta^T$ that contain the simplex $\Delta^{S}$. In other words, let
\begin{equation*}
\Lambda^T_S\coloneq\bigcup_{j\notin S}\Delta^{\hat{\jmath}}.
\end{equation*}
When $T=\{0,\dots,m\}$, we just write $\Lambda^m_S$ for $\Lambda^T_S$.
\end{ntn}
\noindent In this notation, we have
\begin{equation*}
\sigma(N)\cap P_N(k)=\left(\bigcup_{0\leq K<N}\sigma(N)\cap\sigma(K)\right)\cup\left(\bigcup_{j\neq k}\sigma(N)\cap\widetilde{\mathscr{O}}(\Delta^{\hat{\jmath}})^{\op}\right),
\end{equation*}
and our claim is that there is a set $E(N,k)\subset\{0,\dots,m\}$ such that
\begin{equation*}
\sigma(N)\cap P_N(k)=\Lambda^m_{E(N,k)}.
\end{equation*}
We proceed to construct this set now.

\begin{dfn} Suppose $N$ an integer such that $0\leq N\leq 2^{m}-1$, written as $N=\sum_{s=1}^m2^{m-s}d_s$. A \textbf{\emph{jut}} of the completely factored $m$-simplex $\sigma(N)$ is an integer $z\in\{1,\dots,m\}$ such that 
\begin{itemize}
\item $d_{z}=1$, and
\item either $d_{z+1}=0$ or $z=m$.
\end{itemize}
Denote by $Z(N)\subset\{1,\dots,m\}$ the set of juts of $\sigma(N)$.

For any jut $z$, write
\begin{equation*}
K_z=\sum_{s=1}^m2^{m-s}d_{z,s},
\end{equation*}
where
\begin{equation*}
d_{z,s}=\begin{cases}
d_s&\textrm{if }s\notin\{z,z+1\};\\
0&\textrm{if }s=z;\\
1&\textrm{if }s=z+1.
\end{cases}
\end{equation*}
\end{dfn}

\begin{lem} Suppose $N$ an integer such that $0\leq N\leq 2^{m}-1$, written as $N=\sum_{s=1}^m2^{m-s}d_s$.  Then
\begin{equation*}
\bigcup_{0\leq K<N}\sigma(N)\cap\sigma(K)=\bigcup_{z\in Z(N)}\Delta^{\hat{z}}.
\end{equation*}
\begin{proof} It is easy to see that in the poset of simplicial subsets of $\sigma(N)$ of the form $\sigma(N)\cap\sigma(K)$, the maximal elements consist of those subsets of the form $\sigma(N)\cap\sigma(K_z)$, where $z$ in a jut of $\sigma(N)$. Of course $\Delta^{\hat{z}}=\sigma(N)\cap\sigma(K_z)$.
\end{proof}
\end{lem}

\begin{dfn} Suppose $N$ an integer such that $0\leq N\leq 2^{m}-1$ written as $N=\sum_{s=1}^m2^{m-s}d_s$. A \textbf{\emph{crossing}} of $\sigma(N)$ away from $k$ is an integer $x\in\{0,\dots,m-1\}$ such that one of the following holds:
\begin{itemize}
\item $x=0$, $d_1=0$;
\item $x=0$, $d_1=1$, and $k\neq 0$;
\item $x>0$, $d_x=d_{x+1}=1$, and $i_x\neq k$; or
\item $x>0$, $d_x=d_{x+1}=0$ and $j_x\neq k$.
\end{itemize}
Denote by $X(N,k)\subset\{0,\dots,m-1\}$ the set of crossings away from $k$.
\end{dfn}

The crossings away from $k$ are now all we need to complete our computation of the intersections $\sigma(N)\cap P_N(k)$.
\begin{prp} Suppose $N$ an integer such that $0\leq N\leq 2^{m}-1$, written as $N=\sum_{s=1}^m2^{m-s}d_s$. Then
\begin{equation*}
\sigma(N)\cap P_N(k)=\bigcup_{y\in Z(N)\cup X(N,k)}\Delta^{\hat{y}}.
\end{equation*}
\begin{proof} For any crossing $x$ of $\sigma(N)$ away from $k$, it is clear that the corresponding face is given by
\begin{equation*}
\Delta^{\hat{x}}=\begin{cases}
\sigma(N)\cap\widetilde{\mathscr{O}}(\Delta^{\hat{\imath}_x})^{\op}&\textrm{if }d_x=1;\\
\sigma(N)\cap\widetilde{\mathscr{O}}(\Delta^{\hat{\jmath}_x})^{\op}&\textrm{if }d_x=0.
\end{cases}
\end{equation*}
Now for any $j\in\{0,\dots,m\}$, if the set
\begin{equation*}
\{r\in\{0,\dots,m\}\ |\ i_r=j\textrm{ or }j_r=j\}
\end{equation*}
contains more than one element, then it contains a jut $z$, and consequently,
\begin{equation*}
\sigma(N)\cap\widetilde{\mathscr{O}}(\Delta^{\hat{\jmath}})^{\op}\subset\Delta^{\hat{z}}=\sigma(N)\cap\sigma(K_z).\qedhere
\end{equation*}
\end{proof}
\end{prp}

\begin{wrn} Note that this doesn't quite work if $k=m$ (which we expressly excluded): in this case it is just not true that $\sigma(N)\cap P_N(k)$ is a union of faces. For example, in the completely factored $5$-simplex $\sigma(13)$ depicted in Fig. \ref{fig:completelyfactored}, if $k=5$, then the simplicial subset $\sigma(13)\cap P_{13}(5)\subset\sigma(13)$ is the union
\begin{equation*}
\Delta^{\hat{3}}\cup\Delta^{\hat{5}}\cup\Delta^{\{2,3,4,5\}}.
\end{equation*}
\end{wrn}

Let us reformulate what we have shown.
\begin{dfn} Suppose $N$ an integer such that $0\leq N\leq 2^{m}-1$ written as $N=\sum_{s=1}^m2^{m-s}d_s$. Let us call an integer $s\in\{0,\dots,m\}$ an \textbf{\emph{essential vertex}} of $\sigma(N)$ for $k$ if it is neither a jut nor a crossing away from $k$. Denote by
\begin{equation*}
E(N,k)\coloneq\{0,\dots,m\}-(Z(N)\cup X(N,k))
\end{equation*}
the ordered set of essential vertices of $\sigma(N)$ for $k$.
\end{dfn}

We have thus shown that we may write
\begin{equation*}
\sigma(N)\cap P_N(k)=\Lambda^m_{E(N,k)}\subset\Delta^m\cong\sigma(N).
\end{equation*}
Now we want to extend $g$ along each inclusion $\into{P_{N-1}(k)}{P_N(k)}$, which we now write as pushout
\begin{equation*}
\begin{tikzpicture} 
\matrix(m)[matrix of math nodes, 
row sep=4ex, column sep=4ex, 
text height=1.5ex, text depth=0.25ex] 
{\Lambda^m_{E(N,k)}&\Delta^m\\ 
P_{N}(k)&P_{N+1}(k).\\}; 
\path[>=stealth,->,font=\scriptsize] 
(m-1-1) edge (m-1-2) 
edge (m-2-1) 
(m-1-2) edge (m-2-2) 
(m-2-1) edge (m-2-2); 
\end{tikzpicture}
\end{equation*}
For this, we need to determine just what sort of inclusions these ore. For example, if $\into{\Lambda^m_{E(N,k)}}{\Delta^m}$ is inner anodyne, one has the desired extension simply because $p$ is an inner fibration. Let's determine precisely when this does the job.

\begin{lem} Suppose $S\subset\{0,\dots,m\}$ a nonempty ordered subset. Then the inclusion $\into{\Lambda^m_S}{\Delta^m}$ is inner anodyne if the following condition holds.
\begin{itemize}
\item[($\ast$)] there exists elements $a<s<b$ of $\{0,\dots,m\}$ such that $s\in S$, but $a,b\notin S$.
\end{itemize}
\begin{proof} The claim is trivial if either $m=2$ or $S$ has cardinality $1$. For $m\geq3$ and $\# S\geq 2$, assume that the result holds both for all smaller values of $m$ and for subsets $S$ of smaller cardinality.

Choose an element $s\in S$ as follows: if $0\in S$, let $s=0$; otherwise, if $m\in S$, let $s=m$; otherwise choose $s\in S$ arbitrarily. Then the subset
\begin{equation*}
S-\{s\}\subset\{0,\dots,\widehat{s},\dots,m\}
\end{equation*}
satisfies condition ($\ast$) for $m-1$; hence the inclusion
\begin{equation*}
\into{\Lambda^m_S\cap\Delta^{\hat{s}}=\Lambda^{\{0,\dots,\widehat{s},\dots,m\}}_{S-\{s\}}}{\Delta^{\{0,\dots,\widehat{s},\dots,m\}}}
\end{equation*}
is inner anodyne by the inductive hypothesis. The pushout of this edge along the inclusion
\begin{equation*}
\into{\Lambda^m_S\cap\Delta^{\hat{s}}}{\Lambda^m_S}
\end{equation*}
is the inclusion
\begin{equation*}
\into{\Lambda^m_S}{\Lambda^m_{S-\{s\}}},
\end{equation*}
which is thus inner anodyne. Now our claim follows from the observation that the subset $S-\{s\}\subset\{0,\dots,m\}$ also satisfies condition ($\ast$), whence the inclusion
\begin{equation*}
\into{\Lambda^m_{S-\{s\}}}{\Delta^m}
\end{equation*}
is also inner anodyne by the inductive hypothesis.
\end{proof}
\end{lem}

Suppose $N$ an integer such that $0\leq N\leq 2^{m}-1$ written as $N=\sum_{s=1}^m2^{m-s}d_s$. We have shown that if $E(N,k)$ satisfies condition ($\ast$), then $\into{P_N(k)}{P_{N+1}(k)}$ is inner anodyne. If however $E(N,k)$ fails condition ($\ast$), let's refer to $\sigma(N)$ as an \textbf{\emph{exceptional case}}.

Indeed, such an $m$-simplex is quite exceptional: it cannot contain more than one jut, since there must be an essential vertex between any two juts. Consequently, for any $t\in\{0,\dots,m\}$, if
\begin{equation*}
d_s=\begin{cases}
1&\textrm{if }s\leq t;\\
0&\textrm{if }s>t,
\end{cases}
\end{equation*}
then let's write $N_t\coloneq\sum_{s=1}^m2^{m-s}d_s$. Then $\sigma(N_t)$ is the walk
\begin{equation*}
0m\to 1m\to\cdots\to tm\to t(m-1)\cdots\to tt.
\end{equation*}
One thus sees that any exceptional case is of the form $\sigma(N_t)$.

Furthermore, the simplex $\sigma(N_t)$ begins with a crossing if either $k\neq 0$ or $N_t<2^{m-1}$, so $0\in E(N_t,k)$ if and only if both $k=0$ and $N_t\geq 2^{m-1}$. Dually, $\sigma(N_t)$ ends in a jut precisely when $N_t$ is odd, so $m\in E(k)$ if and only if $N_t$ is even. Now a quick analysis of the location of the crossings away from $k$ yields the following classification of all the exceptional cases.
\begin{prp}\label{prp:exceptions} When $k\neq0$, there are only two exceptional cases:
\begin{enumerate}[(\ref{prp:exceptions}.1)]
\item $t=k-1$, in which case $E(N_{k-1},k)=\{m-1,m\}$.
\item $t=k$, in which case $E(N_k,k)=\{m\}$.
\suspend{enumerate}
When $k=0$, there are $t+1$ exceptional cases:
\resume{enumerate}[{[(\ref{prp:exceptions}.1)]}]
\item $t=0$, in which case $E(0,0)=\{m\}$.
\item $0<t<m$, in which case $E(N_t,0)=\{0,m\}$.
\item $t=m$, in which case $E(2^m-1,0)=\{0\}$.
\end{enumerate}
\end{prp}

\begin{figure}
\begin{tikzpicture} 
\matrix(m)[matrix of math nodes, 
row sep={6ex,between origins}, column sep={6ex,between origins}, 
text height=1.5ex, text depth=0.25ex] 
{&&&&&\textcolor{violet}{05}&&&&&\\
&&&&{04}&&\textcolor{violet}{15}&&&&\\
&&&03&&{14}&&\textcolor{violet}{25}&&&\\
&&02&&13&&\textcolor{red}{24}&&\textcolor{blue}{35}&&\\
&01&&12&&\textcolor{red}{23}&&\textcolor{blue}{34}&&45&\\
00&&11&&\textcolor{red}{22}&&\textcolor{blue}{33}&&44&&55\\}; 
\path[>=stealth,<-,font=\scriptsize] 
(m-2-5) edge (m-1-6) 
(m-2-7) edge[color=violet,thick] (m-1-6)
(m-3-4) edge (m-2-5)
(m-3-6) edge (m-2-5)
(m-3-6) edge (m-2-7)
(m-3-8) edge[color=violet,thick] (m-2-7)
(m-4-3) edge (m-3-4)
(m-4-5) edge (m-3-4)
(m-4-5) edge (m-3-6)
(m-4-7) edge (m-3-6)
(m-4-7) edge[color=red,thick] (m-3-8)
(m-4-9) edge[color=blue,thick] (m-3-8)
(m-5-2) edge (m-4-3)
(m-5-4) edge (m-4-3)
(m-5-4) edge (m-4-5)
(m-5-6) edge (m-4-5)
(m-5-6) edge[color=red,thick] (m-4-7)
(m-5-8) edge (m-4-7)
(m-5-8) edge[color=blue,thick] (m-4-9)
(m-5-10) edge (m-4-9)
(m-6-1) edge (m-5-2)
(m-6-3) edge (m-5-2)
(m-6-3) edge (m-5-4)
(m-6-5) edge (m-5-4)
(m-6-5) edge[color=red,thick] (m-5-6)
(m-6-7) edge (m-5-6)
(m-6-7) edge[color=blue,thick] (m-5-8)
(m-6-9) edge (m-5-8)
(m-6-9) edge (m-5-10) 
(m-6-11) edge (m-5-10);
\end{tikzpicture}
\caption{The two exceptional cases $\sigma(N_2)$ and $\sigma(N_3)$ for $k=3$ in $\widetilde{\mathscr{O}}(\Delta^5)^{\op}$, drawn in red and blue (respectively) as walks from $05$ to $22$ and from $05$ to $33$.}\label{fig:exceptional}
\end{figure}
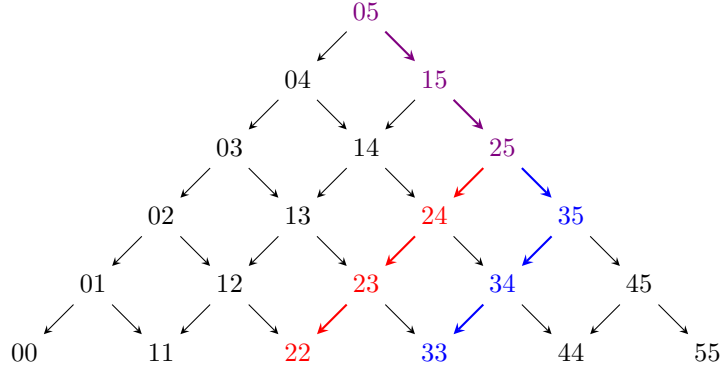

\noindent To illustrate, in Fig. \ref{fig:exceptional}, the two exceptional cases in $\widetilde{\mathscr{O}}(\Delta^5)^{\op}$ for $k=3$ are depicted.

We can now begin the proof of Th. \ref{thm:pain}. Here is the first bit.
\begin{lem} Suppose $0<k<m$. Then there exists a (dotted) lift
\begin{equation*}
\begin{tikzpicture} 
\matrix(m)[matrix of math nodes, 
row sep=4ex, column sep=4ex, 
text height=1.5ex, text depth=0.25ex] 
{\widetilde{\mathscr{O}}(\Lambda^m_k)^{\op}&C\\ 
\widetilde{\mathscr{O}}(\Delta^m)^{\op}&D\\}; 
\path[>=stealth,->,font=\scriptsize] 
(m-1-1) edge node[above]{$g$} (m-1-2) 
edge (m-2-1) 
(m-1-2) edge node[right]{$p$} (m-2-2) 
(m-2-1) edge node[below]{$h$} (m-2-2)
(m-2-1) edge[dotted,inner sep=0.5] node[above left]{$\overline{g}$} (m-1-2); 
\end{tikzpicture}
\end{equation*}
\begin{proof} Let's handle the case $m=2$ separately. Note that $\widetilde{\mathscr{O}}(\Delta^2)^{\op}=(\widetilde{\mathscr{O}}(\Delta^2_1)^{\op})^{\lhd}$, so we may form the desired lift $\overline{g}$ simply by forming the $p$-limit in the sense of \cite[Df. 4.3.1.1]{HTT}.

For $m\geq 3$, of course we will proceed by induction on the filtration
\begin{equation*}
\widetilde{\mathscr{O}}(\Lambda^m_k)^{\op}=P_{0}(k)\subset\cdots\subset P_{2^m}(k)=\widetilde{\mathscr{O}}(\Delta^m)^{\op}.
\end{equation*}
Given a lift
\begin{equation*}
\begin{tikzpicture} 
\matrix(m)[matrix of math nodes, 
row sep=4ex, column sep=8ex, 
text height=1.5ex, text depth=0.25ex] 
{\widetilde{\mathscr{O}}(\Lambda^m_k)^{\op}&C\\ 
P_N(k)&D,\\}; 
\path[>=stealth,->,font=\scriptsize] 
(m-1-1) edge node[above]{$g$} (m-1-2) 
edge (m-2-1) 
(m-1-2) edge node[right]{$p$} (m-2-2) 
(m-2-1) edge node[below]{$h|P_N(k)$} (m-2-2)
(m-2-1) edge[inner sep=0.5] node[below right]{$g_N$} (m-1-2); 
\end{tikzpicture}
\end{equation*}
we seek a (dotted) lift
\begin{equation*}
\begin{tikzpicture} 
\matrix(m)[matrix of math nodes, 
row sep=4ex, column sep=10ex, 
text height=1.5ex, text depth=0.25ex] 
{P_N(k)&C\\ 
P_{N+1}(k)&D.\\}; 
\path[>=stealth,->,font=\scriptsize] 
(m-1-1) edge node[above]{$g_N$} (m-1-2) 
edge (m-2-1) 
(m-1-2) edge node[right]{$p$} (m-2-2) 
(m-2-1) edge node[below]{$h|P_{N+1}(k)$} (m-2-2)
(m-2-1) edge[dotted,inner sep=0.5] node[below right]{$g_{N+1}$} (m-1-2); 
\end{tikzpicture}
\end{equation*}
The only catch will be that we must choose the extensions to the exceptional $m$-simplices and some of their neighbors carefully.

To begin, for $0\leq N<N_{k-1}$, we use the right lifting property with respect to the inner anodyne inclusions $\into{\sigma(N)\cap P_N(k)}{\sigma(N)}$ to obtain the desired lift
\begin{equation*}
g_{N_{k-1}}\colon\fromto{P_{N_{k-1}}(k)}{C}.
\end{equation*}

Now let us call a completely factored $m$-simplex $\sigma(N)$ \textbf{\emph{special}} just in case the corresponding integer $N=\sum_{s=1}^m2^{m-s}d_s$ has the property that $d_s=1$ for every $1\leq s\leq k-1$, and no more than one of $d_k,\dots,d_m$ is equal to $1$. Let $R(k)$ be the collection of those $N$ such that $\sigma(N)$ is special. Note that the exceptional $m$-simplex $\sigma(N_{k-1})$ of (\ref{prp:exceptions}.1) is the first special $m$-simplex, and the exceptional $m$-simplex $\sigma(N_k)$ of (\ref{prp:exceptions}.2) is the last special $m$-simplex. Also observe that for any $N\in R(k)$, one has
\begin{equation*}
\sigma(N)\cap P_N=(\sigma(N)\cap P_{N_{k-1}})\cup\bigcup_{K\in R(k),\ K<N}(\sigma(N)\cap\sigma(K)).
\end{equation*}
Now we have a functor
\begin{equation*}
f_k\colon\fromto{\Delta^1\times\Delta^{m-1}}{\widetilde{\mathscr{O}}(\Delta^m)^{\op}}
\end{equation*}
that is determined on its value on objects:
\begin{equation*}
f_k(u,v)\coloneq(\max\{0,v-k+1\},\max\{m-v-u,k-u\}).
\end{equation*}
The functor $f_k$ restricts to a functor
\begin{equation*}
\fromto{(\Delta^1\times\Lambda^{m-1}_{m-1})\cup^{(\Delta^{\{1\}}\times\Lambda^{m-1}_{m-1})}(\Delta^{\{1\}}\times\Delta^{m-1})}{P_{N_{k-1}}(k)}.
\end{equation*}
Note that (dotted) lifts in the diagram
\begin{equation*}
\begin{tikzpicture} 
\matrix(m)[matrix of math nodes, 
row sep=4ex, column sep=10ex, 
text height=1.5ex, text depth=0.25ex] 
{(\Delta^1\times\Lambda^{m-1}_{m-1})\cup^{(\Delta^{\{1\}}\times\Lambda^{m-1}_{m-1})}(\Delta^{\{1\}}\times\Delta^{m-1})&C\\ 
\Delta^1\times\Delta^{m-1}&D\\}; 
\path[>=stealth,->,font=\scriptsize] 
(m-1-1) edge node[above]{$g_{N_{k-1}}$} (m-1-2) 
edge (m-2-1) 
(m-1-2) edge node[right]{$p$} (m-2-2) 
(m-2-1) edge node[below]{$h|(\Delta^1\times\Delta^{m-1})$} (m-2-2)
(m-2-1) edge[dotted,inner sep=0.5] node[below right]{$g'$} (m-1-2); 
\end{tikzpicture}
\end{equation*}
are in bijection with (dotted) lifts in the adjoint diagram
\begin{equation*}
\begin{tikzpicture} 
\matrix(m)[matrix of math nodes, 
row sep=4ex, column sep=10ex, 
text height=1.5ex, text depth=0.25ex] 
{\Lambda^{m-1}_{m-1}&\Fun(\Delta^1,C)\\ 
\Delta^{m-1}&\Fun(\Delta^{\{1\}},C)\times_{\Fun(\Delta^{\{1\}},D)}\Fun(\Delta^1,D).\\}; 
\path[>=stealth,->,font=\scriptsize] 
(m-1-1) edge (m-1-2) 
edge (m-2-1) 
(m-1-2) edge node[right]{$(t,p)$} (m-2-2) 
(m-2-1) edge node[below]{$\eta$} (m-2-2)
(m-2-1) edge[dotted,inner sep=0.5] (m-1-2); 
\end{tikzpicture}
\end{equation*}
Such a lift exists by \cite[Lm. 6.1.1.1]{HTT}. Consequently, a lift $g'$ exists, and it specifies a family of maps
\begin{equation*}
g'_N\colon\fromto{\sigma(N)}{C}
\end{equation*}
for $N\in R(k)$ with $N>N_{k-1}$.

Now to obtain a (dotted) lift
\begin{equation*}
\begin{tikzpicture} 
\matrix(m)[matrix of math nodes, 
row sep=4ex, column sep=14ex, 
text height=1.5ex, text depth=0.25ex] 
{P_{N_{k-1}}(k)&C\\ 
P_{N_{k-1}+1}(k)&D,\\}; 
\path[>=stealth,->,font=\scriptsize] 
(m-1-1) edge node[above]{$g_{N_{k-1}}$} (m-1-2) 
edge (m-2-1) 
(m-1-2) edge node[right]{$p$} (m-2-2) 
(m-2-1) edge node[below]{$h|P_{N_{k-1}+1}(k)$} (m-2-2)
(m-2-1) edge[dotted,inner sep=0.5] node[below right]{$g_{N_{k-1}+1}$} (m-1-2); 
\end{tikzpicture}
\end{equation*}
we must extend along
\begin{equation*}
\into{\sigma(N_{k-1})\cap P_{N_{k-1}}\cong\Lambda^m_{\{m-1,m\}}}{\Delta^m\cong\sigma(N_{k-1})},
\end{equation*}
which we factor as the composite
\begin{equation*}
\Lambda^m_{\{m-1,m\}}\ \tikz[baseline]\draw[>=stealth,right hook->](0,0.5ex)--(0.5,0.5ex);\ \Lambda^m_{\{m-1,m\}}\cup^{\Lambda_{m-1}^{\{0,\dots,m-1\}}}\Delta^{\{0,\dots,m-1\}}\cong\Lambda^{m}_{m-1}\ \tikz[baseline]\draw[>=stealth,right hook->](0,0.5ex)--(0.5,0.5ex);\ \Delta^m.
\end{equation*}
We extend across the first inclusion by using the restriction of the lift $g'$ we have constructed, and then we extend across the inner horn $\into{\Lambda^m_{m-1}}{\Delta^m}$ using the fact that $p$ is an inner fibration.

Now for any $N>N_{k-1}$, extend across $P_{N+1}$ as follows: if $N\in R(k)$, extend using the chosen map $g'_N$, which since
\begin{equation*}
\sigma(N)\cap P_N=(\sigma(N)\cap P_{N_{k-1}})\cup\bigcup_{K\in R(k),\ K<N}(\sigma(N)\cap\sigma(K)),
\end{equation*}
is compatible with the map $\fromto{P_N}{C}$ constructed so far. If $N\notin R(k)$, simply extend using the fact that $\into{\sigma(N)\cap P_N}{\sigma(N)}$ is an inner horn inclusion. At the end of this procedure, the desired extension
\begin{equation*}
\overline{g}\colon\fromto{P_{2^m}=\widetilde{\mathscr{O}}(\Delta^m)^{\op}}{C}
\end{equation*}
is constructed, and it is ambigressive cartesian by construction.
\end{proof}
\end{lem}

Now let us complete the proof of the theorem.

\begin{lem} Suppose $p$ satisfies conditions (\ref{thm:pain}.1-3) and that $k=0$. If the morphism $\fibto{g(01)}{g(00)}$ is $p$-cartesian and the morphism $\cofto{g(01)}{g(11)}$ is $p$-cocartesian, then there is an ambigressive cartesian (dotted) lift
\begin{equation*}
\begin{tikzpicture} 
\matrix(m)[matrix of math nodes, 
row sep=4ex, column sep=4ex, 
text height=1.5ex, text depth=0.25ex] 
{\widetilde{\mathscr{O}}(\Lambda^m_k)^{\op}&C\\ 
\widetilde{\mathscr{O}}(\Delta^m)^{\op}&D\\}; 
\path[>=stealth,->,font=\scriptsize] 
(m-1-1) edge node[above]{$g$} (m-1-2) 
edge (m-2-1) 
(m-1-2) edge node[right]{$p$} (m-2-2) 
(m-2-1) edge node[below]{$h$} (m-2-2)
(m-2-1) edge[dotted,inner sep=0.5] node[above left]{$\overline{g}$} (m-1-2); 
\end{tikzpicture}
\end{equation*}
\begin{proof} Again let's treat the case $m=2$ separately. In this case, since the morphism $\fibto{g(0,1)}{g(0,0)}$ is $p$-cartesian, we obtain a $2$-simplex
\begin{equation*}
\begin{tikzpicture} 
\matrix(m)[matrix of math nodes, 
row sep=3ex, column sep=0ex, 
text height=1.5ex, text depth=0.25ex] 
{&g(0,1)&\\ 
g(0,2)&&g(0,0).\\}; 
\path[>=stealth,->>,font=\scriptsize] 
(m-2-1) edge (m-1-2) 
edge (m-2-3)
(m-1-2) edge (m-2-3); 
\end{tikzpicture}
\end{equation*}
Now after filling inner horns, we choose a $p$-cocartesian edge $\cofto{g(0,2)}{g(1,2)}$ lying over $\cofto{h(0,2)}{h(1,2)}$, and then by filling the corresponding outer horns, we obtain a diagram
\begin{equation*}
\begin{tikzpicture} 
\matrix(m)[matrix of math nodes, 
row sep=3ex, column sep=0ex, 
text height=1.5ex, text depth=0.25ex] 
{&&g(0,2)&&\\
&g(0,1)&&g(0,2)&\\
g(0,0)&&g(1,1)&&g(2,2)\\}; 
\path[>=stealth,->,font=\scriptsize] 
(m-1-3) edge[->>] (m-2-2) 
edge[>->] (m-2-4)
(m-2-2) edge[->>] (m-3-1)
edge[>->] (m-3-3)
(m-2-4) edge[->>] (m-3-3)
edge[>->] (m-3-5); 
\end{tikzpicture}
\end{equation*}
It follows from conditions (\ref{thm:pain}.1-2) that the morphisms are ingressive or egressive as marked and that the square is ambigressive cartesian.

For $m\geq 3$, we will once again proceed by induction on the filtration
\begin{equation*}
\widetilde{\mathscr{O}}(\Lambda^m_k)^{\op}=P_{0}(k)\subset\cdots\subset P_{2^m}(k)=\widetilde{\mathscr{O}}(\Delta^m)^{\op}.
\end{equation*}
Given a lift
\begin{equation*}
\begin{tikzpicture} 
\matrix(m)[matrix of math nodes, 
row sep=4ex, column sep=8ex, 
text height=1.5ex, text depth=0.25ex] 
{\widetilde{\mathscr{O}}(\Lambda^m_k)^{\op}&C\\ 
P_N(k)&D,\\}; 
\path[>=stealth,->,font=\scriptsize] 
(m-1-1) edge node[above]{$g$} (m-1-2) 
edge (m-2-1) 
(m-1-2) edge node[right]{$p$} (m-2-2) 
(m-2-1) edge node[below]{$h|P_N(k)$} (m-2-2)
(m-2-1) edge[inner sep=0.5] node[below right]{$g_N$} (m-1-2); 
\end{tikzpicture}
\end{equation*}
we seek a (dotted) lift
\begin{equation*}
\begin{tikzpicture} 
\matrix(m)[matrix of math nodes, 
row sep=4ex, column sep=10ex, 
text height=1.5ex, text depth=0.25ex] 
{P_N(k)&C\\ 
P_{N+1}(k)&D.\\}; 
\path[>=stealth,->,font=\scriptsize] 
(m-1-1) edge node[above]{$g_N$} (m-1-2) 
edge (m-2-1) 
(m-1-2) edge node[right]{$p$} (m-2-2) 
(m-2-1) edge node[below]{$h|P_{N+1}(k)$} (m-2-2)
(m-2-1) edge[dotted,inner sep=0.5] node[below right]{$g_{N+1}$} (m-1-2); 
\end{tikzpicture}
\end{equation*}
Once again, one really only has to tiptoe around the exceptional $m$-simplices.

To begin, we may easily extend $g$ along the inclusion
\begin{equation*}
\into{\sigma(0)\cap P_0\cong\Lambda^m_m}{\Delta^m\cong\sigma(0)}
\end{equation*}
(the exceptional $m$-simplex of type (\ref{prp:exceptions}.3)), since the edge $g|\Delta^{\{m-1,m\}}$ is $p$-cartesian.

Now for $0<N<2^m-1$, we have two options for the inclusion
\begin{equation*}
\into{\sigma(N)\cap P_N\cong\Lambda^m_{E(N,0)}}{\Delta^m\cong\sigma(N)}:
\end{equation*}
either it is inner anodyne, in which case it is easy to extend along it, using the fact that $p$ is an inner fibration, or else $N$ is exceptional of type (\ref{prp:exceptions}.4), and hence $N=N_t$ for some integer $0<t<m$, and $E(N_t,0)=\{0,m\}$.

To extend along the inclusion
\begin{equation*}
\into{\sigma(N_t)\cap P_{N_t}\cong\Lambda^m_{\{0,m\}}}{\Delta^m\cong\sigma(N_t)},
\end{equation*}
we factor it as the composite
\begin{equation*}
\Lambda^m_{\{0,m\}}\ \tikz[baseline]\draw[>=stealth,right hook->](0,0.5ex)--(0.5,0.5ex);\ \Lambda^m_{\{0,m\}}\cup^{\Lambda_{0}^{\{0,\dots,m-1\}}}\Delta^{\{0,\dots,m-1\}}\cong\Lambda^{m}_{0}\ \tikz[baseline]\draw[>=stealth,right hook->](0,0.5ex)--(0.5,0.5ex);\ \Delta^m.
\end{equation*}
Extensions along each of these inclusions exists simply because the edge $g|\Delta^{\{0,1\}}$ is $p$-cocartesian.

At the end of this procedure, we are left with an extension $\fromto{P_{2^m-1}}{C}$. To extend over $\sigma(2^m-1)$ (the exceptional $m$-simplex of type (\ref{prp:exceptions}.5)), it suffices just to note that, by assumption, $g|\Delta^{\{0,1\}}$ is $p$-cocartesian, so one may extend over the inclusion
\begin{equation*}
\into{\sigma(2^m-1)\cap P_{2^m-1}\cong\Lambda^m_0}{\Delta^m\cong\sigma(2^m-1)}.
\end{equation*}
The result is the desired extension
\begin{equation*}
\overline{g}\colon\fromto{P_{2^m}=\widetilde{\mathscr{O}}(\Delta^m)^{\op}}{C},
\end{equation*}
which is ambigressive cartesian by construction.
\end{proof}
\end{lem}

The proof of Th. \ref{thm:pain} is complete.

%----------------------------------------------------------------------%

\section{The Burnside Waldhausen bicartesian fibration}\label{sect:BurnWaldbifib} Perhaps the most important Waldhausen bicartesian fibration is the one whose algebraic $K$-theory will be the spectral Burnside Mackey functor $\SS_{(C,C_{\dag},C^{\dag})}$. To describe it, we need some preparatory material.

\begin{ntn} In this subsection, let us fix a disjunctive triple $(C,C_{\dag},C^{\dag})$. 
\end{ntn}

\begin{dfn} Let us say a morphism $\fromto{X}{U}$ of $C$ is a \textbf{\emph{summand inclusion}} if there exists a morphism $\fromto{X'}{U}$ of $C$ that, together with $\fromto{X}{U}$, exhibits $U$ as the coproduct $X\sqcup X'$.

Now if $i\colon\into{X}{U}$ is a summand inclusion, a \textbf{\emph{complement}} of $i$ is a summand inclusion $i'\colon\into{X'}{U}$ such that any square
\begin{equation*}
\begin{tikzpicture} 
\matrix(m)[matrix of math nodes, 
row sep=4ex, column sep=4ex, 
text height=1.5ex, text depth=0.25ex] 
{\varnothing&X\\ 
X'&U\\}; 
\path[>=stealth,right hook->,font=\scriptsize] 
(m-1-1) edge (m-1-2) 
edge (m-2-1) 
(m-1-2) edge node[right]{$i$} (m-2-2) 
(m-2-1) edge node[below]{$i'$} (m-2-2); 
\end{tikzpicture}
\end{equation*}
in which $\varnothing$ is initial in $C$ is a pushout square.
\end{dfn}

\begin{nul} The compatibility of ingressive and egressive morphisms with coproducts implies that summand inclusions are necessarily ingressive and egressive. Note also that the pullback of a summand inclusion along a morphism that is both ingressive and egressive is again a summand inclusion. Furthermore, the pushout of a summand inclusion $i\colon\into{X}{U}$ along any map $f\colon\fromto{X}{Y}$ exists and is again a summand inclusion. Finally, a complement for $i$ is a complement for the pushout $j\colon\into{Y}{V}$ of $i$ along $f$.
\end{nul}

The following lemma will allow us to formulate definitions using complements, as long as we do not use any non-homotopy-invariant constructions.

\begin{lem}\label{lem:compsareunique} Suppose $i\colon\into{X}{U}$ a summand inclusion. If $\mathrm{Compl}(i)\subset C_{/U}$ denotes the full subcategory spanned by the complements of $i$, then the Kan complex $\iota\mathrm{Compl}(i)$ is contractible.
\begin{proof} We show that the diagonal map $\fromto{\iota\mathrm{Compl}(i)}{\iota\mathrm{Compl}(i)\times\iota\mathrm{Compl}(i)}$ is a weak equivalence. To this end, we observe that since $C$ admits all ambigressive pullbacks, it follows that the full subcategory $\mathrm{Sum}(U)\subset C_{/U}$ spanned by the summand inclusions admits all finite products. Consequently, the diagonal functor $\fromto{\mathrm{Sum}(U)}{\mathrm{Sum}(U)\times\mathrm{Sum}(U)}$ admits a right adjoint, which is given informally by the assignment
\begin{equation*}
\goesto{(\into{X'_1}{U},\ \into{X'_2}{U})}{(\into{X'_1\times_UX'_2}{U})}.
\end{equation*}
Our claim is that this right adjoint restricts to a quasi-inverse
\begin{equation*}
\fromto{\iota\mathrm{Compl}(i)\times\iota\mathrm{Compl}(i)}{\iota\mathrm{Compl}(i)}
\end{equation*}
of the diagonal. For this, we must show that if
\begin{equation*}
i'_1\colon\into{X'_1}{U}\textrm{\quad and\quad}i'_2\colon\into{X'_2}{U}
\end{equation*}
are complements of $i$, then the projection maps
\begin{equation*}
\fromto{X'_1\times_UX'_2}{X'_1}\textrm{\quad and\quad}\fromto{X'_1\times_UX'_2}{X'_2}
\end{equation*}
are equivalences, and the morphism $i'_{12}\colon\into{X'_1\times_UX'_2}{U}$ is a complement of $i$. Indeed, the universality of coproducts implies that the projection $\fromto{X'_1\times_UX'_2}{X'_1}$ factors as
\begin{equation*}
\equivto{X'_1\times_UX'_2}{(X'_1\times_UX'_2)\sqcup(X'_1\times_UX)\simeq X'_1\times_U(X'_2\sqcup X)\simeq X'_1\times_UU\simeq X'_1}
\end{equation*}
(and similarly for the projection $\fromto{X'_1\times_UX'_2}{X'_2}$). Now in the cube
\begin{equation*}
\begin{tikzpicture}[cross line/.style={preaction={draw=white, -, 
line width=6pt}}]
\matrix(m)[matrix of math nodes, 
row sep=2ex, column sep=4ex, 
text height=1.5ex, text depth=0.25ex]
{&\varnothing&&\varnothing\\
\varnothing&&X&\\
&X'_1\times_UX'_2&&X'_1\\
X'_2&&U&\\
};
\path[>=stealth,right hook->,font=\scriptsize]
(m-1-2) edge (m-2-1)
edge (m-3-2)
edge (m-1-4)
(m-3-2) edge (m-4-1)
edge (m-3-4)
(m-2-1) edge[cross line] (m-2-3)
edge (m-4-1)
(m-1-4) edge (m-2-3)
edge (m-3-4)
(m-4-1) edge (m-4-3)
(m-3-4) edge (m-4-3)
(m-2-3) edge[cross line] (m-4-3);
\end{tikzpicture}
\end{equation*}
every face is a pullback, and all faces but the top and bottom squares are pushouts, whence $i'_{12}$ is a complement of $i$.
\end{proof}
\end{lem}

\begin{dfn}\label{dfn:RCCdagCdag} We consider the fibration
\begin{equation*}
p\colon\fromto{\Fun(\Delta^2/\Delta^{\{0,2\}},C)\cong\Fun(\Delta^2,C)\times_{\Fun(\Delta^{\{0,2\}},C)}C}{C}.
\end{equation*}
We may think of the objects of the $\infty$-category $\Fun(\Delta^2/\Delta^{\{0,2\}},C)$ as retract diagrams
\begin{equation*}
S_0\to S_1\to S_0;
\end{equation*}
the functor $p$ is given by the assignment
\begin{equation*}
\goesto{[S_0\to S_1\to S_0]}{S_0}.
\end{equation*}
We therefore denote by $C_{S_0/\ /S_0}$ the fiber of $p$ over an object $S_0\in C$.

We consider the full subcategory $\mathscr{R}(C)\subset\Fun(\Delta^2/\Delta^{\{0,2\}},C)$ spanned by those objects $S$ such that the morphism $\fromto{S_0}{S_1}$ is a summand inclusion. We endow $\mathscr{R}(C)$ with the structure of a pair in the following manner. A morphism $\fromto{T}{S}$ will be declared ingressive just in case $\fromto{T_0}{S_0}$ is an equivalence, and $\fromto{T_1}{S_1}$ is a summand inclusion.

Now let $\mathscr{R}(C,C_{\dag},C^{\dag})\subset\mathscr{R}(C)$ be the full subcategory spanned by those objects $S\colon\fromto{\Delta^2/\Delta^{\{0,2\}}}{C}$ such that for any complement $\into{S'_0}{S_1}$ of the summand inclusion $\into{S_0}{S_1}$,
\begin{enumerate}[(\ref{dfn:RCCdagCdag}.1)]
\item the essentially unique morphism $\fromto{S'_0}{\ast}$ to the terminal object of $C$ is egressive, and
\item the composite $S'_0\to S_1\to S_0$ is ingressive.
\end{enumerate}
We endow $\mathscr{R}(C,C_{\dag},C^{\dag})$ with the pair structure induced by $\mathscr{R}(C)$. We will abuse notation by denoting the restriction of the functor $p\colon\fromto{\mathscr{R}(C)}{C}$ to the subcategory $\mathscr{R}(C,C_{\dag},C^{\dag})\subset\mathscr{R}(C)$ again by $p$.
\end{dfn}

We will now show that $p$ is a Waldhausen bicartesian fibration. This claim follows from the following sequence of observations.

\begin{nul} For any object $S_0$ of $C$, the fiber $\mathscr{R}(C)_{S_0}$ can be identified with the full subcategory of $C_{S_{0}/\ /S_{0}}$ spanned by those objects $U$ such that $\into{S}{U}$ is a summand inclusion. A morphism $\fromto{S_1}{S'_1}$ of this $\infty$-category is ingressive just in case it is a summand inclusion. It is an easy consequence of the existence of finite coproducts and the compatibility of the triple structure with these coproducts that the full subcategory $\mathscr{R}(C,C_{\dag},C^{\dag})_{S_0}\subset\mathscr{R}(C)_{S_0}$ is a Waldhausen $\infty$-category.
\end{nul}

\begin{nul} For any ingressive morphism $f\colon\cofto{S_0}{T_0}$ and for any object $S\in\mathscr{R}(C,C_{\dag},C^{\dag})$ over $S_0$, there exists a pushout diagram
\begin{equation*}
\begin{tikzpicture} 
\matrix(m)[matrix of math nodes, 
row sep=4ex, column sep=4ex, 
text height=1.5ex, text depth=0.25ex] 
{S_0&T_0\\ 
S_1&T_1\\
S_0&T_0,\\}; 
\path[>=stealth,->,font=\scriptsize] 
(m-1-1) edge[>->] (m-1-2) 
edge[right hook->] (m-2-1) 
(m-1-2) edge[right hook->] (m-2-2) 
(m-2-1) edge[>->] (m-2-2)
edge (m-3-1)
(m-2-2) edge (m-3-2)
(m-3-1) edge[>->] (m-3-2); 
\end{tikzpicture}
\end{equation*}
hence a $p$-cocartesian edge covering $f$. The compatibility of the triple structure with coproducts ensures that this defines a functor
\begin{equation*}
f_!\colon\fromto{\mathscr{R}(C,C_{\dag},C^{\dag})_{S_0}}{\mathscr{R}(C,C_{\dag},C^{\dag})_{T_0}}.
\end{equation*}
The functors $f_!$ are exact because they are all left adjoints, which preserve any colimits that exist.
\end{nul}

\begin{nul} Dually, for any egressive morphism $f\colon\fibto{S_0}{T_0}$ and for any object $T\in\mathscr{R}(C,C_{\dag},C^{\dag})$ over $T_0$, there exists a pullback diagram
\begin{equation*}
\begin{tikzpicture}
\matrix(m)[matrix of math nodes, 
row sep=4ex, column sep=4ex, 
text height=1.5ex, text depth=0.25ex] 
{S_0&T_0\\ 
S_1&T_1\\
S_0&T_0.\\}; 
\path[>=stealth,->,font=\scriptsize] 
(m-1-1) edge[->>] (m-1-2) 
edge[right hook->] (m-2-1) 
(m-1-2) edge[right hook->] (m-2-2) 
(m-2-1) edge[->>] (m-2-2)
edge (m-3-1)
(m-2-2) edge (m-3-2)
(m-3-1) edge[->>] (m-3-2); 
\end{tikzpicture}
\end{equation*}
We claim that the functor $f^{\star}\colon\fromto{\mathscr{R}(C)_{T_0}}{\mathscr{R}(C)_{S_0}}$ is given informally by the assignment $\goesto{T_1}{T_1\times_{T_0}S_0}$. This follows from the fact that for any cofibration $\cofto{T'_0}{T_0}$ such that the morphism $\fibto{T'_0}{\ast}$ is egressive, the pullback $\cofto{T'_0\times_{T_0}S_0}{S_0}$ is a cofibration and the morphism $\fromto{T'_0\times_{T_0}S_0}{\ast}$ is egressive. This follows from the fact that the pullback of an egressive map is egressive, and the pullback of an ingressive map along an egressive map is ingressive. The universality of finite coproducts in $C$ ensures that the functors $f^{\star}$ preserve finite coproducts; in particular the functors $f^{\star}$ preserve summand inclusions and pushouts along summand inclusions.
\end{nul}

\begin{nul} Now suppose $I$ a finite set, and suppose $\{X_i\ |\ i\in I\}$ a collection of objects of $C$ indexed by the elements of $I$ with coproduct $X$. We claim that the functor
\begin{equation*}
\equivto{C_{/X}}{\prod_{i\in I}C_{/X_i}}
\end{equation*}
induced by the inclusions $\into{X_i}{X}$ induce an equivalence
\begin{equation*}
\equivto{\mathscr{R}(C,C_{\dag},C^{\dag})_{X}}{\prod_{i\in I}\mathscr{R}(C,C_{\dag},C^{\dag})_{X_i}}.
\end{equation*}
One need only note that both this functor, which is given by pullbacks along summand inclusions (which are egressive), and its left adjoint, which is given by coproduct, preserve the desired subcategories and restrict to adjoint equivalences.
\end{nul}

\begin{nul} Finally, the base change condition for $\mathscr{R}(C,C_{\dag},C^{\dag})$ states that for any ambigressive pullback square
\begin{equation*}
\begin{tikzpicture} 
\matrix(m)[matrix of math nodes, 
row sep=4ex, column sep=4ex, 
text height=1.5ex, text depth=0.25ex] 
{S_0&S'_0\\ 
T_0&T'_0,\\}; 
\path[>=stealth,->,font=\scriptsize] 
(m-1-1) edge[>->] node[above]{$i$} (m-1-2) 
edge[->>] node[left]{$q$} (m-2-1) 
(m-1-2) edge[->>] node[right]{$q'$} (m-2-2) 
(m-2-1) edge[>->] node[below]{$j$} (m-2-2); 
\end{tikzpicture}
\end{equation*}
of $C$ and for any object $T'_1$ of $C$ over $T_0$, the base change morphism
\begin{equation*}
\fromto{((T_0\sqcup T'_1)\times_{T_0}S_0)\cup^{S_0}S'_0}{((T_0\sqcup T'_1)\cup^{T_0}T'_0)\times_{T'_0}S'_0}
\end{equation*}
is an equivalence. This follows from the identifications
\begin{eqnarray}
((T_0\sqcup T'_1)\times_{T_0}S_0)\cup^{S_0}S'_0&\simeq&(S_0\sqcup(T'_1\times_{T_0}S_0))\cup^{S_0}S'_0\nonumber\\
&\simeq&S'_0\sqcup(T'_1\times_{T_0}S_0)\nonumber\\
&\simeq&S'_0\sqcup(T'_1\times_{T'_0}S'_0)\nonumber\\
&\simeq&(T'_0\sqcup T'_1)\times_{T'_0}S'_0\nonumber\\
&\simeq&((T_0\sqcup T'_1)\cup^{T_0}T'_0)\times_{T'_0}S'_0,\nonumber
\end{eqnarray}
which all follow from the universality of finite coproducts in $C$ and the equivalence $S_0\simeq T_0\times_{T'_0}S'_0$.
\end{nul}

We thus conclude the following.
\begin{thm} For any disjunctive triple $(C,C_{\dag},C^{\dag})$ that is either left or right complete, the functor
\begin{equation*}
p\colon\fromto{\mathscr{R}(C,C_{\dag},C^{\dag})}{C}
\end{equation*}
is a Waldhausen bicartesian fibration over $(C,C_{\dag},C^{\dag})$.
\end{thm}

Now let's unfurl this Waldhausen bicartesian fibration to obtain a Waldhausen cocartesian fibration
\begin{equation*}
\Upsilon(p)\colon\fromto{\Upsilon(\mathscr{R}(C,C_{\dag},C^{\dag}))}{A^{\eff}(C,C_{\dag},C^{\dag})},
\end{equation*}
whence we obtain a Mackey functor
\begin{equation*}
\mathscr{M}_p\colon\fromto{A^{\eff}(C,C_{\dag},C^{\dag})}{\Wald_{\infty}}.
\end{equation*}

Note that the assignment
\begin{equation*}
\goesto{[1\ \tikz[baseline]\draw[>=stealth,<<-](0,0.5ex)--(0.5,0.5ex);\ U\ \tikz[baseline]\draw[>=stealth,>->](0,0.5ex)--(0.5,0.5ex);\ X]}{[X\ \tikz[baseline]\draw[>=stealth,right hook->](0,0.5ex)--(0.5,0.5ex);\ X\sqcup U\ \tikz[baseline]\draw[>=stealth,->](0,0.5ex)--(0.5,0.5ex);\ X]}
\end{equation*}
defines a functor
\begin{equation*}
\fromto{A^{\eff}(C,C_{\dag},C^{\dag})_{1/}}{\iota_{A^{\eff}(C,C_{\dag},C^{\dag})}\mathscr{R}(C,C_{\dag},C^{\dag})}
\end{equation*}
of left fibrations over $A^{\eff}(C,C_{\dag},C^{\dag})$, and it follows from Lm. \ref{lem:compsareunique} that it is a fiberwise equivalence. Consequently, we deduce that the functor $\iota\circ\mathscr{M}_p$ is naturally equivalent to the functor represented by the terminal object $1$.

Now, almost by definition, the Waldhausen $\infty$-category $\mathscr{M}_p(S)\simeq\mathscr{R}(C,C_{\dag},C^{\dag})_S$ is $\iota$-split. We therefore conclude the following.
\begin{thm}\label{thm:BurnMackWaldbifib} For any disjunctive triple $(C,C_{\dag},C^{\dag})$ that is either left or right complete, the functor
\begin{equation*}
\KK\circ\mathscr{M}_p\colon\fromto{A^{\eff}(C,C_{\dag},C^{\dag})}{\Sp}
\end{equation*}
is the Burnside Mackey functor $\SS_{(C,C_{\dag},C^{\dag})}$ --- i.e., the Mackey functor represented by the terminal object $1\in C$ (Df. \ref{dfn:BurnsideMack}).
\end{thm}

%----------------------------------------------------------------------%

\appendix

\section{Coherent {$n$}-topoi and the Segal--tom Dieck splitting}\label{sect:coherenttopoi} The similarities between the axioms for a disjunctive $\infty$-category and the Giraud axioms for $n$-topoi \cite[Th. 6.1.0.6(3)]{HTT} suggest a deep relationship between the two notions. Here, we briefly describe a way for higher topoi to give rise to a disjunctive $\infty$-category.

\begin{exm} Any subcategory of an $n$-topos ($1\leq n\leq\infty$) that is stable under finite limits and finite coproducts is obviously disjunctive. All our examples in this paper are ultimately of this kind. In fact, one can show that every disjunctive $\infty$-category arises (possibly after a change of universe) in this manner. (If $\tau$ is a strongly inaccessible uncountable cardinal, then any $\tau$-small disjunctive $\infty$-category $C$ can be embedded in the full subcategory $\mathscr{X}\subset\Fun(C^{\op},\Kan(\tau))$ spanned by the functors that preserve products. This is an accessible localization of the $\infty$-category $\Fun(C^{\op},\Kan(\tau))$, and one can show that the localization functor $\fromto{\Fun(C^{\op},\Kan(\tau))}{\mathscr{X}}$ is left exact, whence $\mathscr{X}$ is an $\infty$-topos.)
\end{exm}

\begin{exm}\label{exm:cohntopos} More particularly, the full subcategory $\mathscr{X}^{\coh}$ of a coherent $\infty$-topos $\mathscr{X}$ spanned by the coherent objects \cite[Df. 3.12]{DAGVII} is thus a disjunctive $\infty$-category. Furthermore, for any natural number $n$, the full subcategory $\tau_{\leq n}\mathscr{X}^{\coh}$ spanned by the $n$-truncated objects is a full subcategory of an $(n+1)$-topos that is closed under finite coproducts and finite limits; hence it too is a disjunctive $\infty$-category.

For $1\leq n<\infty$, let us say that an object $U$ of an $n$-topos $\mathscr{X}$ is \textbf{\emph{coherent}} if for any $n$-localic $\infty$-topos $\mathscr{Y}$ and any equivalence $\phi\colon\equivto{\mathscr{X}}{\tau_{\leq n-1}\mathscr{Y}}$, the object $\phi(U)$ is coherent. If $\mathscr{Y}$ is coherent, then we will say that $\mathscr{X}$ is \textbf{\emph{coherent}}, and the full subcategory $\mathscr{X}^{\coh}\subset\mathscr{X}$ by the coherent objects is a disjunctive $\infty$-category.
\end{exm}

\begin{exm} In particular, the $\infty$-category $\Kan^{\mathit{coh}}$ of Kan simplicial sets all of whose homotopy groups are finite is a disjunctive $\infty$-category, and for any $n\geq 0$, the truncation $\tau_{\leq n}\Kan^{\mathit{coh}}$ (whose objects may be called \textbf{\emph{finite $n$-groupoids}}) is a disjunctive $\infty$-category.

The effective Burnside $\infty$-category $A^{\eff}(\tau_{\leq n}\Kan^{\coh})$ is an $(n+1)$-category in the sense of \cite[Df. 2.3.4.1]{HTT}. The homotopy category of $A^{\eff}(\tau_{\leq 1}\Kan^{\coh})=A^{\eff}(\FF)$ (where $\FF$ is the nerve of the ordinary category of finite sets) is the ordinary effective Burnside category for the trivial group.
\end{exm}

The following proposition, whose proof we leave to the reader, can be summarized by saying that the $\infty$-category $A^{\eff}(\FF)$ as the free $\infty$-category with direct sums generated by a single object (the terminal object $1$ of $\FF$).
\begin{prp} For any $\infty$-category $E$ that admits direct sums, evaluation at the terminal object $1$ induces an equivalence of $\infty$-categories
\begin{equation*}
\equivto{\Mack(\FF,E)}{E}.
\end{equation*}
\end{prp}

Let us turn now to a general version of the Segal--tom Dieck splitting theorem.
\begin{ntn} To this end, note that if $\mathscr{X}$ is a coherent $n$-topos ($1\leq n\leq\infty$), then there is a functor $\delta\colon\fromto{\FF}{\mathscr{X}^{\mathit{coh}}}$ informally given by the assignment
\begin{equation*}
\goesto{S}{S\otimes 1\simeq\coprod_{s\in S}1}.
\end{equation*}
See \cite[\S 4.4.4]{HTT}.
\end{ntn}
 
\begin{dfn} A coherent $n$-topos $\mathscr{X}$ will be said to be \textbf{\emph{locally connected}} if the functor $\delta$ admits a left adjoint $\pi$.
\end{dfn}

The idea here is of course that for any coherent object $X$, the unit morphism $\fromto{X}{\delta\pi X}$ will decompose $X$ into finitely many summands, each of which will be ``connected,'' in a unique manner. Indeed, for any object $X$, one may exhibit $X$ as a canonical coproduct
\begin{equation*}
X\simeq\coprod_{\alpha\in\pi X}(X\times_{\delta\pi X}\delta\{\alpha\}).
\end{equation*}
These objects $X\times_{\delta\pi X}\delta\{\alpha\}$ are now connected in the following sense.

\begin{dfn}\label{dfn:connectedobject} Suppose $\mathscr{X}$ a coherent $n$-topos. A coherent object $X$ of $\mathscr{X}$ will be said to be \textbf{\emph{connected}} if the functor $\fromto{\mathscr{X}^{\mathit{coh}}}{\Kan}$ it corepresents preserves finite coproducts. Equivalently, $X$ is connected just in case $\pi X$ is a one-point set. Denote by $\mathscr{X}^{\mathit{conn}}\subset\mathscr{X}^{\mathit{coh}}$ the full subcategory spanned by the connected objects.
\end{dfn}

\begin{wrn} It is not necessarily the case that the terminal object $\ast$ of a locally connected coherent $n$-topos $\mathscr{X}$ is connected. If it is, then $\mathscr{X}$ is said to be \textbf{\emph{connected}}, and in this case $\delta$ is fully faithful.
\end{wrn}

It turns out that the spectral Burnside ring of a locally connected coherent $n$-topos can be identified with a suspension spectrum.

\begin{thm}[Segal--tom Dieck, \protect{\cite{MR0436177}}]\label{thm:tomDieck} Suppose $\mathscr{X}$ a locally connected coherent $n$-topos. Then there is a natural equivalence
\begin{equation*}
\SS_{\mathscr{X}^{\mathit{coh}}}(1)\simeq\Sigma^{\infty}_{+}\iota\mathscr{X}^{\mathit{conn}}\simeq\bigvee_{X}\Sigma_+^{\infty}B\Aut(X),
\end{equation*}
where the wedge is taken over all equivalence classes of connected objects $X$, and $\Aut(X)$ denotes the space of auto-equivalences of $X$.
\begin{proof} It follows from Th. \ref{thm:BurnMackWaldbifib} that $\SS_{\mathscr{X}^{\mathit{coh}}}(1)$ can be identified with the group completion of the object $\iota\mathscr{X}^{\mathit{coh}}\in\mathrm{CAlg}(\Kan)$, where the commutative algebra structure is given by coproduct. It therefore suffices to show that $\iota\mathscr{X}^{\mathit{coh}}$ is the free commutative algebra generated by the space $\iota\mathscr{X}^{\mathit{conn}}$ in the sense of \cite[Ex. 3.1.3.12]{HA}.

For this, note that the (homotopy) fiber of the map $\iota(\pi)\colon\fromto{\iota\mathscr{X}^{\mathit{coh}}}{\iota N\FF}$ over a finite set $I$ of cardinality $n$ may be described as the space of pairs $(X,f)$ consisting of an object $X\in\iota\mathscr{X}^{\mathit{coh}}$ and an isomorphism $\equivto{\pi(X)}{I}$, which can in turn be identified with the product $(\iota\mathscr{X}^{\mathit{conn}})^n$. We therefore obtain an identification
\begin{equation*}
\mathrm{Sym}^n(\iota\mathscr{X}^{\mathit{conn}})\simeq B\Sigma_n\times^h_{\iota N\FF}\iota\mathscr{X}^{\mathit{coh}}.
\end{equation*}

Since one has $\iota N\FF\simeq\coprod_{n\geq 0}B\Sigma_n$, we find that
\begin{equation*}
\iota\mathscr{X}^{\mathit{coh}}\simeq\coprod_{n\geq 0}\mathrm{Sym}^n(\iota\mathscr{X}^{\mathit{conn}}).
\end{equation*}
The proof is thus complete by \cite[Ex. 3.1.3.11]{HA}.
\end{proof}
\end{thm}

%----------------------------------------------------------------------%

\section{Equivariant spectra for a profinite group}\label{exm:profgroups} Certain topological groups give a subexample of Ex. \ref{sect:coherenttopoi}.

\begin{dfn} A topological group is \textbf{\emph{coherent}} if for every open subgroup $H$ of $G$, there exist only finitely many subsets of the form $HgH$ for $g\in G$.
\end{dfn}

\begin{nul} Suppose $G$ a coherent topological group. Suppose $B_{G}$ the classifying $1$-topos of $G$; that is, $B_{G}$ is the nerve of the ordinary category of sets equipped with a continuous action of $G$. Then $B_{G}$ is the nerve of a coherent topos \cite[\S D3.4]{MR2063092}, and the full subcategory $B_G^{\mathit{fin}}$ spanned by the coherent objects is a disjunctive $\infty$-category.
\end{nul}

\begin{nul} A pro-discrete group $G$ is coherent just in case it is profinite. In this case, the coherent objects of $B_{G}$ are simply those continuous $G$-sets with only finitely many orbits. Hence the full subcategory $B_{G}^{\mathit{fin}}\subset B_{G}$ spanned by finite $G$-sets with open stabilizers is a disjunctive $\infty$-category. The effective Burnside $\infty$-category of $B_{G}^{\mathit{fin}}$ will be denoted, abusively, $A^{\eff}(G)$.
\end{nul}

\begin{nul} For any finite group $G$, the $\infty$-category $A^{\eff}(G)$ is in fact a $2$-category whose homotopy category $hA^{\eff}(G)$ is the ordinary effective Burnside category for $G$; the ordinary Burnside category for $G$ is the local group completion (obtained by forming the Grothendieck group of each of the Hom-sets under direct sum).
\end{nul}

\begin{ntn} When $C=NB_G^{\mathit{fin}}$ for some profinite group $G$ and $E$ is some additive $\infty$-category, we will write
\begin{equation*}
\Mack_G(E)\coloneq\Mack(NB_G^{\mathit{fin}},E).
\end{equation*}
\end{ntn}

\begin{exm} When $G$ is finite and $E=N\mathbf{Ab}$ is the nerve of the ordinary category of abelian groups, the $\infty$-category $\Fun(A^{\eff}(NB_G^{\mathit{fin}}),N\mathbf{Ab})$ is naturally equivalent to the nerve of the ordinary category of functors $\Fun(hA^{\eff}(NB_G^{\mathit{fin}}),\mathbf{Ab})$. Using the fact that the homotopy category $hA^{\eff}(NB_G^{\mathit{fin}})$ is the ordinary effective Burnside category \eqref{nul:ordinaryBurnside}, we conclude that the full subcategory $\Mack_G(N\mathbf{Ab})$ is equivalent to the full subcategory of $\Fun(hA^{\eff}(NB_G^{\mathit{fin}}),\mathbf{Ab})$ spanned by Mackey functors in the classical sense.

When $G$ is finite and $E=\Sp$ is the $\infty$-category of spectra, work of Guillou and May \cite{guillou2} show that this $\infty$-category is equivalent to the underlying $\infty$-category of the relative category of genuine $G$-spectra in the sense of Lewis--May--Steinberger \cite{MR866482}, Mandell--May \cite{MR1922205}, and Hill--Hopkins--Ravenel \cite{MR2757358,MR2906370,HHRarxiv}. For general profinite groups $G$, we call $\Mack_G(\Sp)$ the \textbf{\emph{$\infty$-category of $G$-equivariant spectra}}.
\end{exm}

\begin{nul} Suppose $G$ a profinite group, and suppose $H$ a closed normal subgroup of $G$. Then the natural functor $\varphi\colon\fromto{G}{G/H}$ induces a morphism of topoi
\begin{equation*}
\adjunct{\varphi^{\star}}{B_{G/H}}{B_{G}}{\varphi_{\star}}.
\end{equation*}
Both functors preserve coherent objects, finite coproducts and pullbacks. Hence we obtain functors
\begin{equation*}
\varphi^{\star}\colon\fromto{NB_{G/H}^{\mathit{fin}}}{NB_{G}^{\mathit{fin}}}\textrm{\quad and\quad}\varphi_{\star}\colon\fromto{NB_{G}^{\mathit{fin}}}{NB_{G/H}^{\mathit{fin}}}.
\end{equation*}
For any presentable additive $\infty$-category $E$, we obtain adjunctions
\begin{equation*}
\adjunct{\varphi^{\star}_{!}}{\Mack_{G/H}(E)}{\Mack_{G}(E)}{\varphi^{\star\star}}
\end{equation*}
and
\begin{equation*}
\adjunct{\varphi_{\star!}}{\Mack_{G}(E)}{\Mack_{G/H}(E)}{\varphi_{\star}^{\star}}.
\end{equation*}

We write $\Psi^H\coloneq\varphi^{\star\star}$ and $\Phi^H\coloneq\varphi_{\star!}$. When $G$ is finite, one may use the equivalence of Guillou--May \cite{guillou2} to interpret $\Psi^{H}$ and $\Phi^{H}$ as functors on the $\infty$-categories of $G$- and $G/H$-equivariant spectra. One may show that these functors agree (up to equivalence) with the \emph{Lewis--May fixed points} $\Psi^H$ and the \emph{geometric fixed points} $\Phi^H$ constructed by Mandell--May \cite{MR1922205}. 
\end{nul}

For any profinite group $G$, the $1$-topos $NB_{G}^{\mathit{fin}}$ is locally connected (connected, in fact). This fact permits us to use our Segal--tom Dieck Theorem \ref{thm:tomDieck} to compute the value of the Burnside spectral Mackey functor
\begin{equation*}
\SS_G\coloneq\SS_{NB_{G}^{\mathit{fin}}}
\end{equation*}
on the terminal $G$-set $[G/G]$. Indeed, the connected objects of $NB_{G}^{\mathit{fin}}$ are precisely the finite $G$-sets with open stabilizers that are transitive. Up to equivalence, these are classified by conjugacy classes of open subgroups of $G$. The space of autoequivalences of the transitive $G$-set $[G/H]$ is equivalent to the quotient $N_{G}H/H$, whence we obtain the traditional Segal--tom Dieck splitting \cite{MR0436177}, now for profinite groups:
\begin{prp}\label{prp:SegtomDieckprofinite} For any profinite group $G$, one has
\begin{equation*}
\SS_G([G/G])\simeq\bigvee_{H}\Sigma^{\infty}_{+}B(N_{G}H/H),
\end{equation*}
where the wedge is indexed by conjugacy classes of open subgroups $H\leq G$.
\end{prp}

%----------------------------------------------------------------------%

\section{$A$-theory and upside-down $A$-theory of $\infty$-topoi}\label{sect:Athy} In this section, we introduce two dual disjunctive triple structures on an $\infty$-topos, where the ingressive or egressive morphisms are defined by means of a finiteness condition. We use these structures to construct both the $A$-theory and the $\upsidedown{A}$-theory of $\infty$-topoi together with all their functorialities.

To begin, if we consider an $\infty$-topos as a disjunctive $\infty$-category, then there is a natural Waldhausen bicartesian fibration that lies over it. Let us investigate this now.

The following is an easy analogue of \cite[Lm. 6.1.1.1]{HTT}.
\begin{lem}\label{lem:Eoverandunder} If $C$ is an $\infty$-category that admits both pullbacks and pushouts, then the functor
\begin{equation*}
\fromto{\Fun(\Delta^2,C)}{\Fun(\Delta^{\{0,2\}},C)}
\end{equation*}
is a both a cartesian fibration and a cocartesian fibration.
\end{lem}

\begin{ntn} For the remainder of this section, suppose $\mathscr{X}$ an $\infty$-topos. We consider the fibration
\begin{equation*}
p\colon\fromto{\Fun(\Delta^2/\Delta^{\{0,2\}},\mathscr{X})\cong\Fun(\Delta^2,\mathscr{X})\times_{\Fun(\Delta^{\{0,2\}},\mathscr{X})}\mathscr{X}}{\mathscr{X}},
\end{equation*}
which is both cartesian and cocartesian. We may think of the objects of the $\infty$-category $\Fun(\Delta^2/\Delta^{\{0,2\}},\mathscr{X})$ as retract diagrams
\begin{equation*}
X\to X'\to X;
\end{equation*}
the functor $p$ is given by the assignment
\begin{equation*}
\goesto{[X\to X'\to X]}{X}.
\end{equation*}
We therefore denote $\mathscr{X}_{X/\ /X}$ the fiber of $p$ over an object $X$ of $\mathscr{X}$. For any morphism $f\colon\fromto{X}{Y}$, the corresponding functors
\begin{equation*}
f_{!}\colon\fromto{\mathscr{X}_{X/\ /X}}{\mathscr{X}_{Y/\ /Y}}\textrm{\quad and\quad}f^{\star}\colon\fromto{\mathscr{X}_{Y/\ /Y}}{\mathscr{X}_{X/\ /X}}
\end{equation*}
are given informally by the assignments $\goesto{X'}{Y\cup^XX'}$ and $\goesto{Y'}{Y'\times_YX}$.
\end{ntn}

\begin{prp} Endow the $\infty$-category $\Fun(\Delta^2/\Delta^{\{0,2\}},\mathscr{X})$ with the pair structure in which a morphism $f$ is ingressive just in case $p(f)$ is an equivalence. Then $p\colon\fromto{\Fun(\Delta^2/\Delta^{\{0,2\}},\mathscr{X})}{\mathscr{X}}$ is a Waldhausen bicartesian fibration.
\begin{proof} The only nontrivial points are the following. For any morphism $f\colon\fromto{X}{Y}$, the functor $f^{\star}$ preserves pushout squares by the universality of colimits in $\mathscr{X}$. The admissibility of $p$ also follows from the universality of colimits in $\mathscr{X}$. Finally, if $I$ is a finite set and $\{X_i\}_{i\in I}$ a family of objects of $\mathscr{X}$ with coproduct $X$, then pullback along the various inclusions $\into{X_i}{X}$ defines an equivalence of $\infty$-categories
\begin{equation*}
\equivto{\mathscr{X}_{X/\ /X}}{\prod_{i\in I}\mathscr{X}_{X_i/\ /X_i}}.\qedhere
\end{equation*}
\end{proof}
\end{prp}

This Waldhausen bicartesian fibration is not so interesting from the point of view of algebraic $K$-theory, as each $\infty$-category $\mathscr{X}_{X/\ /X}$ has vanishing $K$-theory. To make it more interesting, we must restrict attention to objects with a finiteness condition. We thus turn to the study of finiteness conditions on objects of the $\infty$-topos $\mathscr{X}$. We begin with the following result.
\begin{prp}\label{prp:compactovercat} Suppose $X\in\mathscr{X}$ an object. An object $X'\in\mathscr{X}_{/X}$ is compact if and only if it is compact as an object of $\mathscr{X}$.
\begin{proof} It is easy to see that if $X'$ is compact in $\mathscr{X}$, then it is compact in $\mathscr{X}_{/X}$. Conversely, the forgetful functor $\fromto{\mathscr{X}_{/X}}{\mathscr{X}}$ admits a right adjoint that, thanks to the universality of colimits, preserves colimits. It follows from \cite[Pr. 5.5.7.2(1)]{HTT} that if $X'$ is compact as an object of $\mathscr{X}_{/X}$, then it is compact as an object of $\mathscr{X}$.
\end{proof}
\end{prp}

\noindent The following is now immediate.
\begin{cor}\label{cor:compactgenovercat} If the $\infty$-topos $\mathscr{X}$ is compactly generated, then for any object $X\in\mathscr{X}$, the $\infty$-topos $\mathscr{X}_{/X}$ is compactly generated as well.
\end{cor}

In the $\infty$-category of spaces, those maps whose fibers are (homotopy) retracts of finite CW complexes play a special role. In a more general $\infty$-topos, this role is played by the \emph{relatively compact} morphisms, which we define now.
\begin{dfn}[Lurie, \protect{\cite[Df. 6.1.6.4]{HTT}}] A morphism $\fromto{Y}{X}$ of $\mathscr{X}$ is said to be \textbf{\emph{relatively compact}} just in case, for any compact object $S$ of $\mathscr{X}$ and any morphism $\eta\colon\fromto{S}{X}$ and any pullback square
\begin{equation*}
\begin{tikzpicture} 
\matrix(m)[matrix of math nodes, 
row sep=4ex, column sep=4ex, 
text height=1.5ex, text depth=0.25ex] 
{Y_{\eta}&Y\\ 
S&X,\\}; 
\path[>=stealth,->,font=\scriptsize] 
(m-1-1) edge (m-1-2) 
edge (m-2-1) 
(m-1-2) edge (m-2-2) 
(m-2-1) edge node[below]{$\eta$} (m-2-2); 
\end{tikzpicture}
\end{equation*}
the object $Y_{\eta}$ is compact.
\end{dfn}

\begin{prp}\label{nul:Xcompactfibers} Suppose $\mathscr{G}\subset\mathscr{X}$ a full subcategory whose objects are compact that generates $\mathscr{X}$ under colimits (so that $\mathscr{X}$ is compactly generated). Then a morphism $\fromto{Y}{X}$ of $\mathscr{X}$ is relatively compact just in case, for any object $T\in\mathscr{G}$ and any morphism $\xi\colon\fromto{T}{X}$ and any pullback square
\begin{equation*}
\begin{tikzpicture} 
\matrix(m)[matrix of math nodes, 
row sep=4ex, column sep=4ex, 
text height=1.5ex, text depth=0.25ex] 
{Y_{\xi}&Y\\ 
T&X,\\}; 
\path[>=stealth,->,font=\scriptsize] 
(m-1-1) edge (m-1-2) 
edge (m-2-1) 
(m-1-2) edge (m-2-2) 
(m-2-1) edge node[below]{$\xi$} (m-2-2); 
\end{tikzpicture}
\end{equation*}
the object $Y_{\xi}$ is compact.
\begin{proof} To prove this claim, it suffices to show that if
\begin{enumerate}[(\ref{nul:Xcompactfibers}.1)]
\item $S\in\mathscr{X}^{\omega}$,
\item $Z\in\mathscr{X}$, and
\item $\fromto{Z}{S}$ is a morphism such that for any pullback square
\begin{equation*}
\begin{tikzpicture} 
\matrix(m)[matrix of math nodes, 
row sep=4ex, column sep=4ex, 
text height=1.5ex, text depth=0.25ex] 
{Z_{\xi}&Z\\ 
T&S,\\}; 
\path[>=stealth,->,font=\scriptsize] 
(m-1-1) edge (m-1-2) 
edge (m-2-1) 
(m-1-2) edge (m-2-2) 
(m-2-1) edge node[below]{$\xi$} (m-2-2); 
\end{tikzpicture}
\end{equation*}
of $\mathscr{X}$ in which $T\in\mathscr{G}$, the object $Z_{\xi}$ is compact,
\end{enumerate}
then the object $Z$ is also compact. To see this, we first argue that $\mathscr{X}^{\omega}$ is generated under finite colimits and retracts by $\mathscr{G}$; indeed, if $\mathscr{X}^f\subset\mathscr{X}^{\omega}$ denotes the full subcategory generated by $\mathscr{G}$ under finite colimits, then \cite[Pr. 5.3.5.11]{HTT} implies that the colimit preserving functor $\fromto{\Ind(\mathscr{X}^f)}{\Ind(\mathscr{X}^{\omega})=\mathscr{X}}$ corresponding to the inclusion $\into{\mathscr{X}^f}{\mathscr{X}^{\omega}}$ is an equivalence, whence $\mathscr{X}^{\omega}$ is the idempotent completion of $\mathscr{X}^f$ thanks to \cite[Pr. 5.5.7.8]{HTT}. We therefore may write $S$ as a retract of a finite colimit of objects of $\mathscr{G}$, and employing the universality of colimits, we obtain a presentation of $Z$ as a retract of a finite colimit of compact objects.
\end{proof}
\end{prp}

\begin{exm} In particular, we deduce that the $\infty$-category $\Kan$ admits a disjunctive triple structure in which every morphism is ingressive, and a morphism is egressive just in case its fibers are compact (or, in the parlance of \cite[p. 3]{MR1982793}, \emph{homotopy finitely dominated}).
\end{exm}

\begin{prp} If $\mathscr{X}_{\mathrm{rc}}\subset\mathscr{X}$ denotes the subcategory of relatively compact morphisms, then the triples $(\mathscr{X},\mathscr{X},\mathscr{X}_{\mathrm{rc}})$ and $(\mathscr{X},\mathscr{X}_{\mathrm{rc}},\mathscr{X})$ are disjunctive.
\begin{proof} It is obvious that relatively compact morphisms are closed under pullback. The universality of colimits and the fact that compact objects of $\mathscr{X}$ are closed under finite coproducts implies that the class of relatively compact morphisms is compatible with coproducts. The other axioms all follow directly from the universality of colimits and the disjointness of finite coproducts in $\mathscr{X}$.
\end{proof}
\end{prp}

\begin{ntn}\label{ntn:IXandJX} Now let us restrict the Waldhausen bicartesian fibration
\begin{equation*}
p\colon\fromto{\Fun(\Delta^2/\Delta^{\{0,2\}},\mathscr{X})}{\mathscr{X}}.
\end{equation*}
Let us write $\mathscr{I}(\mathscr{X})$ for the full subcategory of $\Fun(\Delta^2/\Delta^{\{0,2\}},\mathscr{X})$ spanned by those retract diagrams
\begin{equation*}
X\to X'\to X
\end{equation*}
such that the object $X'\in\mathscr{X}_{X/\ /X}$ is compact. Dually, let us write $\mathscr{J}(\mathscr{X})$ for the full subcategory of $\Fun(\Delta^2/\Delta^{\{0,2\}},\mathscr{X})$ spanned by those retract diagrams
\begin{equation*}
X\to X'\to X
\end{equation*}
such that the morphism $\fromto{X'}{X}$ is relatively compact.
\end{ntn}

\begin{thm} Assume that $\mathscr{X}$ is compactly generated, and assume that the terminal object $1\in\mathscr{X}$ is compact. The functor $p$ restricts to Waldhausen bicartesian fibrations
\begin{equation*}
p_{\!\mathscr{I}(\mathscr{X})}\colon\fromto{\mathscr{I}(\mathscr{X})}{\mathscr{X}}
\end{equation*}
over $(\mathscr{X},\mathscr{X},\mathscr{X}_{\mathrm{rc}})$ and
\begin{equation*}
p_{\!\mathscr{J}(\mathscr{X})}\colon\fromto{\mathscr{J}(\mathscr{X})}{\mathscr{X}}
\end{equation*}
over $(\mathscr{X},\mathscr{X}_{\mathrm{rc}},\mathscr{X})$.
\begin{proof} The only points left to be shown are the following for a morphism $f\colon\fromto{X}{Y}$ of $\mathscr{X}$.
\begin{enumerate}[(A)]
\item If $X'\in\mathscr{X}_{X/\ /X}$ is compact, then $X'\cup^XY\in\mathscr{X}_{Y/\ /Y}$ is also compact.
\item If $\fromto{Y'}{Y}$ is relatively compact, then the pullback $\fromto{Y'\times_YX}{X}$ is also relatively compact.
\item If $f$ is relatively compact and $Y'\in\mathscr{X}_{Y/\ /Y}$ is compact, then $Y'\times_YX\in\mathscr{X}_{X/\ /X}$ is also compact.
\item If $f$ is relatively compact and $\fromto{X'}{X}$ is relatively compact, then the pushout $\fromto{X'\cup^XY}{Y}$ is also relatively compact.
\end{enumerate}

To prove (1), note that the functor $\fromto{\mathscr{X}_{X/\ /X}}{\mathscr{X}_{Y/\ /Y}}$ given by the assignment $\goesto{X'}{X'\cup^XY}$ can be identified with the tensor product
\begin{equation*}
u\otimes\id\colon\fromto{\mathscr{X}_{/X}\otimes\Kan_{\ast}}{\mathscr{X}_{/Y}\otimes\Kan_{\ast}}
\end{equation*}
of presentable $\infty$-categories (in the sense of \cite[\S 6.3.1]{HA}) of the forgetful functor $u\colon\fromto{\mathscr{X}_{/X}}{\mathscr{X}_{/Y}}$ with the identity functor \cite[Pr. 6.3.2.11]{HA}. By Pr. \ref{prp:compactovercat} and Cor. \ref{cor:compactgenovercat}, $u$ is an $\omega$-good functor between compactly generated $\infty$-categories in the sense of \cite[Nt. 6.3.7.8]{HA}. Hence by \cite[Lm. 6.3.7.11]{HA}, $u\otimes\id$ is $\omega$-good as well. In particular, it preserves  compact objects.

Assertion (2) is clear.

To prove (3), note that the functor $v\colon\fromto{\mathscr{X}_{/Y}}{\mathscr{X}_{/X}}$ given by the assignment $\goesto{Y'}{Y'\times_YX}$ preserves colimits by the universality of colimits in $\mathscr{X}$, and it preserves compact objects thanks to Pr. \ref{prp:compactovercat}; we thus conclude that $v$ is $\omega$-good. Once again we may identify the functor $\fromto{\mathscr{X}_{Y/\ /Y}}{\mathscr{X}_{X/\ /X}}$ given by the assignment $\goesto{Y'}{Y'\times_YX}$ with $v\otimes\id$, and once again we may appeal to \cite[Lm. 6.3.7.11]{HA} to conclude that $v\otimes\id$ is $\omega$-good.

Finally, to prove (4), assume that $S\in\mathscr{X}$ is a compact object and $\fromto{S}{Y}$ a morphism. The universality of colimits implies that
\begin{equation*}
(X'\cup^XY)\times_YS\simeq(X'\times_YS)\cup^{(X\times_YS)}S.
\end{equation*}
Since both $\fromto{X'}{X}$ and $\fromto{X}{Y}$ are relatively compact, it follows that both $X'\times_YS$ and $X\times_YS$ are compact. Since compact objects are closed under finite colimits, it follows that $(X'\cup^XY)\times_YS$ is compact as well, and this completes the proof that $\fromto{X'\cup^XY}{Y}$ is relatively compact.
\end{proof}
\end{thm}

\begin{cor} Assume that $\mathscr{X}$ is compactly generated, and assume that the terminal object $1\in\mathscr{X}$ is compact. The unfurlings
\begin{equation*}
\fromto{\Upsilon(p_{\!\mathscr{I}(\mathscr{X})})}{A^{\eff}(\mathscr{X},\mathscr{X},\mathscr{X}_{\mathrm{rc}})}\textrm{\quad and\quad}\fromto{\Upsilon(p_{\!\mathscr{J}(\mathscr{X})})}{A^{\eff}(\mathscr{X},\mathscr{X}_{\mathrm{rc}},\mathscr{X})}
\end{equation*}
classify Mackey functors
\begin{equation*}
\mathscr{M}_{\!\mathscr{I}(\mathscr{X})}\colon\fromto{A^{\eff}(\mathscr{X},\mathscr{X},\mathscr{X}_{\mathrm{rc}})}{\Wald_{\infty}}
\end{equation*}
and
\begin{equation*}
\mathscr{M}_{\!\mathscr{J}(\mathscr{X})}\colon\fromto{A^{\eff}(\mathscr{X},\mathscr{X}_{\mathrm{rc}},\mathscr{X})}{\Wald_{\infty}}.
\end{equation*}
\end{cor}

\begin{ntn}\label{ntn:AandupsidedownA} Assume that $\mathscr{X}$ is compactly generated, and assume that the terminal object $1\in\mathscr{X}$ is compact. Write $\AA$ and $\upsidedown{\AA}$ for the composite Mackey functors
\begin{equation*}
\AA_{\mathscr{X}}\coloneq\KK\circ\mathscr{M}_{\!\mathscr{I}(\mathscr{X})}\textrm{\quad and\quad}\upsidedown{\AA}_{\!\mathscr{X}}\coloneq\KK\circ\mathscr{M}_{\!\mathscr{J}(\mathscr{X})}.
\end{equation*}
For any object $X\in\mathscr{X}$, the spectrum $\AA_{\mathscr{X}}(X)$ is the algebraic $K$-theory of the full subcategory $\mathscr{I}(\mathscr{X})_X\subset\mathscr{X}_{X/\ /X}$ spanned by the compact objects, and $\upsidedown{\AA}_{\!\mathscr{X}}(X)$ is the algebraic $K$-theory of the full subcategory $\mathscr{J}(\mathscr{X})_X\subset\mathscr{X}_{X/\ /X}$ spanned by those retract diagrams $X\to X'\to X$ such that $\fromto{X'}{X}$ is relatively compact.
\end{ntn}

\begin{nul}\label{nul:Aassembly} Assume that $\mathscr{X}$ is compactly generated, and assume that the terminal object $1\in\mathscr{X}$ is compact. We have equivalences
\begin{equation*}
\Map_{\Mack(\mathscr{X},\mathscr{X},\mathscr{X}_{\mathrm{rc}};\Sp)}(\SS_{(\mathscr{X},\mathscr{X},\mathscr{X}_{\mathrm{rc}})},\AA_{\mathscr{X}})\simeq\AA_{\mathscr{X}}(1)
\end{equation*}
and
\begin{equation*}
\Map_{\Mack(\mathscr{X},\mathscr{X}_{\mathrm{rc}},\mathscr{X};\Sp)}(\SS_{(\mathscr{X},\mathscr{X}_{\mathrm{rc}},\mathscr{X})},\upsidedown{\AA}_{\!\mathscr{X}})\simeq\upsidedown{\AA}_{\!\mathscr{X}}(1),
\end{equation*}
where $1$ denotes the terminal object of $\mathscr{X}$. The object $1\sqcup 1\in\mathscr{X}_{1/\ /1}$ lies in both $\mathscr{I}(\mathscr{X})_1$ and $\mathscr{J}(\mathscr{X})_1$, whence its classes in the corresponding $K$-theories specify  morphisms of Mackey functors
\begin{equation*}
\fromto{\SS_{(\mathscr{X},\mathscr{X},\mathscr{X}_{\mathrm{rc}})}}{\AA_{\mathscr{X}}}\textrm{\quad and\quad}\fromto{\SS_{(\mathscr{X},\mathscr{X}_{\mathrm{rc}},\mathscr{X})}}{\upsidedown{\AA}_{\!\mathscr{X}}}.
\end{equation*}
We thus obtain, for any object $X\in\mathscr{X}$, assembly morphisms
\begin{equation*}
\fromto{\SS_{(\mathscr{X},\mathscr{X},\mathscr{X}_{\mathrm{rc}})}(X)}{\AA_{\mathscr{X}}(X)}\textrm{\quad and\quad}\fromto{\SS_{(\mathscr{X},\mathscr{X}_{\mathrm{rc}},\mathscr{X})}(X)}{\upsidedown{\AA}_{\!\mathscr{X}}(X)}.
\end{equation*}
We in turn obtain, for any three objects $U,V,X\in\mathscr{X}$ such that $\fromto{U}{1}$ is relatively compact, assembly morphisms
\begin{equation*}
\fromto{\Sigma^{\infty}_{+}\Map_{\mathscr{X}}(U,X)}{\AA_{\mathscr{X}}(X)}\textrm{\quad and\quad}\fromto{\Sigma^{\infty}_{+}\Map_{\mathscr{X}_{\mathrm{rc}}}(V,X)}{\upsidedown{\AA}_{\!\mathscr{X}}(X)}.
\end{equation*}

As in Df. \ref{dfn:assembly}, we obtain, for any object $X\in\mathscr{X}$, morphisms of Mackey functors
\begin{equation*}
\fromto{\SS^X_{(\mathscr{X},\mathscr{X},\mathscr{X}_{\mathrm{rc}})}}{F(\AA_{\mathscr{X}}(X),\AA_{\mathscr{X}})}\textrm{\quad and\quad}\fromto{\SS^X_{(\mathscr{X},\mathscr{X}_{\mathrm{rc}},\mathscr{X})}}{F(\upsidedown{\AA}_{\!\mathscr{X}}(X),\upsidedown{\AA}_{\!\mathscr{X}})},
\end{equation*}
which induce, for any object $Y\in\mathscr{X}$, assembly morphisms
\begin{equation*}
\fromto{\SS^X_{(\mathscr{X},\mathscr{X},\mathscr{X}_{\mathrm{rc}})}(Y)\wedge\AA_{\mathscr{X}}(X)}{\AA_{\mathscr{X}}(Y)}
\end{equation*}
and
\begin{equation*}
\fromto{\SS^X_{(\mathscr{X},\mathscr{X}_{\mathrm{rc}},\mathscr{X})}(Y)\wedge\upsidedown{\AA}_{\!\mathscr{X}}(X)}{\upsidedown{\AA}_{\!\mathscr{X}}(Y)},
\end{equation*}
We in turn obtain, for any objects $U,V\in\mathscr{X}$ such that $\fromto{U}{1}$ is relatively compact, assembly morphisms
\begin{equation*}
\fromto{\Sigma^{\infty}_{+}\Map_{\mathscr{X}}(U,X)}{F(\upsidedown{\AA}_{\!\mathscr{X}}(X),\upsidedown{\AA}_{\!\mathscr{X}}(1))}
\end{equation*}
and
\begin{equation*}
\fromto{\Sigma^{\infty}_{+}\Map_{\mathscr{X}_{\mathrm{rc}}}(V,X)}{F(\AA_{\mathscr{X}}(X),\AA_{\mathscr{X}}(1))}.
\end{equation*}
\end{nul}

\begin{exm} If $\mathscr{X}=\Kan$, then the functors $\AA$ and $\upsidedown{\AA}$ are fully functorial version of $A$-theory and $\upsidedown{A}$-theory as considered by Waldhausen in \cite[\S 2.1]{MR86m:18011} (modulo the small point that here we deal with finitely dominated spaces in place of finite spaces). The assembly morphisms above described above are
\begin{equation*}
\fromto{\Sigma^{\infty}_{+}\Map(U,X)}{\AA(X)}\textrm{\quad and\quad}\fromto{\Sigma^{\infty}_{+}\Map_{\mathrm{rc}}(V,X)}{\upsidedown{\AA}(X)},
\end{equation*}
any three objects $U,V,X\in\mathscr{X}$ such that $\fromto{U}{\ast}$ is relatively compact, where $\Map_{\mathrm{rc}}(V,X)\subset\Map(V,X)$ is the union of the connected components corresponding to maps $\fromto{V}{X}$ with finitely dominated (homotopy) fibers. When $U=\ast$, the first of these morphisms is the usual assembly morphism $\fromto{\Sigma^{\infty}_{+}X}{\AA(X)}$. Dually, we also have co-assembly morphisms
\begin{equation*}
\fromto{\Sigma^{\infty}_{+}\Map(U,X)}{F(\upsidedown{\AA}(X),\AA(\ast))}
\end{equation*}
and
\begin{equation*}
\fromto{\Sigma^{\infty}_{+}\Map_{\mathrm{rc}}(V,X)}{F(\AA(X),\AA(\ast))}.
\end{equation*}
(Observe that $\AA(\ast)\simeq\upsidedown{\AA}(\ast)$.) When $U=\ast$, the first of these morphisms seems to have been studied independently by Cary Malkiewich. When composed with the morphism induced by the trace
\begin{equation*}
\fromto{\AA(\ast)\simeq\KK(\SS)}{\mathrm{THH}(\SS)\simeq\SS},
\end{equation*}
we obtain morphisms
\begin{equation*}
\fromto{\Sigma^{\infty}_{+}\Map(U,X)}{D\upsidedown{\AA}(X)}\textrm{\quad and\quad}\fromto{\Sigma^{\infty}_{+}\Map_{\mathrm{rc}}(V,X)}{D\AA(X)}.
\end{equation*}
\end{exm}

%----------------------------------------------------------------------%

\section{Algebraic $K$-theory of derived stacks}\label{exm:stk} We turn to a more geometric setting. Over suitable schemes, complexes of quasicoherent sheaves admit both a pushforward and pullback functor. The pullback carries perfect complexes to perfect complexes, but the pushforward does not, unless some heavy constraints are placed on the morphisms of schemes involved. In this subsection, we discuss such restrictions for a very general class of derived stacks.

There are two classes of examples in which we will be interested: {\it (i)} spectral Deligne--Mumford stacks and {\it (ii)} arbitrary sheaves of spaces for the flat site. In each case we will construct algebraic $K$-theory as a spectral Mackey functor relative to certain disjunctive triple structures.

\begin{nul} For simplicity, in this section we work \emph{absolutely}, i.e., over the sphere spectrum. Nothing here uses that fact in a nontrivial way, and all the results of this section can obviously be adapted to work over more general bases.
\end{nul}

\begin{ntn} We consider the $\infty$-category $\mathbf{CAlg}^{\mathrm{cn}}$ of connective $E_{\infty}$ rings and the huge $\infty$-topos
\begin{equation*}
\mathbf{Shv}_{\mathrm{flat}}\subset\Fun(\mathbf{CAlg}^{\mathrm{cn}},\Kan(\kappa_1))
\end{equation*}
of large sheaves on $\mathbf{CAlg}^{\mathrm{cn},\op}$ for the flat topology \cite[Pr. 5.4]{DAGVII}.
\end{ntn}

Let us begin by identifying two sources of disjunctive triples within the $\infty$-category $\mathbf{Shv}_{\mathrm{flat}}$. Since our aim is to study categories of modules that are contravariantly functorial in all morphisms but only covariantly functorial in certain classes of morphisms, these disjunctive triples will have the property that every morphism will be egressive, but the ingressives will be heavily restricted.

\begin{ntn}\label{exm:relscall} The full subcategory $\mathbf{DM}\subset\mathbf{Shv}_{\mathrm{flat}}$ spanned by those sheaves that are representable by spectral Deligne--Mumford stacks is closed under finite coproducts and pullbacks, whence it is a disjunctive $\infty$-category.

For any class $\mathscr{P}$ of morphisms of $\mathbf{DM}$ that is stable under pullback and compatible with coproducts in the sense of (\ref{dfn:adequate}.\ref{item:compatiblewithcoprods}), one obtains a disjunctive triple
\begin{equation*}
(\mathbf{DM},\mathbf{DM}_{\mathscr{P}},\mathbf{DM}),
\end{equation*}
where $\mathbf{DM}_{\mathscr{P}}\subset\mathbf{DM}$ is the subcategory that contains all the objects, in which the morphisms lie in $\mathscr{P}$.
\end{ntn}

\begin{ntn} Alternately, one may opt to keep all the objects of $\mathbf{Shv}_{\mathrm{flat}}$. Then for any class $\mathscr{P}$ of morphisms of $\mathbf{Shv}_{\mathrm{flat}}$ that is stable under pullback and compatible with coproducts in the sense of (\ref{dfn:adequate}.\ref{item:compatiblewithcoprods}), one obtains a disjunctive triple
\begin{equation*}
(\mathbf{Shv}_{\mathrm{flat}},\mathbf{Shv}_{\mathrm{flat},\mathscr{P}},\mathbf{Shv}_{\mathrm{flat}}),
\end{equation*}
where $\mathbf{Shv}_{\mathrm{flat},\mathscr{P}}\subset\mathbf{Shv}_{\mathrm{flat}}$ is the subcategory that contains all the objects, in which the morphisms lie in $\mathscr{P}$.
\end{ntn}

Our ultimate interest will be in the study of perfect modules over these sorts of objects, but let us first consider the larger $\infty$-category of quasicoherent modules.

\begin{ntn} We let
\begin{equation*}
\begin{tikzpicture} 
\matrix(m)[matrix of math nodes, 
row sep=5ex, column sep=4ex, 
text height=1.5ex, text depth=0.25ex] 
{\mathbf{Mod}&\QCoh\\ 
\mathbf{CAlg}^{\mathrm{cn}}&\mathbf{Shv}_{\mathrm{flat}}^{\op}\\}; 
\path[>=stealth,->,font=\scriptsize] 
(m-1-1) edge (m-1-2) 
edge node[left]{$q$} (m-2-1) 
(m-1-2) edge node[right]{$p$} (m-2-2) 
(m-2-1) edge[right hook->] (m-2-2); 
\end{tikzpicture}
\end{equation*}
be a pullback square in which $q$ is the cocartesian fibration of \cite[Df. 4.4.1.1]{HA}, and $p$ is a cocartesian fibration classified by the right Kan extension of the functor that classifies $q$. The objects of $\QCoh$ can be thought of as pairs $(X,M)$ consisting of a sheaf $X\colon\fromto{\mathbf{CAlg}^{\mathrm{cn}}}{\Kan(\kappa_1)}$ for the flat topology and a quasicoherent module $M$ over $X$.
\end{ntn}

\begin{nul} Let us note that by \cite[Pr. 2.7.17(1)]{DAGVIII}, the fibers of $p$ are all stable $\infty$-categories. Furthermore, since the functor $\fromto{\mathbf{Shv}_{\mathrm{flat}}^{\op}}{\Cat_{\infty}(\kappa_1)}$ that classifies $p$ preserves limits \cite[Pr. 2.7.14]{DAGVIII}, it follows that for any finite set $I$ and any collection $\{X_i\ |\ i\in I\}$ of sheaves, the natural functor
\begin{equation*}
\fromto{\QCoh_{\coprod_{i\in I}X_i}}{\prod_{i\in I}\QCoh_{X_i}}
\end{equation*}
is an equivalence.
\end{nul}

\begin{wrn} Note that $p$ is not a cartesian fibration. In general, the functor
\begin{equation*}
f^{\star}\colon\fromto{\QCoh_{Y}}{\QCoh_{X}}
\end{equation*}
induced by a natural transformation $f\colon\fromto{X}{Y}$ will preserve all small colimits \cite[Pr. 2.7.17(2)]{DAGVIII}, but the $\infty$-categories may not be presentable unless one knows that $X$ and $Y$ are in some sense ``small.''
\end{wrn}

This smallness is ensured if, for example, $X$ and $Y$ are represented by spectral Deligne--Mumford stacks (\cite[Prs. 2.3.13 and 2.7.18]{DAGVIII}). We therefore conclude the following.
\begin{lem} The pulled back functor
\begin{equation*}
\fromto{\QCoh\times_{\mathbf{Shv}_{\mathrm{flat}}^{\op}}\mathbf{DM}^{\op}}{\mathbf{DM}^{\op}}
\end{equation*}
is both a cocartesian fibration and a cartesian fibration.
\end{lem}

\begin{ntn} In the $\infty$-category $\mathbf{DM}$, consider the class $\mathscr{RS}$ of \textbf{\emph{relatively scalloped}} morphisms in the sense of \cite[Df. 2.5.10]{DAGVIII}. Then since relatively scalloped morphisms are stable under pullback and compatible with coproducts in the sense of (\ref{dfn:adequate}.\ref{item:compatiblewithcoprods}), we obtain a disjunctive triple
\begin{equation*}
(\mathbf{DM},\mathbf{DM}_{\mathscr{RS}},\mathbf{DM}).
\end{equation*}
\end{ntn}

\begin{prp} The functor
\begin{equation*}
\fromto{\QCoh^{\op}\times_{\mathbf{Shv}_{\mathrm{flat}}}\mathbf{DM}}{\mathbf{DM}}
\end{equation*}
is a Waldhausen bicartesian fibration for the left complete disjunctive triple
\begin{equation*}
(\mathbf{DM},\mathbf{DM}_{\mathscr{RS}},\mathbf{DM}).
\end{equation*}
\begin{proof} By \cite[Pr. 2.5.14]{DAGVIII}, the relevant base change functors are all equivalences.
\end{proof}
\end{prp} 

\begin{ntn}\label{exm:repqaff} In the bigger $\infty$-category of all flat sheaves, we consider the class $\mathscr{QA}$ in $\mathbf{Shv}_{\mathrm{flat}}$ of \textbf{\emph{quasi-affine representable}} morphisms --- that is, those morphisms $\fromto{X}{Y}$ such that for any connective $E_{\infty}$ ring $R$ and any $R$-point $\eta\colon\fromto{\Spec R}{Y}$, the pullback $X\times_Y\Spec R$ is representable by a spectral Deligne--Mumford stack $X_{\eta}$, and the map $\fromto{X_{\eta}}{\Spec R}$ is quasi-affine \cite[Df. 3.1.24]{DAGVIII}. Then since quasi-affine representable morphisms are stable under pullback and compatible with coproducts in the sense of (\ref{dfn:adequate}.\ref{item:compatiblewithcoprods}), we have a disjunctive triple
\begin{equation*}
(\mathbf{Shv}_{\mathrm{flat}},\mathbf{Shv}_{\mathrm{flat},\mathscr{QA}},\mathbf{Shv}_{\mathrm{flat}}).
\end{equation*}
\end{ntn}

The following is immediate from \cite[Pr. 3.2.5]{DAGVIII}.
\begin{lem} The pulled back functor
\begin{equation*}
\fromto{\QCoh\times_{\mathbf{Shv}_{\mathrm{flat}}^{\op}}\mathbf{Shv}_{\mathrm{flat},\mathscr{QA}}^{\op}}{\mathbf{Shv}_{\mathrm{flat},\mathscr{QA}}^{\op}}
\end{equation*}
is both a cocartesian fibration and a cartesian fibration.
\end{lem}

Furthermore, we deduce from  \cite[Cor. 3.2.6(2)]{DAGVIII} that the relevant base change functors are all equivalences, whence we conclude the following.
\begin{prp} The functor
\begin{equation*}
p\colon\fromto{\QCoh^{\op}}{\mathbf{Shv}_{\mathrm{flat}}}
\end{equation*}
is a Waldhausen bicartesian fibration for the left complete disjunctive triple
\begin{equation*}
(\mathbf{Shv}_{\mathrm{flat}},\mathbf{Shv}_{\mathrm{flat},\mathscr{QA}},\mathbf{Shv}_{\mathrm{flat}}).
\end{equation*}
\end{prp}

\begin{ntn} Now denote by $\Perf\subset\QCoh$ the full subcategory spanned by those pairs $(X,M)$ in which $M$ is a perfect quasicoherent module on $X$ --- i.e., a strongly dualizable object of the symmetric monoidal $\infty$-category $\QCoh_X$. Endow $\Perf^{\op}$ with its fiberwise maximal pair structure, so that
\begin{equation*}
\Perf^{\dag}\coloneq\Perf\times_{\mathbf{Shv}_{\mathrm{flat}}^{\op}}\iota\mathbf{Shv}_{\mathrm{flat}}^{\op}.
\end{equation*}
Note that the restricted functor
\begin{equation*}
\fromto{\Perf^{\op}}{\mathbf{Shv}_{\mathrm{flat}}}
\end{equation*}
remains a Waldhausen cartesian fibration, since pullbacks are symmetric monoidal.
\end{ntn}

We now wish to shrink our classes of ingressives on $\mathbf{DM}$ and $\mathbf{Shv}_{\mathrm{flat}}$ to ensure that perfect modules define a supported Waldhausen bicartesian fibration. For this, it's enough to find a class of morphisms such that the pushforward preserves perfection.
\begin{dfn} Suppose $f\colon\fromto{X}{Y}$ a morphism of $\mathbf{Shv}_{\mathrm{flat}}$ such that the functor $f^{\star}\colon\fromto{\QCoh_Y}{\QCoh_X}$ admits a right adjoint $f_{\star}$. Call $f$ \textbf{\emph{perfect}} if $f_{\star}\colon\fromto{\QCoh_X}{\QCoh_Y}$ carries perfect objects \cite[Df. 2.7.21]{DAGVIII} to perfect objects.
\end{dfn}

Between spectral Deligne--Mumford stacks, Lurie identifies an important class of perfect morphisms.
\begin{prp}[Lurie, \protect{\cite[Pr. 3.3.20]{DAGXII}}] If a morphism $f\colon\fromto{X}{Y}$ of spectral Deligne--Mumford stacks is strongly proper \cite[Df. 3.1.1]{DAGXII}, of finite Tor-amplitude \cite[Df. 3.3.9]{DAGXII}, and locally almost of finite presentation \cite[Df. 8.16]{DAGIX}, then it is perfect.
\end{prp}

\begin{ntn}\label{ntn:properness} Denote by $\mathscr{F\!P}$ the class of strongly proper morphisms of finite Tor-amplitude \cite[Df. 3.3.9]{DAGXII} and locally almost of finite presentation. To see that the class $\mathscr{F\!P}$ is stable under pullbacks, combine \cite[Rk. 3.1.5 and Pr. 3.3.16]{DAGXII} and \cite[Pr. 8.25]{DAGIX}. It is easy to see that this class is compatible with coproducts in the sense of (\ref{dfn:adequate}.\ref{item:compatiblewithcoprods}).
\end{ntn}

\begin{prp} The functor
\begin{equation*}
p_{\mathbf{DM}}\colon\fromto{\Perf^{\op}\times_{\mathbf{Shv}_{\mathrm{flat}}}\mathbf{DM}}{\mathbf{DM}}
\end{equation*}
is a Waldhausen bicartesian fibration for the left complete disjunctive triple
\begin{equation*}
(\mathbf{DM},\mathbf{DM}_{\mathscr{F\!P}},\mathbf{DM}).
\end{equation*}
\end{prp}

\begin{cor} The unfurling
\begin{equation*}
\fromto{\Upsilon(p_{\mathbf{DM}})}{A^{\eff}(\mathbf{DM},\mathbf{DM}_{\mathscr{F\!P}},\mathbf{DM})}
\end{equation*}
is classified by a Mackey functor
\begin{equation*}
\mathscr{M}_{\mathbf{DM}}\colon\fromto{A^{\eff}(\mathbf{DM},\mathbf{DM}_{\mathscr{F\!P}},\mathbf{DM})}{\Wald_{\infty}}.
\end{equation*}
\end{cor}

\begin{ntn}\label{ntn:KofDMstacks} As a result, we obtain a Mackey functor
\begin{equation*}
\KK\circ\mathscr{M}_{\mathbf{DM}}\colon\fromto{A^{\eff}(\mathbf{DM},\mathbf{DM}_{\mathscr{F\!P}},\mathbf{DM})}{\Sp}.
\end{equation*}
The algebraic $K$-theory of any Deligne--Mumford stack is an $E_{\infty}$ ring spectrum, whence it has a canonical $K$-theory class, given by the unit. Correspondingly, we obtain a morphism of Mackey functors
\begin{equation*}
\fromto{\SS_{(\mathbf{DM},\mathbf{DM}_{\mathscr{F\!P}},\mathbf{DM})}}{\KK\circ\mathscr{M}_{\mathbf{DM}}}.
\end{equation*}
We thus obtain, for any Deligne--Mumford stacks $U$ and $X$, an assembly morphism
\begin{equation*}
\fromto{\Sigma^{\infty}_{+}\Mor_{\mathscr{F\!P}}(U,X)}{\KK(X)},
\end{equation*}
where $\Mor_{\mathscr{F\!P}}$ is the space of strongly proper morphisms of finite Tor-amplitude and locally almost of finite presentation.

As in Df. \ref{dfn:assembly}, we obtain a morphism of spectral Mackey functors
\begin{equation*}
\fromto{\SS^X_{(\mathbf{DM},\mathbf{DM}_{\mathscr{F\!P}},\mathbf{DM})}}{F(\KK(X),\KK\circ\mathscr{M}_{\mathbf{DM}})},
\end{equation*}
which induces assembly morphisms 
\begin{equation*}
\fromto{\SS^X_{(\mathbf{DM},\mathbf{DM}_{\mathscr{F\!P}},\mathbf{DM})}(Y)\wedge\KK(X)}{\KK(Y)}.
\end{equation*}
We in turn obtain, for any spectral Deligne--Mumford stack $U$ such that the morphism $\fromto{U}{\Spec^{\et}(\SS)}$ is strongly proper, of finite Tor-amplitude, and locally almost of finite presentation, an assembly morphism
\begin{equation*}
\fromto{\Sigma^{\infty}_{+}\Map(U,X)}{F(\KK(X),\KK(\SS))}.
\end{equation*}
When composed with the morphism induced by the trace
\begin{equation*}
\fromto{\KK(\SS)}{\mathrm{THH}(\SS)\simeq\SS},
\end{equation*}
we obtain a morphism
\begin{equation*}
\fromto{\Sigma^{\infty}_{+}\Map(U,X)}{D\KK(X)}.
\end{equation*}
\end{ntn}

Now among the more general class of flat sheaves, let us characterize perfect morphisms that are also quasi-affine representable.
\begin{prp}\label{prp:perf} The following are equivalent for a quasi-affine representable morphism $f\colon\fromto{X}{Y}$ of $\mathbf{Shv}_{\mathrm{flat}}$.
\begin{enumerate}[(\ref{prp:perf}.1)]
\item The morphism $f$ is perfect.
\item The quasicoherent module $f_{\star}\mathscr{O}_X$ is perfect.
\end{enumerate}
\begin{proof} That the first condition implies the second is obvious. To prove the converse, let us first reduce to the case where $Y$ is affine. Indeed, for any quasicoherent module $M$ over $X$, the quasicoherent module $f_{\star}M$ is perfect just in case, for any point $\eta\colon\fromto{\Spec R}{Y}$, the $R$-module $\eta^{\star}f_{\star}M$ is perfect. Now consider the pullback
\begin{equation*}
\begin{tikzpicture} 
\matrix(m)[matrix of math nodes, 
row sep=4ex, column sep=4ex, 
text height=1.5ex, text depth=0.25ex] 
{X_{\eta}&\Spec R\\ 
X&Y.\\}; 
\path[>=stealth,->,font=\scriptsize] 
(m-1-1) edge node[above]{$f'$} (m-1-2) 
edge node[left]{$\varepsilon$} (m-2-1) 
(m-1-2) edge node[right]{$\eta$} (m-2-2) 
(m-2-1) edge node[below]{$f$} (m-2-2); 
\end{tikzpicture}
\end{equation*}
One now has $\eta^{\star}f_{\star}M\simeq f'_{\star}\varepsilon^{\star}M$ \cite[Cor. 3.2.6(2)]{DAGVIII}, and if $M$ is perfect, so is $\varepsilon^{\star}M$.

So we now assume that $Y=\Spec R$. Now let $A\coloneq f_{\star}\mathscr{O}_X$, an $R$-module. In light of \cite[Pr. 3.2.5]{DAGVIII}, we obtain an equivalence $\QCoh_X\simeq\Mod_A$, under which the functor $f_{\star}$ can be identified with the forgetful functor
\begin{equation*}
U\colon\fromto{\Mod_A}{\Mod_R}.
\end{equation*}
Now $U^{-1}(\Perf_R)$ is a stable subcategory that is closed under retracts, and by assumption it contains $A$. Hence $\Perf_A\subset U^{-1}(\Perf_R)$.
\end{proof}
\end{prp}

\begin{cor} The collection $\mathscr{QP}$ of quasi-affine representable and perfect morphisms of $\mathbf{Shv}_{\mathrm{flat}}$ is stable under pullback.
\begin{proof} Suppose
\begin{equation*}
\begin{tikzpicture} 
\matrix(m)[matrix of math nodes, 
row sep=4ex, column sep=4ex, 
text height=1.5ex, text depth=0.25ex] 
{X'&Y'\\ 
X&Y\\}; 
\path[>=stealth,->,font=\scriptsize] 
(m-1-1) edge node[above]{$f'$} (m-1-2) 
edge node[left]{$\alpha$} (m-2-1) 
(m-1-2) edge node[right]{$\beta$} (m-2-2) 
(m-2-1) edge node[below]{$f$} (m-2-2); 
\end{tikzpicture}
\end{equation*}
a pullback square in $\mathbf{Shv}_{\mathrm{flat}}$, and suppose $f$ quasi-affine representable and perfect. Then the quasicoherent module
\begin{equation*}
f'_{\star}\mathscr{O}_{X'}\simeq f'_{\star}\alpha^{\star}\mathscr{O}_X\simeq\beta^{\star}f_{\star}\mathscr{O}_X
\end{equation*}
is perfect as well.
\end{proof}
\end{cor}

It is a straightforward matter to see that $\mathscr{QP}$ is compatible with coproducts in the sense of (\ref{dfn:adequate}.\ref{item:compatiblewithcoprods}). We thus conclude the following.

\begin{prp} The functor
\begin{equation*}
\fromto{\Perf^{\op}}{\mathbf{Shv}_{\mathrm{flat}}}
\end{equation*}
is a Waldhausen bicartesian fibration for the left complete disjunctive triple
\begin{equation*}
(\mathbf{Shv}_{\mathrm{flat}},\mathbf{Shv}_{\mathrm{flat},\mathscr{QP}},\mathbf{Shv}_{\mathrm{flat}}).
\end{equation*}
\end{prp}

\begin{cor} The unfurling
\begin{equation*}
\fromto{\Upsilon(p_{\mathbf{Shv}_{\mathrm{flat}}})}{A^{\eff}(\mathbf{Shv}_{\mathrm{flat}},\mathbf{Shv}_{\mathrm{flat},\mathscr{QP}},\mathbf{Shv}_{\mathrm{flat}})}
\end{equation*}
is classified by a Mackey functor
\begin{equation*}
\mathscr{M}_{\mathbf{Shv}_{\mathrm{flat}}}\colon\fromto{A^{\eff}(\mathbf{Shv}_{\mathrm{flat}},\mathbf{Shv}_{\mathrm{flat},\mathscr{QP}},\mathbf{Shv}_{\mathrm{flat}})}{\Wald_{\infty}}.
\end{equation*}
\end{cor}

\begin{ntn}\label{ntn:KasMackonflatsheaves} As a result, we obtain a Mackey functor
\begin{equation*}
\KK\circ\mathscr{M}_{\mathbf{Shv}_{\mathrm{flat}}}\colon\fromto{A^{\eff}(\mathbf{Shv}_{\mathrm{flat}},\mathbf{Shv}_{\mathrm{flat},\mathscr{QP}},\mathbf{Shv}_{\mathrm{flat}})}{\Sp}.
\end{equation*}
The algebraic $K$-theory of any flat sheaf is an $E_{\infty}$ ring spectrum, whence it has a canonical $K$-theory class, given by the unit. Correspondingly, we obtain a morphism of Mackey functors
\begin{equation*}
\fromto{\SS_{(\mathbf{Shv}_{\mathrm{flat}},\mathbf{Shv}_{\mathrm{flat},\mathscr{QP}},\mathbf{Shv}_{\mathrm{flat}})}}{\KK\circ\mathscr{M}_{\mathbf{Shv}_{\mathrm{flat}}}}.
\end{equation*}
We thus obtain, for any flat sheaves $U$ and $X$, an assembly morphism
\begin{equation*}
\fromto{\Sigma^{\infty}_{+}\Mor_{\mathscr{QP}}(U,X)}{\KK(X)},
\end{equation*}
where $\Mor_{\mathscr{QP}}$ is the space of quasi-affine representable and perfect morphisms.

Dually, as in Df. \ref{dfn:assembly}, we obtain a morphism of Mackey functors
\begin{equation*}
\fromto{\SS^X_{(\mathbf{Shv}_{\mathrm{flat}},\mathbf{Shv}_{\mathrm{flat},\mathscr{QP}},\mathbf{Shv}_{\mathrm{flat}})}}{F(\KK(X),\KK\circ\mathscr{M}_{\mathbf{Shv}_{\mathrm{flat}}})},
\end{equation*}
which induces assembly morphisms
\begin{equation*}
\fromto{\SS^X_{(\mathbf{Shv}_{\mathrm{flat}},\mathbf{Shv}_{\mathrm{flat},\mathscr{QP}},\mathbf{Shv}_{\mathrm{flat}})}(Y)\wedge\KK(X)}{\KK(Y)}
\end{equation*}
for any flat sheaf $Y$. We in turn obtain, for any flat sheaf $U$ such that $\fromto{U}{\Spec^{\mathrm{f}}(\SS)}$ is quasiaffine representable and perfect, an assembly morphism
\begin{equation*}
\fromto{\Sigma^{\infty}_{+}\Map(U,X)}{F(\KK(X),\KK(\SS))}.
\end{equation*}
When composed with the morphism induced by the trace
\begin{equation*}
\fromto{\AA(\ast)\simeq\KK(\SS)}{\mathrm{THH}(\SS)\simeq\SS},
\end{equation*}
we obtain a morphism
\begin{equation*}
\fromto{\Sigma^{\infty}_{+}\Map(U,X)}{D\KK(X)}.
\end{equation*}
\end{ntn}

Finally, let us restrict Ex. \ref{exm:stk} to relate it to Ex. \ref{exm:profgroups}.

\begin{ntn} Suppose $X$ a spectral Deligne--Mumford stack. We denote by $\mathbf{F\acute{E}t}(X)$ the subcategory of $\mathbf{DM}_{/X}$ whose objects are finite \cite[Df. 3.2.4]{DAGXII} and \'etale morphisms $\fromto{Y}{X}$ and whose morphisms are finite and \'etale morphisms over $X$. We will abuse notation and write $A^{\eff}(X)$ for the effective Burnside $\infty$-category of $\mathbf{F\acute{E}t}(X)$.
\end{ntn}

\begin{exm}\label{exm:GaloisequivariantK} Suppose $X$ a connected, noetherian scheme, suppose $x$ a geometric point of $X$. Then if $\pi_1^{\mathrm{\acute{e}t}}(X,x)$ is the \'etale fundamental group of $X$, then by Grothendieck's Galois duality \cite[Exp. V, \S 7]{MR50:7129}, the $\infty$-category $\mathbf{F\acute{E}t}(X)$ is canonically equivalent to the $\infty$-category $NB_{\pi_1^{\mathrm{\acute{e}t}}(X,x)}^{\mathit{fin}}$ of Ex. \ref{exm:profgroups}, whence we have a canonical equivalence
\begin{equation*}
A^{\eff}(\pi_1^{\mathrm{\acute{e}t}}(X,x))\simeq A^{\eff}(X).
\end{equation*}
(We can, of course, relax the condition that $X$ be connected by passing to the \'etale fundamental groupoid.)

Restricting the Mackey functor $\KK\circ\mathscr{M}_{\mathbf{DM}}$ of Nt. \ref{ntn:KofDMstacks} to
\begin{equation*}
A^{\eff}(X)\subset A^{\eff}(\mathbf{DM},\mathbf{DM}_{\mathscr{F\!P}},\mathbf{DM}),
\end{equation*}
we obtain a $\pi_1^{\et}(X,x)$-equivariant $K$-theory spectrum
\begin{equation*}
\KK_{\pi_1^{\et}(X,x)}(X)\colon\fromto{A^{\eff}(\pi_1^{\et}(X,x))}{\Sp},
\end{equation*}
whose value on the $\pi_1^{\et}(X,x)$-set $[\pi_1^{\et}(X,x)/H]$ (with $H\leq\pi_1^{\et}(X,x)$ open) is the algebraic $K$-theory of the \'etale cover $\fromto{X'}{X}$ corresponding to $H$. Corresponding to the unit of $\KK(X)$, we obtain a morphism of Mackey functors
\begin{equation*}
\fromto{\SS_{\pi_1^{\et}(X,x)}}{\KK_{\pi_1^{\et}(X,x)}(X)}.
\end{equation*}
The Segal--tom Dieck splitting thus provides a collection of maps
\begin{equation*}
\fromto{\Sigma_{+}^{\infty}B(N_{\pi_1^{\et}(X,x)}H/H)}{\KK(X)},
\end{equation*}
one for each conjugacy class of open subgroups $H\leq\pi_1^{\et}(X,x)$.
\end{exm}

%----------------------------------------------------------------------%
%----------------------------------------------------------------------%
%----------------------------------------------------------------------%

\bibliographystyle{amsplain}
\bibliography{kthy}

\end{document}